\documentclass[reqno,a4paper,11pt]{amsart}
\usepackage[all,poly]{xy}
\usepackage{amsfonts}
\usepackage{amssymb}
\usepackage{amsmath}
\usepackage{color}
\usepackage[pagebackref,colorlinks]{hyperref}
\usepackage{enumerate}
\usepackage{comment}

\DeclareFontFamily{U}{mathx}{\hyphenchar\font45}
\DeclareFontShape{U}{mathx}{m}{n}{
      <5> <6> <7> <8> <9> <10>
      <10.95> <12> <14.4> <17.28> <20.74> <24.88>
      mathx10
      }{}
\DeclareSymbolFont{mathx}{U}{mathx}{m}{n}
\DeclareFontSubstitution{U}{mathx}{m}{n}
\DeclareMathSymbol{\bigplus}        {1}{mathx}{"90}
\DeclareMathSymbol{\bigtimes}       {1}{mathx}{"91}
     
\theoremstyle{plain}
\newtheorem {lemma}{Lemma}[subsection]
\newtheorem {proposition}[lemma]{Proposition}
\newtheorem {theorem}[lemma]{Theorem}
\newtheorem {corollary}[lemma]{Corollary}

\theoremstyle{definition}
\newtheorem {definition}[lemma]{Definition}
\newtheorem {remark}[lemma]{Remark}
\newtheorem {example}[lemma]{Example}

\newenvironment{proofs}{%
  \proof}{\endproof}

\newcommand{\R}{\mathbb{R}}
\newcommand{\N}{\mathbb{N}}
\newcommand{\Z}{\mathbb{Z}}
\newcommand{\X}{\langle X \rangle}
\newcommand{\GKdim}{\operatorname{GKdim}}
\newcommand{\nod}{\operatorname{nod}}
\newcommand{\reg}{\operatorname{reg}}
\newcommand{\st}{\operatorname{st}}
\newcommand{\fin}{\operatorname{fin}}
\newcommand{\Mat}{\operatorname{\mathbb{M}}}
\newcommand{\im}{\operatorname{im}}
\newcommand{\id}{\operatorname{id}}
\newcommand{\Id}{\operatorname{Id}}
\newcommand{\coker}{\operatorname{coker}}
\newcommand{\diag}{\operatorname{diag}}
\newcommand{\ann}{\operatorname{ann}}
\newcommand{\Span}{\operatorname{span}}
\newcommand{\NF}{\operatorname{NF}}
\newcommand{\supp}{\operatorname{supp}}
\newcommand{\G}{\operatorname{{\bf WG}}}
\newcommand{\M}{\operatorname{{\bf MON}}}
\newcommand{\A}{\operatorname{{\bf ALG}}}
\newcommand{\Gr}{\operatorname{{\bf GR}}}
\newcommand{\Mod}{\operatorname{{\bf MOD}}}
\newcommand{\proj}{\operatorname{proj}}
\newcommand{\gr}{\operatorname{gr}}
\newcommand{\cX}{\mathcal{X}}
\newcommand{\cY}{\mathcal{Y}}
\newcommand{\Path}{\operatorname{Path}}
\newcommand{\Ob}{\operatorname{Ob}}

\DeclareMathOperator{\RG}{\mathbf{RG}}

\newcommand{\nb}{\operatorname{red}}
\newcommand{\End}{\operatorname{End}}
\newcommand{\tg}{\operatorname{tag}}
\newcommand{\mot}{\operatorname{st}}
\newcommand{\op}{\operatorname{op}}
\newcommand{\Hom}{\operatorname{Hom}}
\newcommand{\Mf}{\mathfrak{M}}
\newcommand{\V}{\mathcal{V}}
\newcommand{\cM}{\mathcal{M}}

\title[Weighted Leavitt path algebras]{Weighted Leavitt path algebras - an overview}

\author{Raimund Preusser}
\email{raimund.preusser@gmx.de}
\date{}

\makeatletter
\@namedef{subjclassname@2020}{%
  \textup{2020} Mathematics Subject Classification}
\makeatother


\begin{document}

\begin{abstract} 
Weighted Leavitt path algebras were introduced in 2013 by Roozbeh Hazrat. These algebras generalise simultaneously the usual Leavitt path algebras and William Leavitt's algebras $L(m,n)$. In this paper we try to give an overview of what is known about the weighted Leavitt path algebras. We also prove some new results (in particular on the graded K-theory of weighted Leavitt path algebras) and mention open problems.
\end{abstract}

\maketitle

\setcounter{tocdepth}{1}
\tableofcontents

\section{Introduction}
Leavitt path algebras are algebras associated to directed graphs. They were introduced by G. Abrams and G. Aranda Pino in 2005 \cite{aap05} and independently by P. Ara, M. Moreno and E. Pardo in 2007 \cite{Ara_Moreno_Pardo}. The Leavitt path algebras turned out to be a very rich and interesting class of algebras, whose studies so far have comprised over 150 research papers and counting. A comprehensive treatment of the subject can be found in the book~\cite{abrams-ara-molina}. %These algebras have created a stir in the algebraic community, on arXiv one can already find 170 papers containing ``Leavitt path algebra'' in the title. %Even the fields medal winner E. Zelmanov contributed two papers to the topic \cite{zelmanov1, zelmanov2}.\\
The definition of the Leavitt path algebras was inspired by the algebras $L(m,n)$ studied by W. Leavitt in the 1950's and 60's \cite{vitt56, vitt57, vitt62, vitt65}. Recall that for positive integers $m<n$, the Leavitt algebra $L(m,n)$ is universal with the property that $L(m,n)^m\cong L(m,n)^{n}$ as left $L(m,n)$-modules. 
% Recall that for a field $K$ and positive integers $m<n$, the Leavitt algebra $L_K(m,n)$ is the $K$-algebra $A$ universal with the property that $A^m\cong A^{n}$ as left $A$-modules. 
%with a universal left $A$-module isomorphism $A^m\rightarrow A^{n}$. %Recall that for positive integers $m<n$, the Leavitt algebra $L(m,n)$ is universal with the property that $L(m,n)^m\cong L(m,n)^{n}$ as left $L(m,n)$-modules. 

Weighted Leavitt path algebras are algebras associated to {\it weighted graphs}, i.e. directed graphs with a weight map associating to each edge a positive integer. These algebras were introduced by R. Hazrat in 2013 \cite{hazrat13}. If all the weights are equal to $1$, then the weighted Leavitt path algebras reduce to the usual Leavitt path algebras. While their unweighted cousins only include the Leavitt algebras $L(1,n)~(n>1)$ as special cases, the weighted Leavitt path algebras embrace {\it all} Leavitt algebras (see \cite{Raimund4}). 

In \cite{hazrat-preusser} Gr\"obner bases for weighted Leavitt path algebras were obtained. Using these bases it was shown that a simple or graded simple weighted Leavitt path algebra is isomorphic to an unweighted Leavitt path algebra. Moreover, the weighted Leavitt path algebras that are domains were classified, and it was shown that there is a huge class of weighted Leavitt path algebras having a ``local'' valuation.

Recall that for a unital ring $R$, the {\it $\V$-monoid} $\V(R)$ of $R$ is the set of isomorphism classes of finitely generated projective left $R$-modules, which becomes an abelian monoid by defining $[P]+[Q]:=[P\oplus Q]$ for any $[P],[Q]\in \V(R)$. This definition can be extended to nonunital rings. For a ring $R$ with local units the Grothendieck group $K_0(R)$ is the group completion of $\V(R)$. In \cite{preusser-1} the $\V$-monoid and the Grothendieck group of a weighted Leavitt path algebra were computed.

In \cite{preusser} the Gelfand-Kirillov dimension of a weighted Leavitt path algebra was determined. Moreover, it was shown that a finite-dimensional weighted Leavitt path algebra is isomorphic to an unweighted Leavitt path algebra.

Recall that a group graded $K$-algebra $A=\bigoplus_{g\in G} A_g$ is called {\it locally finite} if $\dim_K A_g < \infty$ for every $g\in G$. In \cite{preusser1} the weighted graphs $(E,w)$ for which the weighted Leavitt path algebra $L(E,w)$ is locally finite with respect to its standard grading were characterised. It was also shown that the locally finite weighted Leavitt path algebras are precisely the Noetherian ones, and that a locally finite weighted Leavitt path algebra is isomorphic to an unweighted Leavitt path algebra.

In \cite{Raimund2} Leavitt path algebras of hypergraphs were introduced, which include as special cases the {\it vertex-weighted Leavitt path algebras}, i.e. weighted Leavitt path algebras of weighted graphs having the property that any two edges emitted by the same vertex have the same weight. Moreover, the graded $\V$-monoid (see Section 10) of a Leavitt path algebra of a hypergraph was computed, covering the case of a vertex-weighted Leavitt path algebra.

In \cite{Raimund4} the weighted graphs $(E,w)$ having the property that the weighted Leavitt path algebra $L(E,w)$ is isomorphic to an unweighted Leavitt path algebra were characterised. Moreover, it was shown  that if a weighted Leavitt path algebra is Artinian, or von Neumann regular, or has finite Gelfand-Kirillov dimension, then it is isomorphic
to an unweighted Leavitt path algebra. 

In \cite{hazrat-preusser-shchegolev} modules for weighted Leavitt path algebras were obtained by introducing the notion of a representation graph for a weighted graph. It was shown that each connected component of the category $\RG(E,w)$ of representation graphs for a weighted graph $(E,w)$ contains a unique object yielding a simple $L(E,w)$-module, and a ``universal'' object, yielding an indecomposable $L(E,w)$-module. It was also shown that specialising to unweighted graphs, one recovers the simple modules for the usual Leavitt path algebras constructed by X. Chen via infinite paths \cite{C}.

In the remainder of this paper we depict the above-mentioned results on weighted Leavitt path algebras in more detail. We sometimes state results without giving a proof, or we only give a {\it Sketch of Proof}. In either case we provide a reference to a formal proof. The paper also contains some new results, for which we of course provide formal proofs. In Section 2 we recall some graph-theoretical notions and moreover the definition of the Leavitt algebras $L(m,n)$. In Section 3 we recall the definitions of the unweighted and weighted Leavitt path algebras. In Section 4 we present the Gr\"obner basis result from \cite{hazrat-preusser} and the computation of the Gelfand-Kirillov dimension from \cite{preusser}. Section 5 contains the characterisation of weighted Leavitt path algebras that are isomorphic to unweighted Leavitt path algebras from \cite{Raimund4}. In Section 6 we give graph-theoretical criteria for finite-dimensionality, Noetherianess and von Neumann regularity of weighted Leavitt path algebras. The criteria for finite-dimensionality and Noetherianess were obtained in \cite{preusser1}, the criterion for von Neumann regularity is new. In Section 7 we introduce the notion of a generalised corner skew Laurent polynomial ring, generalising the notion of a corner skew Laurent polynomial ring (see e.g. \cite{arabrucom}). The definition of a generalised corner skew Laurent polynomial ring already appeared in an early version of \cite{hazrat-preusser}, but is not contained in the final version of that article. We prove that a weighted Leavitt path algebra of a finite weighted graph without sinks is graded isomorphic to a generalised corner skew Laurent polynomial ring. In Section 8 we describe the local valuations for weighted Leavitt path algebras found in \cite{hazrat-preusser}. Section 9 contains the computation of the $V$-monoid and the Grothendieck group of a weighted Leavitt path algebra from \cite{preusser-1}. In Section 10, we compute the graded $V$-monoid and the graded Grothendieck group of a weighted Leavitt path algebra using the Leavitt path algebras of bi-separated graphs, which were recently introduced by R. Mohan and B. Suhas \cite{mohan-suhas}. By doing so, we generalise the graded $\V$-monoid result for vertex-weighted Leavitt path algebras obtained in \cite{Raimund2} to arbitrary weighted Leavitt path algebras. In Section 11 we present the modules for weighted Leavitt path algebras found in \cite{hazrat-preusser-shchegolev}. In Section 12 we list some open problems.

\section{Preliminaries}

Throughout the paper $K$ denotes a fixed field. Rings and algebras are associative but not necessarily commutative or unital. By an ideal we mean a two-sided ideal. $\N$ denotes the set of positive integers, $\N_0$ the set of nonnegative integers, $\Z$ the set of integers and $\R_+$ the set of positive real numbers.

\subsection{Graphs}

A {\it (directed) graph} is a quadruple $E=(E^0,E^1,s,r)$ where $E^0$ and $E^1$ are sets and $s,r:E^1\rightarrow E^0$ maps. The elements of $E^0$ are called {\it vertices} and the elements of $E^1$ {\it edges}. If $e$ is an edge, then $s(e)$ is called its {\it source} and $r(e)$ its {\it range}. If $v$ is a vertex and $e$ an edge, we say that $v$ {\it emits} $e$ if $s(e)=v$, and $v$ {\it receives} $e$ if $r(e)=v$. A vertex is called a {\it source} if it receives no edges,
a {\it sink} if it emits no edges, an {\it infinite emitter} if it emits infinitely many edges and {\it regular} if it is neither a sink nor an infinite emitter. The subset of $E^0$ consisting of all regular vertices is denoted by $E^0_{\reg}$. The graph $E$ is called {\it row-finite} if no vertex is an infinite emitter, {\it finite} if $E^0$ and $E^1$ are finite sets, and {\it empty} if $E^0=E^1=\emptyset$.

Let $E$ and $F$ be graphs. A {\it graph homomorphism} $\phi: E\to F$ consists of two maps $\phi^0 : E^0\to F^0$ and $\phi^1 : E^1\to F^1$ such that $s(\phi^1(e)) = \phi^0(s(e))$ and $r(\phi^1(e)) = \phi^0(r(e))$ for any $e\in E^1$. A graph $G$ is called a {\it subgraph} of a graph $E$ if $G^0\subseteq E^0$, $G^1\subseteq E^1$, $s_G=s_E|_{G^0}$ and $r_G=r_E|_{G^0}$.

Let $E$ be a graph. The graph $E_d=(E_d^0, E_d^1, s_d, r_d)$ where $E_d^0=E^0$, $E_d^1=\{e,e^*\mid e\in E^1\}$, and $s_d(e)=s(e),~r_d(e)=r(e),~s_d(e^*)=r(e),~r_d(e^*)=s(e)$ for any $e\in E^1$ is called the {\it double graph} of $E$. Sometimes the edges in the graph $E$ are called {\it real edges} and the additional edges in $E_d$ {\it ghost edges}. 

A {\it path} in a graph $E$ is a finite, nonempty word $p=x_1\dots x_n$ over the alphabet $E^0\cup E^1$ such that either $x_i\in E^1~(i=1,\dots,n)$ and $r(x_i)=s(x_{i+1})~(i=1,\dots,n-1)$ or $n=1$ and $x_1\in E^0$. By definition, the {\it length} $|p|$ of $p$ is $n$ in the first case and $0$ in the latter case. $p$ is called {\it nontrivial} if $|p|>0$ and {\it trivial} otherwise. We set $s(p):=s(x_1)$ and $r(p):=r(x_n)$ using the convention $s(v)=r(v)=v$ for any $v\in E^0$. A {\it closed path (based at $v$)} is a nontrivial path $p$ such that $s(p)=r(p)=v$. A {\it cycle (based at $v$)} is a closed path $p=x_1\dots x_n$ based at $v$ such that $s(x_i)\neq s(x_j)$ for any $i\neq j$. An edge $e\in E^1$ is called an {\it exit} of a cycle $x_1\dots x_n$ if there is an $i\in \{1,\dots,n\}$ such that $s(e)=s(x_i)$ and $e\neq x_i$.

If $u,v\in E^0$ and there is a path $p$ such that $s(p)=u$ and $r(p)=v$, then we write $u\geq v$.  %Clearly $\geq$ is a preorder on $E^0$.
If $u\in E^0$, then $T(u):=\{v\in E^0 \mid u\geq v\}$ is called the {\it tree} of $u$. If $X\subseteq E^0$, we define $T(X):=\bigcup_{v\in X}T(v)$. A subset $H\subseteq E^0$ is called {\it hereditary} if $T(H)\subseteq H$. %Clearly any tree of a vertex is a hereditary subset of $E^0$.
Two edges $e,f\in E^1$ are called {\it in line} if $e=f$ or $r(e)\geq s(f)$ or $r(f)\geq s(e)$. The graph $E$ is called {\it connected} if $u\geq v$ in $E_d$ for any $u,v\in E_d^0=E^0$.

Let $E$ be a graph. The $K$-algebra $P(E)$ presented by the generating set $E^0\cup E^1$ and the relations 
\begin{enumerate}[(i)]
\item $uv=\delta_{uv}u$ for any $u,v\in E^0$ and
\item $s(e)e=e=er(e)$ for any $e\in E^1$
\end{enumerate}
is called the {\it path algebra} of $E$. The paths in $E$ form a linear basis for $P(E)$.

In order to simplify the exposition we make the following assumption.\\
\\
{\bf STANDING ASSUMPTION:} In this paper all graphs are assumed to be nonempty, row-finite and connected.

\subsection{Weighted graphs}
A {\it weighted graph} is a pair $(E,w)$ where $E$ is a graph and $w:E^1\rightarrow \N$ is a map. If $e\in E^1$, then $w(e)$ is called the {\it weight} of $e$. An edge $e\in E^{1}$ is called {\it unweighted} if $w(e)=1$ and {\it weighted} otherwise. The subset of $E^1$ consisting of all unweighted edges is denoted by $E_{uw}^1$ and the subset consisting of all weighted edges by $E_{w}^1$. We set $w(v):=\max\{w(e)\mid e\in s^{-1}(v)\}$ for any regular vertex $v$ and $w(v):=0$ for any sink $v$. % The weighted subgraph $(E_{T(E^{0}_w)},w_{T(E^{0}_w)})$ of $(E,w)$ defined by the hereditary subset $T(E^{0}_w)\subseteq E^0$ is called the {\it weighted part of $(E,w)$}.
A weighted graph $(E,w)$ is called {\it finite} if the graph $E$ is finite.

A homomorphism of weighted graphs is a graph homomorphism that preserves the weights. A weighted graph $(G,w_G)$ is called a {\it weighted subgraph} of a weighted graph $(E,w_E)$ if $G$ is a subgraph of $E$ and $w_G(g)=w_E(g)$ for any $g\in G^1$.

Let $(E,w)$ be a weighted graph. The graph $\hat E=(\hat E^0, \hat E^1, \hat s, \hat r)$ where $\hat E^0=E^0$, $\hat E^1:=\{e_1,\dots,e_{w(e)}\mid e\in E^1\}$, $\hat s(e_i)=s(e)$ and $\hat r(e_i)=r(e)$ for any $e\in E^1$ and $1\leq i\leq w(e)$ is called the {\it unweighted graph associated to $(E,w)$}.\\
\\
{\bf WARNING:} Let $(E,w)$ be a weighted graph. In \cite{hazrat13} and \cite{hazrat-preusser}, the set $E^1$ was denoted by $E^{\st}$. The set $\hat E^1=\{e_1,\dots,e_{w(e)}\mid e\in E^1\}$ was denoted by $E^1$.

\subsection{Leavitt algebras}

Fix positive integers $m<n$. If $F$ is a field, then $F^m\not\cong F^n$ as $F$-vector spaces (by the dimension theorem for vector spaces). One can ask if there is a unital ring $R$ such that $R^m\cong R^n$ as left $R$-modules. Suppose we have found a unital ring $R$ such that there are matrices $X\in \Mat_{m\times n}(R)$ and $Y\in\Mat_{n \times m}(R)$ such that 
\begin{equation}
XY=\Id_m\text{ and }YX=\Id_n.
\end{equation} 
Then $X$ and $Y$ will define left $R$-module homomorphisms $R^m\to R^n$ respectively $R^n\to R^m$ which are inverse to each other, whence $R^m\cong R^n$. It is easy to construct such a ring $R$: Let $R$ be the unital $K$-algebra presented by the generating set $\{x_{ij},y_{ji}\mid 1\leq i\leq m, 1\leq j\leq n\}$ and the relations 
\begin{equation}
\sum_{j=1}^nx_{kj}y_{jl}=\delta_{kl}~(1\leq k,l\leq m)~\text{  and  }~\sum_{i=1}^my_{ki}x_{il}=\delta_{kl}~(1\leq k,l\leq n).\end{equation}
Let $X\in \Mat_{m\times n}(R)$ be the matrix whose entry at position $(i,j)$ is $x_{ij}$, and $Y\in\Mat_{n \times m}(R)$ the matrix whose entry at position $(j,i)$ is $y_{ji}$. It follows from the relations (2) that $X$ and $Y$ satisfy (1), and therefore we have $R^m\cong R^n$. The $K$-algebra $R$ constructed above is called the {\it Leavitt algebra of type $(m,n)$} and is denoted by $L(m,n)$. 

\begin{comment}
In \cite{vitt62} Leavitt claimed (without giving a proof) that $L(m,n)$ has no zero divisors if $K=\mathbb{F}_2$ and $m\geq 2$. In \cite{vitt65} he proved that $L(m,n)$ is a simple ring if $m=1$. In \cite{cohn66} Cohn showed that $L(m,n)$ has no zero divisors if $m\geq 2$ (covering Leavitt's claim mentioned above). Moreover, he proved that $L(m,n)$ has module type $(m,n)$ if $K=\Z$. Recall that a unital ring $R$ has the {\it Invariant Basis Number (IBN) property} if there are no $m\neq n\in \N$ such that $R^m\cong R^n$ as left $R$-modules. A unital ring $R$ not having the IBN property has {\it module type $(m,n)$} if $(m,n)=\min\{(m',n')\in\N\times \N\mid m'<n',~R^m\cong R^n\}$ with respect to the lexicographical order. 
\end{comment}

\section{Unweighted and weighted Leavitt path algebras}

\subsection{Unweighted Leavitt path algebras}

%\begin{definition}[{\sc Connected components, connected weigthed graph}]
%If $u,v\in E^0$ and there is a path $p$ in the double graph $E_d$ such that $s_d(p)=u$ and $r_d(p)=v$, then we write $u\geq_d v$. Clearly $\geq_d$ is an equivalence relation on $E^0$. The equivalence classes of $\geq_d$ are called {\it connected components}. $(E,w)$ is called {\it connected} if there is only one connected component. 
%\end{definition}

%\begin{definition}[{\sc Weighted subgraph defined by hereditary vertex set}]
%Let $H\subseteq E^0$ be a hereditary subset. Set $E_H^0:=H$, $E_H^1:=\{e\in E^1\mid s(e)\in H\}$, $r_H:=r|_{E_H^1}$, $s_H=s|_{E_H^1}$ and $w_H:=w|_{E_H^1}$. Then $E_H:=(E_H^0,E_H^1,s_H,r_H)$ is a directed graph and $(E_H,w_H)$ a weighted graph. We call $(E_H,w_H)$ the {\it weighted subgraph of $(E,w)$ defined by $H$}.   
%\end{definition}

\begin{definition}\label{deflpa}
Let $E$ be a graph. The $K$-algebra $L(E)=L_K(E)$ presented by the generating set $\{v,e,e^*\mid v\in E^0,e\in E^1\}$ and the relations
\begin{enumerate}[(i)]
\item $uv=\delta_{uv}u\quad(u,v\in E^0)$,
\medskip
\item $s(e)e=e=er(e),~r(e)e^*=e^*=e^*s(e)\quad(e\in E^1)$,
\medskip
\item $\sum_{e\in s^{-1}(v)}ee^*= v\quad(v\in E_{\reg}^0)$ and
\medskip
\item $e^*f= \delta_{ef}r(e)\quad(v\in E_{\reg}^0, e,f\in s^{-1}(v))$
\end{enumerate}
is called the {\it (unweighted) Leavitt path algebra} of $E$. 
\end{definition}

\begin{remark}\label{remlpa1}
Recall from Section 2 that $E_d$ denotes the double graph of $E$. The Leavitt path algebra $L(E)$ is isomorphic to the quotient of the path algebra $P(E_d)$ by the ideal generated by the relations (iii) and (iv) in Definition \ref{deflpa}.
\end{remark}

\begin{remark}\label{remlpa2}
The relations (iii) and (iv) in Definition \ref{deflpa} can be expressed using matrices: Let $K\X$ be the free $K$-algebra on the set $X=\{v,e,e^*\mid v\in E^0,e\in E^1\}$. For any $v\in E^0_{\reg}$ write $s^{-1}(v)=\{e^{v,1},\dots,e^{v,n(v)}\}$ and define the matrices
\[X_v:=\begin{pmatrix}e^{v,1}&\dots&e^{v,n(v)}\end{pmatrix}\in \Mat_{1\times n(v)}(K\X)\]
and
\[Y_v:=\begin{pmatrix}(e^{v,1})^*\\\vdots\\(e^{v,n(v)})^*\end{pmatrix}\in \Mat_{n(v)\times 1}(K\X).\]
Then 
\[(iii)\Leftrightarrow  \Big(X_vY_v=v~\forall v\in E_{\reg}^0\Big)\]
and 
\[(iv)\Leftrightarrow \Big(Y_vX_v=\diag(r(e^{v,1}),\dots,r(e^{v,n(v)}))~\forall v\in E_{\reg}^0\Big). \]
\end{remark}

\begin{example}\label{lpapp}
Let $n>1$ and $E$ be the graph 
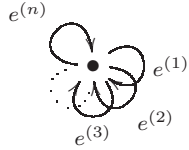
\[
\xymatrix{
& \bullet\ar@{.}@(l,d) \ar@(ur,dr)^{e^{(1)}} \ar@(r,d)^{e^{(2)}} \ar@(dr,dl)^{e^{(3)}} \ar@(l,u)^{e^{(n)}}& 
}
\]
with one vertex and $n$ loops. Then there is an isomorphism from the Leavitt algebra $L(1,n)$ to the Leavitt path algebra $L(E)$ mapping $x_{1j}$ to $e^{(j)}$ and $y_{j1}$ to $(e^{(j)})^*$ for any $1\leq j\leq n$. 
\end{example}

\begin{comment}
Let $E$ be a graph and $A$ a $K$-algebra. An {\it $E$-family} in $A$ is a subset $X=\{\alpha_v, \beta_{e}, \gamma_{e}\mid v\in E^0, e\in E^1\}\subseteq A$ such that\\
\vspace{-0.1cm}
\begin{enumerate}[(i)]
\item the $\alpha_v$'s are pairwise orthogonal idempotents, 
\medskip
\item
$\alpha_{s(e)}\beta_{e}=\beta_{e}=\beta_{e}\alpha_{r(e)},~\alpha_{r(e)}\gamma_{e}=\gamma_{e}=\gamma_{e}\alpha_{s(e)}\quad(e\in E^1)$,
\medskip
\item $\gamma_{e}\beta_{f}= \delta_{ef}\alpha_{r(e)}\quad(e,f\in E^1)$ and
\medskip
\item $\sum_{e\in s^{-1}(v)}\beta_{e}\gamma_{e}= \alpha_{v}\quad(v\in E^0_{\reg})$.
\end{enumerate}
By the relations defining $L(E)$, there exists a unique $K$-algebra homomorphism $\phi: L(E)\rightarrow A$ such that $\phi(v)=\alpha_v$, $\phi(e)=\beta_{e}$ and $\phi(e^*)=\gamma_{e}$ for all $v\in E^0$ and $e\in E^1$. We will refer to this as the {\it Universal Property of $L(E)$}.
\end{comment}

\subsection{Weighted Leavitt path algebras}

\begin{definition}\label{def3}
Let $(E,w)$ be a weighted graph.  The $K$-algebra $L(E,w)=L_K(E,w)$ presented by the generating set $\{v,e_i,e_i^*\mid v\in E^0, e\in E^1, 1\leq i\leq w(e)\}$ and the relations
\begin{enumerate}[(i)]
\item $uv=\delta_{uv}u\quad(u,v\in E^0)$,
\medskip
\item $s(e)e_i=e_i=e_ir(e),~r(e)e_i^*=e_i^*=e_i^*s(e)\quad(e\in E^1, 1\leq i\leq w(e))$,
\medskip
\item 
$\sum_{e\in s^{-1}(v)}e_ie_j^*= \delta_{ij}v\quad(v\in E^0_{\reg},1\leq i, j\leq w(v))$ and
\medskip 
\item $\sum_{1\leq i\leq w(v)}e_i^*f_i= \delta_{ef}r(e)\quad(v\in E^0_{\reg}, e,f\in s^{-1}(v))$
\end{enumerate}
is called the {\it weighted Leavitt path algebra} of $(E,w)$. In relations (iii) and (iv) we set $e_i$ and $e_i^*$ zero whenever $i > w(e)$. 
\end{definition}

\begin{remark}\label{remwlpa1}
Recall from Section 2 that $\hat E$ denotes the unweighted graph associated to $(E,w)$. We denote by $\hat E_d$ the double graph of $\hat E$. The weighted Leavitt path algebra $L(E,w)$ is isomorphic to the quotient of the path algebra $P(\hat E_d)$ by the ideal generated by the relations (iii) and (iv) in Definition \ref{def3}.
\end{remark}

\begin{remark}\label{remwlpa2}
The relations (iii) and (iv) in Definition \ref{def3} can be expressed using matrices: Let $K\X$ be the free $K$-algebra on the set $X=\{v,e_i,e_i^*\mid v\in E^0, e\in E^1, 1\leq i\leq w(e)\}$. For any $v\in E^0_{\reg}$ write $s^{-1}(v)=\{e^{v,1},\dots,e^{v,n(v)}\}$ and define the matrices
\[X_v:=\begin{pmatrix}e^{v,1}_1&\dots&e^{v,n(v)}_1\\\vdots&&\vdots\\e^{v,1}_{w(v)}&\dots&e^{v,n(v)}_{w(v)}\end{pmatrix}\in \Mat_{w(v)\times n(v)}(K\X)\]
and
\[Y_v:=\begin{pmatrix}(e^{v,1}_1)^*&\dots&(e^{v,1}_{w(v)})^*\\\vdots&&\vdots\\(e^{v,n(v)}_{1})^*&\dots&(e^{v,n(v)}_{w(v)})^*\end{pmatrix}\in \Mat_{n(v)\times w(v)}(K\X).\]
Here we set $e^{v,j}_i$ and $(e^{v,j}_i)^*$ zero whenever $i > w(e^{v,j})$.
Then 
\[(iii) \Leftrightarrow \Big(X_vY_v=\diag(v,\dots,v)~\forall v\in E_{\reg}^0\Big)\]
and 
\[(iv)\Leftrightarrow \Big(Y_vX_v=\diag(r(e^{v,1}),\dots,r(e^{v,n(v)}))~\forall v\in E_{\reg}^0\Big). \]
\end{remark}

\begin{example}\label{exex1}
If $(E,w)$ is a weighted graph such that $w(e)=1$ for all $e \in E^{1}$, then $L(E,w)$ is isomorphic to the unweighted Leavitt path algebra $L(E)$. 
\end{example}

\begin{example}\label{wlpapp}
Let $1\leq m<n$ and $(E,w)$ be the weighted graph 
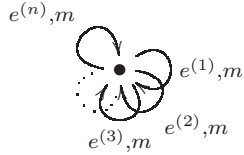
\begin{equation*}
\xymatrix{
& \bullet \ar@{.}@(l,d) \ar@(ur,dr)^{e^{(1)},m} \ar@(r,d)^{e^{(2)},m} \ar@(dr,dl)^{e^{(3)},m} \ar@(l,u)^{e^{(n)},m}& 
}
\end{equation*}
with one vertex and $n$ loops $e^{(1)},\dots,e^{(n)}$ each of which has weight $m$. Then there is an isomorphism from the Leavitt algebra $L(m,n)$ to the weighted Leavitt path algebra $L(E,w)$ mapping $x_{ij}$ to $e_i^{(j)}$ and $y_{ji}$ to $(e^{(j)}_i)^*$ for any $1\leq i\leq m$ and $1\leq j\leq n$. 
\end{example}

\begin{comment}
Let $(E,w)$ be a weighted graph and $A$ a $K$-algebra. An {\it $(E,w)$-family} in $A$ is a subset $X=\{\alpha_v, \beta_{e,i}, \gamma_{e,i}\mid v\in E, e\in E^1, 1\leq i\leq w(e)\}\subseteq A$ such that\\
\vspace{-0.1cm}
\begin{enumerate}[(i)]
\item the $\alpha_v$'s are pairwise orthogonal idempotents, 
\medskip
\item
$\alpha_{s(e)}\beta_{e,i}=\beta_{e,i}=\beta_{e,i}\alpha_{r(e)},~\alpha_{r(e)}\gamma_{e,i}=\gamma_{e,i}=\gamma_{e,i}\alpha_{s(e)}\quad(e\in E^1, 1\leq i\leq w(e))$,
\medskip
\item $\sum_{1\leq i\leq \max\{w(e),w(f)\}}\gamma_{e,i}\beta_{f,i}= \delta_{ef}\alpha_{r(e)}\quad(e,f\in E^1)$ and
\medskip
\item $\sum_{e\in s^{-1}(v)}\beta_{e,i}\gamma_{e,j}= \delta_{ij}\alpha_{v}\quad(v\in E_{\reg}^0,1\leq i, j\leq w(v))$.
\end{enumerate}
In relations (iii) and (iv) we set $\beta_{e,i}$ and $\gamma_{e,i}$ zero whenever $i > w(e)$. By the relations defining $L(E,w)$, there exists a unique $K$-algebra homomorphism $\phi: L(E,w)\rightarrow A$ such that $\phi(v)=\alpha_v$, $\phi(e_{i})=\beta_{e,i}$ and $\phi(e^*_{i})=\gamma_{e,i}$ for all $v\in E^0$, $e\in E^1$ and $1\leq i\leq w(e)$. We will refer to this as the {\it Universal Property of $L(E,w)$}. The proof of the Universal Property for a weighted Leavitt path algebra is analogous to the proof in the unweighted case.
\end{comment}

Recall that a ring $R$ is said to have a {\it set of local units} $X$ in case $X$ is a set of idempotents in $R$ having the property that for each finite subset $S\subseteq R$ there exists an $x\in X$ such that $xsx=s$ for any $s\in S$. If $(E,w)$ is a weighted graph, then $L(E,w)$ has a set of local units, namely the set of all finite sums of distinct elements of $E^0$. If $E$ is finite, then $L(E,w)$ is a unital ring with $\sum_{v\in E^0} v=1$.

There is an involution $*$ on $L(E,w)$ mapping $k\mapsto k$, $v\mapsto v$, $e_i\mapsto e_i^*$ and $e_i^*\mapsto e_i$ for any $k\in K$, $v\in E^0$, $e\in E^1$ and $1\leq i\leq w(e)$. Moreover, one can define a grading on $L(E,w)$ as follows. Set $\lambda:=\sup\{w(e) \mid e \in E^{1}\}$ if this supremum is finite and otherwise $\lambda:=\omega$ where $\omega$ is the smallest infinite ordinal. Let $\mathbb Z^\lambda$ denote the sum of $\lambda$-many copies of $\mathbb Z$. There is a $\mathbb Z^\lambda$-grading on $L(E,w)$ such that $\deg(v):=0$, $\deg(e_i):=\alpha_i$ and $\deg(e_i^*):=-\alpha_i$ for any $v\in E^0$, $e \in E^{1}$ and $1\leq i\leq w(e)$. Here $\alpha_i$ denotes the element of $\mathbb Z^\lambda$ whose $i$-th component is $1$ and whose other components are $0$. We will refer to this grading as the {\it standard grading} of $L(E,w)$. For proofs of the claims made in this and the previous paragraph see \cite[Proposition 5.7]{hazrat13}. 

\section{Linear bases and the Gelfand-Kirillov dimension}

\subsection{Linear bases for weighted Leavitt path algebras}
Let $(E,w)$ be a weighted graph and %We recall the basis result of \cite{hazrat-preusser}.
set $X:=\{v,e_i,e_i^*\mid v\in E^0,e\in E^1,1\leq i\leq w(e)\}$. Let $\X$ denote the set of all finite, nonempty words over $X$ and set $\overline{\X}:=\X\cup\{\text{empty word}\}$. Together with juxtaposition of words $\X$ becomes a semigroup and $\overline{\X}$ a monoid. If $A,B\in \overline{\X}$, then $B$ is called a {\it prefix of $A$} if there is a $D \in \overline{\X}$ such that $A=BD$, a {\it suffix of $A$} if there is a $C \in \overline{\X}$ such that $A=CB$, and a {\it subword of $A$} if there are $C,D \in\overline{\X}$ such that $A=CBD$. 

For any $v\in E_{\reg}^0$ we fix an edge $e^v\in s^{-1}(v)$ such that $w(e^v)=w(v)$. The edges $e^v~(v\in E_{\reg}^0)$ are called {\it special edges}. The words $e^v_i(e^v_j)^*~(v\in E_{\reg}^0,1\leq i,j\leq w(v))$ and $e^*_1f_1~(v\in E_{\reg}^0, e,f\in s^{-1}(v))$ in $\X$ are called {\it forbidden}. Recall from Section 2 that $\hat E$ denotes the unweighted graph associated to $(E,w)$ and $\hat E_d$ the double graph of $\hat E$. We call a path in $\hat E_d$ a {\it d-path}. A {\it normal d-path} or {\it nod-path} is a d-path such that none of its subwords is forbidden.

Let $K\X$ be the free $K$-algebra on $X$ (i.e. the semigroup $K$-algebra of $\X$). Then $L(E,w)$ is the quotient $K\X/I$ where $I$ is the ideal generated by the relations (i)-(iv) in Definition \ref{def3}. Let $K\X_{\nod}$ be the linear subspace of $K\X$ spanned by the nod-paths.

\begin{theorem}[{\cite[Theorem 16]{hazrat-preusser}}] \label{thmbasis}
The canonical map $K\X_{\nod}\rightarrow L(E,w)$ is an isomorphism of $K$-vector spaces. In particular the images of the nod-paths under this map form a linear basis for $L(E,w)$.
\end{theorem} 
\begin{proofs}
Consider the relations (i')-(v') below.
\begin{enumerate}[(i')]
\item For any $u,v \in E^0$, \[uv = \delta_{uv}u.\]
\medskip

\item For any $v \in E^0$, $e\in E^{1}$ and $1\leq i \leq w(e)$,
\begin{align*}
ve_i&=\delta_{vs(e)}e_i,\\
e_iv&=\delta_{vr(e)}e_i,\\
ve_i^*&=\delta_{vr(e)}e_i^*,\\
e_i^*v&=\delta_{vs(e)}e_i^*.
\end{align*}

\medskip

\item For any $e,f\in E^{1}$, $1\leq i \leq w(e)$ and $1\leq j \leq w(f)$,

\begin{align*}
e_if_j&=0\text{ if  } r(e)\neq s(f),\\
e_i^*f_j&=0\text{ if  } s(e)\neq s(f),\\
e_if_j^*&=0\text{ if  } r(e)\neq r(f),\\
e_i^*f_j^*&=0\text{ if  } s(e)\neq r(f).\\
\end{align*}

\medskip

\item For any $v\in E_{\reg}^0$ and $1\leq i, j\leq w(v)$, 
\[e^v_i(e^v_j)^*= \delta_{ij}v-\sum_{e\in s^{-1}(e)\setminus\{e^v\}}e_ie_j^*.\]

\medskip 

\item For any $v\in E_{\reg}^0$ and $e,f\in s^{-1}(v)$, 
\[e_1^*f_1= \delta_{ef}r(e)-\sum_{2\leq i\leq w(v)}e_i^*f_i.\]

\end{enumerate}
In relations (iv') and (v'), we set $e_i$ and $e_i^*$ zero whenever $i > w(e)$. Clearly the relations (i')-(v') generate the same ideal $I$ of $K\X$ as the relations (i)-(iv) in Definition \ref{def3}. Let $S$ be the set of all pairs $\sigma=(W_\sigma,f_\sigma)$ where $W_\sigma$ equals the left hand side of an equation in (i')-(v') and $f_\sigma$ the corresponding right hand side. For any $\sigma\in S$ and $A,B \in\overline{\X}$ let $r_{A\sigma B}$ denote the endomorphism of the $K$-vector space $K\X$ that maps $AW_\sigma B$ to $Af_\sigma B$ and fixes all other elements of $\X$. The set $S$ is called a {\it reduction system} and the maps $r_{A\sigma B}$ are called {\it reductions}. %A finite sequence of reductions $r_1,\dots, r_n$ is called {\it final} on $a\in K\X$ if $r_n\dots r_1(a)\in K\X_{\nod}$.

One can show that there is a semigroup partial ordering on $\X$ compatible with $S$ and having the descending chain condition. Moreover, all ambiguities of $S$ are resolvable (cf. \cite[Section 1]{bergman78} or \cite[Section 2]{hazrat-preusser}). Now G. Bergman's Diamond Lemma implies the following. Let $a\in L(E,w)=K\X/I$ and $b\in K\X$ a representative of $a$. If one successively applies reductions to $b$ (not leaving the given element of $K\X$ invariant), then eventually one will obtain an element $c\in K\X_{\nod}$. It follows from the Diamond Lemma that the element $c$ does not depend on the choice of the representative or the reductions. We call $\NF(a):=c$ the {\it normal form} of $a$. The map $\NF:L(E,w)\to K\X_{\nod}$ is an isomorphism of $K$-vector spaces. Its inverse is the canonical map $K\X_{\nod}\rightarrow L(E,w)$. 
\end{proofs}

\begin{comment}
\begin{remark}\label{rembasis}
Let $\NF:L(E,w)\to K\X_{\nod}$ be the inverse of the isomorphism $K\X_{\nod}\rightarrow L(E,w)$. If $a\in L(E,w)=K\X/I$, then $\NF(a)$ is called the {\it normal form} of $a$. Let $S$ denote the reduction system defined relations (i')-(v'). The normal form $\NF(a)$ can be obtained as follows:\\
\\
{\bf Step 1} Choose an arbitrary representative $b$ of $a$ in $K\X$.\\
\\
{\bf Step 2} Apply successively reductions to $b$ corresponding to the relations (i')-(v')  until arriving at an element of $K\X_{\nod}$ (i.e. replace forbidden words $W_\sigma$ by the corresponding right hand sides $f_\sigma$ until it is no longer possible).  \\
\\
\end{remark}
\end{comment}

\subsection{The Gelfand-Kirillov dimension of a weighted Leavitt path algebra}
First we recall some facts about the growth of algebras. Let $A$ be a nonzero, finitely generated $K$-algebra and $V$ a finite-dimensional subspace of $A$ that generates $A$ as an algebra. For $n\geq 1$ let $V^n$ denote the linear span of the set $\{v_1\dots v_k\mid k\leq n, v_1,\dots,v_k\in V\}$. Then 
\[V =V^1\subseteq V^2\subseteq V^3\subseteq \dots, \quad A =\bigcup_{n\in \N}V^n\text{ and }d_V(n):=\dim V^n<\infty.\] 
Given functions $f, g:\N\rightarrow \R^+$, we write $f\preccurlyeq g$ if there is a $c\in\N$ such that $f(n)\leq cg(cn)$ for all $n$. If $f\preccurlyeq g$ and $g\preccurlyeq f$, then the functions $f, g$ are called {\it asymptotically equivalent} and we write $f\sim g$. If $W$ is another finite-dimensional subspace that generates $A$, then $d_V\sim d_W$. The {\it Gelfand-Kirillov dimension} or {\it GK dimension} of $A$ is defined as
\[\GKdim A := \limsup_{n\rightarrow \infty}\log_nd_V(n).\]
The definition of the GK dimension does not depend on the choice of the finite-dimensional generating space $V$. If $d_V\preccurlyeq n^m$ for some $m\in \N$, then $A$ is said to have {\it polynomial growth} and we have $\GKdim A \leq m$. If $d_V\sim a^n$ for some real number $a>1$, then $A$ is said to have {\it exponential growth} and we have $\GKdim A =\infty$. If $A$ is not finitely generated as a $K$-algebra, then the GK dimension of $A$ is defined as
\[\GKdim(A) := \sup\{\GKdim(B)\mid B \text{ is a finitely generated subalgebra of }A\}.\]

%If $X$ is a set, we denote by $\X$ the set of all words over $X$ including the empty word. Together with juxtaposition $\X$ is a monoid. If $A,B\in \X$, then we write $A|B$ if there is a $C\in\X$ such that $AC=B$.
% If $p\overset{\emptyset}{\Longrightarrow} q$ and $q\overset{\emptyset}{\Longrightarrow} p$ (resp. $p\overset{\nod}{\Longrightarrow} q$ and $q\overset{\nod}{\Longrightarrow} p$ or $p\Longrightarrow q$ and $q\Longrightarrow p$), we write $p\overset{\emptyset}{\Longleftrightarrow} q$ (resp. $p\overset{\nod}{\Longleftrightarrow} q$ or $p
%Let $A$ be a word over some alphabet. We call a subword $B$ of $A$ {\it proper} if $B\neq A$.

Let $(E,w)$ be a weighted graph. A {\it nod$^2$-path} in $(E,w)$ is a nod-path $p$ such that $p^2$ is a nod-path. A {\it quasicycle} is a nod$^2$-path $p$ such that none of the subwords of $p^2$ of length $<|p|$ is a nod$^2$-path. If $p = x_1\dots x_n$ is a quasicycle, then $x_i\neq x_j$ for all $i\neq j$ by \cite[Remark 16(a)]{preusser}. It follows that there is only a finite number of quasicycles if $(E,w)$ is finite. Note that if $p=x_1\dots x_n$ is a nod$^2$-path (resp. quasicycle), then $p^*:=x_n^*\dots x_1^*$ is a nod$^2$-path (resp. quasicycle).

In general it is not so easy to read off the quasicycles from a weigthed graph. But there is the following algorithm for finding all the quasicycles in a finite weigthed graph: For any vertex $v$ list all the d-paths $p=x_1\dots x_n$ starting and ending at $v$ and having the property that $x_i\neq x_j$ for any $i\neq j$ (there are only finitely many d-paths with this property). Now delete from that list any $p$ such that $p^2$ is not a nod-path. Next delete from the list any $p$ such that $p^2$ has a subword $q$ of length $|q|<|p|$ such that $q^2$ is a nod-path. The remaining d-paths on the list are precisely the quasicycles starting (and ending) at $v$.

\begin{example}\label{exqu}
Suppose $(E,w)$ is the weighted graph 
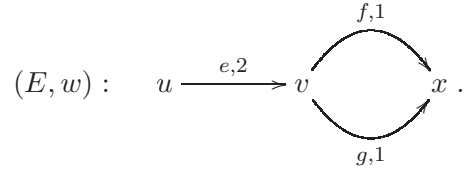
\[
(E,w):\quad\xymatrix@C+15pt{
u\ar[r]^{e,2}&  v\ar@/^1.7pc/[r]^{f,1}\ar@/_1.7pc/[r]_{g,1}  & x  
}.
\]
Then the associated unweighted graph $\hat E$ and its double graph $\hat E_d$ are given by
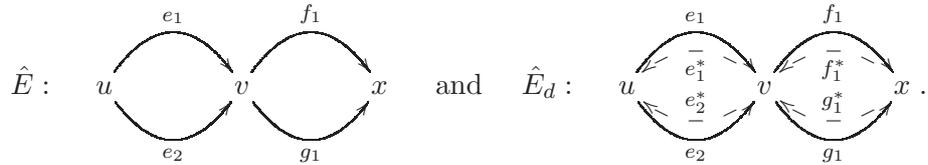
\[
\hat E:\quad\xymatrix@C+15pt{
u\ar@/^1.7pc/[r]^{e_1}\ar@/_1.7pc/[r]_{e_2}&  v\ar@/^1.7pc/[r]^{f_1}\ar@/_1.7pc/[r]_{g_1}  & x 
}\quad\text{ and }\quad \hat E_d:\quad\xymatrix@C+15pt{
u\ar@/^1.7pc/[r]^{e_1}\ar@/_1.7pc/[r]_{e_2}&  v\ar@{-->}@/_1.1pc/[l]^{e_1^*}\ar@{-->}@/^1.1pc/[l]_{e_2^*}\ar@/^1.7pc/[r]^{f_1}\ar@/_1.7pc/[r]_{g_1}  & x \ar@{-->}@/_1.1pc/[l]^{f_1^*}\ar@{-->}@/^1.1pc/[l]_{g_1^*} 
}.
\]
One checks easily that $p:=e_2f_1g_1^*e_2^*$ and $q:=e_2f_1g_1^*e_1^*$ are quasicycles independent of the choice of the special edge $e^v$. 
\end{example}

If $x_1\dots x_n$ is a finite, nonempty word over some alphabet, then we call the words $x_{m+1}\dots x_{n}x_1\dots x_{m}~(1\leq m\leq n)$ {\it shifts} of $x_1\dots x_n$. One checks easily that any shift of a quasicycle is a quasicycle. If $p$ and $q$ are quasicycles, then we write $p\approx q$ iff $q$ is a shift of $p$. Then $\approx$ is an equivalence relation on the set of all quasicycles.

Let $p$ and $q$ be nod-paths. If there is a nod-path $o$ such that $p$ is not a prefix of $o$ and $poq$ is a nod-path, then we write $p\overset{\nod}{\Longrightarrow} q$. If $pq$ is a nod-path or $p\overset{\nod}{\Longrightarrow} q$, then we write $p\Longrightarrow q$. A nod-path $p$ is called {\it selfconnected} if $p\overset{\nod}{\Longrightarrow}p$. A sequence $p_1,\dots,p_k$ of quasicycles such that $p_i\not\approx p_j$ for any $i\neq j$ is called a {\it chain of length $k$} if $p_1 \Longrightarrow p_2\Longrightarrow \dots \Longrightarrow p_k$. 

If $(E,w)$ is a weighted graph, we denote by $E'$ the subset of $\{e_i,e_i^*\mid e\in E^1, 1\leq i\leq w(e)\}$ consisting of all the elements which are not a letter of a quasicycle. We denote by $P'$ the set of all nod-paths which are composed from elements of $E'$. If $(E,w)$ is finite, then $|P'|<\infty$ by \cite[Lemma 21]{preusser}.

\begin{lemma}[{\cite[Lemma 22]{preusser}}]\label{lemnewgk}
Let $(E,w)$ be a weighted graph. If there is no selfconnected quasicycle, then any nontrivial nod-path $\alpha$ can be written as
\begin{equation*}
\alpha=o_1p_1^{l_1}q_1o_2p_2^{l_2}q_2\dots o_kp_k^{l_k}q_ko_{k+1}
\end{equation*}
where $k\geq 0$, $o_i$ is the empty word or $o_i\in P'~(1\leq i \leq k+1)$, $p_1,\dots,p_k$ is a chain of quasicycles, $l_i$ is a nonnegative integer $(1\leq i \leq k)$, and $q_i\neq p_i$ is a prefix of $p_i~(1\leq i \leq k)$.
\end{lemma}

\begin{theorem}[{\cite[Lemma 19, Theorem 23]{preusser}}]\label{thmgk}
Let $(E,w)$ be a finite weighted graph. Then:
\begin{enumerate}[(i)]
\item If there is a selfconnected quasicycle, then $L(E,w)$ has exponential growth and hence \[\GKdim L(E,w)=\infty.\]
\item If there is no selfconnected quasicycle, then $L(E,w)$ has polynomial growth. In this case \[\GKdim L(E,w)=d\]
where $d$ is the maximal length of a chain of quasicycles.
\end{enumerate}
\end{theorem}
\begin{proofs}
Let $V$ denote the finite-dimensional subspace of $L(E,w)$ spanned by $\{v,e_i,e_i^*\mid v\in E^0,e\in E^1,1\leq i\leq w(e)\}$. By Theorem \ref{thmbasis} and its proof, the nod-paths of length $\leq n$ form a basis for $V^n$ (since reductions shorten the length of words). Hence $d_V(n)=\dim V^n=\#\{\text{nod-paths of length}\leq n\}$.

First suppose that there is a selfconnected quasicycle $p$. Let $o$ be a nod-path such that $p$ is not a prefix of $o$ and $pop$ is a nod-path. For a fixed $n\in \N$ consider the nod-paths
\begin{equation}
p^{i_1}op^{i_2}\dots op^{i_k}
\end{equation}
where $k,i_1,\dots,i_k\in \N$ satisfy
\begin{equation}
(i_1+\dots+i_k)|p|+(k-1)|o|\leq n.
\end{equation}
One can show that different solutions $(k,i_1,\dots,i_k)$ of inequality (4) define different nod-paths in (3). Since the number of solutions of (4) is $\sim 2^n$, we obtain $2^n\preccurlyeq d_V$. Thus $L(E,w)$ has exponential growth.  

Suppose now that there is no selfconnected quasicycle. Fix an $n\in\N$. By Lemma \ref{lemnewgk} we can write any nontrivial nod-path $\alpha$ of length $\leq n$ as
\begin{equation}
\alpha=o_1p_1^{l_1}q_1o_2p_2^{l_2}q_2\dots o_kp_k^{l_k}q_ko_{k+1}
\end{equation}
where $k\geq 0$, $o_i$ is the empty word or $o_i\in P'~(1\leq i \leq k+1)$, $p_1,\dots,p_k$ is a chain of quasicycles, $l_i$ is a nonnegative integer $(1\leq i \leq k)$, and $q_i\neq p_i$ is a prefix of $p_i~(1\leq i \leq k)$. Clearly 
\begin{equation}
l_1|p_1| +\dots+l_k|p_k| \leq n
\end{equation}
since $|\alpha|\leq n$. Now fix a chain $p_1,\dots,p_k$ of quasicycles and further $o_i$'s and $q_i$'s as above. The number of solutions $(l_1,\dots,l_k)$ of inequality (6) is $\sim n^k$. This implies that the number of nod-paths $\alpha$ of length $\leq n$ that can be written as in (5) (corresponding to the choice of the $p_i$'s, $o_i$'s and $q_i$'s) is $\preccurlyeq n^k \leq n^d$. Since there are only finitely many quasicycles and finitely many choices for the $o_i$'s and $q_i$'s (note that $|P'|<\infty$), we obtain $d_V\preccurlyeq n^d$.

It remains to show that $n^d \preccurlyeq d_V$. Choose a chain $p_1,\dots,p_d$ of length $d$. Then $p_1o_1p_2\dots o_{d-1}p_d$ is a nod-path for some $o_1,\dots,o_{d-1}$ such that for any $i\in\{1,\dots,d-1\}$, $o_i$ is either the empty word or a nod-path such that $p_i$ is not a prefix of $o_i$. Consider the nod-paths
\begin{equation}
p_1^{l_1}o_1p_2^{l_2}\dots o_{d-1}p_d^{l_d}
\end{equation}
where $l_1,\dots,l_d\in \N$ satisfy
\begin{equation}
l_1|p_1| +\dots+l_d|p_d|+|o_1|+\dots+|o_{d-1}|\leq n.
\end{equation}
One can show that different solutions $(l_1,\dots,l_d)$ of inequality (8) define different nod-paths in (7). The number of solutions of (8) is $\sim n^d$ and thus $n^d \preccurlyeq d_V$.
\end{proofs}

\begin{example}
Let $1\leq m<n$ and $(E,w)$ be the weighted graph 
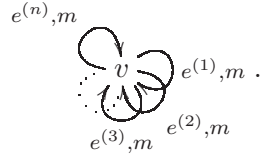
\begin{equation*}
\xymatrix{
& v \ar@{.}@(l,d) \ar@(ur,dr)^{e^{(1)},m} \ar@(r,d)^{e^{(2)},m} \ar@(dr,dl)^{e^{(3)},m} \ar@(l,u)^{e^{(n)},m}& 
}.
\end{equation*}
As mentioned in Example \ref{wlpapp}, $L(E,w)$ is isomorphic to the Leavitt algebra $L(m,n)$. Clearly $p:=e^{(1)}_1$ is a quasicycle. Since $e^{(1)}_1e^{(2)}_1e^{(1)}_1$ is a nod-path, $p$ is selfconnected. Hence $\GKdim L(E,w)=\infty$ by Theorem \ref{thmgk}.
\end{example}

\begin{example}
Consider the weighted graph 
\[
(E,w):\quad\vcenter{\vbox{
\xymatrix@C+15pt@R+15pt{
u\ar[rd]_{g,1}\ar[r]^{e,2}&  v\ar[d]^{f,1}  \\&\hspace{0.2cm}x~.
}}}
\]
We leave it to the reader to check that up to $\approx$-equivalence the only quasicycles are $p=e_2f_1g_1^*$ and $p^*=g_1f_1^*e_2^*$. Moreover, there is no selfconnected quasicycle and the maximal length of a chain of quasicycles is $2$, see \cite[Example 31]{preusser}. Thus $\GKdim L(E,w)=2$ by Theorem \ref{thmgk}.
\end{example}

Next we want to determine the GK dimensions of weighted Leavitt path algebras of non-finite weighted graphs. In order to that we need the following definition.

\begin{definition}\label{defcompl}
A homomorphism $\phi:(\tilde E,\tilde w)\to (E,w)$ of weighted graphs is called {\it complete} if $\phi^0$ is injective and $\phi^1|_{\tilde s^{-1}(v)}: \tilde s^{-1}(v)\to s^{-1}(\phi^0(v))$ is a bijection for any $v\in \tilde E^0_{\reg}$. A weighted subgraph $(\tilde E,\tilde w)$ of a weighted graph $(E,w)$ is called {\it complete} if $\tilde s^{-1}(v)=s^{-1}(v)$ for any $v\in \tilde E^0_{\reg}$.
\end{definition}

Note that if $(\tilde E,\tilde w)$ is a complete weighted subgraph of a weighted graph $(E,w)$, then the inclusion homomorphism $\phi:(\tilde E,\tilde w)\to (E,w)$ is complete. We denote by $\G$ the category whose objects are the weighted graphs and whose morphisms are the complete weighted graph homomorphisms. If $\phi:(\tilde E,\tilde w)\to (E,w)$ is a morphism in $\G$, then it induces an algebra homomorphism $L(\phi):L(\tilde E,\tilde w)\to L(E,w)$. Since one can choose the special edges in $(\tilde E,\tilde w)$ and $(E,w)$ in such a way that $L(\phi)$ maps distinct nod-paths in $L(\tilde E,\tilde w)$ to distinct nod-paths in $L(E,w)$, the homomorphism $L(\phi)$ is injective. Let $\A$ denote the category of $K$-algebras. Then we obtain a functor $L:\G\to \A$ which commutes with direct limits.

\begin{lemma}[{\cite[Lemma 5.19]{hazrat13}}]\label{lemdirlim}
Let $(E,w)$ be a weighted graph and $\{(E_i,w_i)\}$ the direct system of all finite complete weighted subgraphs of $(E,w)$. Then $(E,w)=\varinjlim_i (E_i,w_i)$ and therefore $L(E,w)=\varinjlim_i L(E_i,w_i)$.
\end{lemma}

We need the following result by J. Moreno-Fernandez and M. Molina.

\begin{theorem}[{\cite[Theorem 3.1]{moreno-molina}}]\label{thmmomo}
Let $A=\bigcup_{i}A_i$ be a directed union of algebras. Then $\GKdim A=\sup_i \GKdim A_i$.
\end{theorem}

The theorem below follows from Lemma \ref{lemdirlim} and Theorem \ref{thmmomo}.

\begin{theorem}[{\cite[Remark 24]{preusser}}]\label{thmgkinf}
Let $(E,w)$ be a weighted graph and $\{(E_i,w_i)\}$ the direct system of all finite complete weighted subgraphs of $(E,w)$. Then $\GKdim L(E,w)=\sup_i \GKdim L(E_i, w_i)$.
\end{theorem}

\section{Weighted Leavitt path algebras that are isomorphic to unweighted Leavitt path algebras}
Throughout this section $(E,w)$ denotes a weighted graph.

\subsection{Condition (LPA)}
Recall from Section 2 that an edge in $(E,w)$ is unweighted if its weight is $1$ and weighted otherwise. $E^1_w$ denotes the set of all weighted edges in $(E,w)$. If $X$ is a set of vertices in $(E,w)$, then $T(X)$ denotes the union of all the trees of the elements of $X$. Two edges $e$ and $f$ are in line if they are equal or there is a path from $r(e)$ to $s(f)$ or there is a path from $r(f)$ to $s(e)$. Consider the following conditions: 

\begin{enumerate}[({LPA}1)]
\item Any vertex $v\in E^0$ emits at most one weighted edge.
\smallskip
\item Any vertex $v\in T(r(E^1_w))$ emits at most one edge.
\smallskip
\item If two weighted edges $e,f\in E^1_w$ are not in line, then $T(r(e))\cap T(r(f))=\emptyset$.
\smallskip
\item If $e\in E^1_w$ and $c$ is a cycle based at some vertex $v\in T(r(e))$, then $e$ belongs to $c$.
\end{enumerate}
Each of the conditions above ``forbids" a certain constellation in the weighted graph $(E,w)$. The pictures below illustrate these forbidden constellations. Symbols above or below edges indicate the weight. A dotted arrow stands for a path.
\begin{enumerate}[({LPA}1)]
\item \[\xymatrix@R-1.5pc{& \bullet\\\bullet \ar_{>1}[dr] \ar^{>1}[ur]& \\& \bullet.}\]
\item \[\xymatrix@R-1.5pc{& & & \bullet\\\bullet \ar^{>1}[r] & \bullet \ar@{.>}[r] & \bullet \ar[dr] \ar[ur]& \\& & & \bullet.}\]
\item \[\xymatrix{\bullet \ar^{>1}[r] & \bullet \ar@{.>}[r] & \bullet &\bullet \ar@{.>}[l] & \bullet. \ar_{>1}[l]}\]
\item \[\xymatrix@R-1pc@C-0.5pc{
		&&&\bullet \ar[r]&\bullet \ar[rd]&\\
		\bullet \ar^{>1}[r] & \bullet \ar@{.>}[r] &\bullet \ar[ru]&&&\bullet. \ar[ld]\\
		&&&\bullet \ar[lu]&\bullet \ar@{.}[l]}	
		\]
\end{enumerate}
We say that $(E,w)$ {\it satisfies Condition (LPA)} if it satisfies Conditions (LPA1)-(LPA4). 

\subsection{Presence of Condition (LPA)}

\begin{theorem}[{\cite[Theorem 1]{Raimund4}}]\label{thmm}
If $(E,w)$ satisfies Condition (LPA), then there is a graph $F$ such that $L(E,w)\cong L(F)$ as $K$-algebras.
\end{theorem}
\begin{proofs}
In \cite{Raimund4} the theorem was proved in two steps. \\
\\
{\bf Step 1} Suppose $(E,w)$ satisfies Condition (LPA) and set $Z:=T(r(E^1_w))$. Let $(\tilde E,\tilde w)$ be the weighted graph one obtains by replacing in $(E,w)$ each edge $e\in s^{-1}(Z)$ by $w(e)$ unweighted edges with reversed orientation. More formally, $(\tilde E,\tilde w)$ is the weighted graph defined by
\begin{align*}
\tilde E^0&=E^0,\\
\tilde E^1&=\tilde E^1_Z\sqcup \tilde E^1_{Z^c}\text{ where }\\
\tilde E^1_Z&=\{e^{(1)},\dots,e^{(w(e))}\mid e\in E^1, s(e)\in Z\},\\
\tilde E^1_{Z^c}&=\{e\mid e\in E^1,s(e)\not\in Z\},\\
\tilde s(e^{(i)})&=r(e),~\tilde r(e^{(i)})=s(e),~\tilde w(e^{(i)})=1\text{ for any }e^{(i)}\in \tilde E^1_Z,\\
\tilde s(e)&=s(e), \tilde r(e)=r(e), \tilde w(e)=w(e)\text{ for any }e\in \tilde E^1_{Z^c}.
\end{align*} 
The weighted graph $(\tilde E,\tilde w)$ has the property that ranges of weighted edges are sinks and no vertex emits or receives two distinct weighted edges. Moreover, there is a $K$-algebra isomorphism $L(E,w)\to L(\tilde E, \tilde w)$ mapping 
\begin{alignat*}{2}
v&\mapsto v&\hspace{1cm}&(v\in E^0),\\
e_i&\mapsto (e^{(i)}_{1})^*&&(e\in E^1, s(e)\in Z,1\leq i\leq w(e)),\\
e_i^*&\mapsto e^{(i)}_{1}&&(e\in E^1, s(e)\in Z,1\leq i\leq w(e)),\\
e_i&\mapsto e_i&&(e\in E^1, s(e)\not\in Z,1\leq i\leq w(e)),\\
e_i^*&\mapsto e^*_{i}&&(e\in E^1, s(e)\not\in Z,1\leq i\leq w(e)).
\end{alignat*}
For details see \cite[Proof of Lemma 9]{Raimund4}. For example suppose that $(E,w)$ is the weighted graph 
\[(E,w):\quad\xymatrix@C=1.5cm{u&v\ar[l]_{e,1}\ar[r]^{f,2} & x\ar[r]^{g,1} &y\ar[r]^{h,3}& z}\]
which satisfies Condition (LPA). Then $Z=T(r(E^1_w))=\{x,y,z\}$. If one replaces any edge $i$ in $s^{-1}(Z)$ by $w(i)$ unweighted edges with reversed orientation, one obtains the weighted graph
\[(\tilde E,\tilde w):\quad\xymatrix@C=1.5cm{u&v\ar[l]_{e,1}\ar[r]^{f,2} & x &y\ar[l]_{g^{(1)},1}& z\ar@/_1.7pc/[l]_{h^{(1)},1}\ar[l]_{h^{(2)},1}\ar@/^1.7pc/[l]_{h^{(3)},1}}\]
which has the property that ranges of weighted edges are sinks and no vertex emits or receives two distinct weighted edges. There is a $K$-algebra isomorphism $L(E,w)\to L(\tilde E, \tilde w)$ mapping $u\mapsto u,\dots, z\mapsto z$, $e_1\mapsto e_1$, $e_1^*\mapsto e_1^*$, $f_i\mapsto f_i$ and $f_i^*\mapsto f_i^*$ for any $i\in \{1,2\}$, $g_1\mapsto (g_1^{(1)})^*$, $g_1^*\mapsto g_1^{(1)}$, $h_i\mapsto (h_1^{(i)})^*$ and $h_i^*\mapsto h_1^{(i)}$ for any $i\in \{1,2,3\}$. \\
\\
{\bf Step 2.} Suppose now that $(E,w)$ is a weighted graph having the property that ranges of weighted edges are sinks and no vertex emits or receives two distinct weighted edges. Let $v\in r(E^1_w)$. Then there is a unique edge $g^v\in E^1_w$ such that $r(g^v)=v$. We replace $v$ by $w(g^v)$ vertices $v^{(1)},\dots,v^{(w(g^v))}$. We replace $g^v$ by $w(g^v)$ unweighted edges $(g^v)^{(1)},\dots, (g^v)^{(w(g^v))}$ such that $(g^v)^{(1)}$ starts in $s(g^v)$ and ends in $v^{(1)}$, and $(g^v)^{(i)}$ starts in $v^{(i)}$ and ends in $s(g^v)$ for any $2\leq i\leq w(g^v)$. Moreover, we replace any unweighted edge $e$ such that $r(e)=v$ by $w(g^v)$ unweighted edges $e^{(1)},\dots,e^{(w(g^{v}))}$ such that $e^{(i)}$ starts in $s(e)$ and ends in $v^{(i)}$ for any $1\leq i\leq w(g^{v})$. By doing this for any $v\in r(E^1_w)$ we obtain an unweighted graph $\tilde E$. More formally, $\tilde E$ is defined by
\begin{align*}
\tilde E^0&=M\sqcup N \text{ where }\\
M&=E^0\setminus r(E^1_w),\\
N&=\{v^{(1)},\dots,v^{(w(g^v))}\mid v\in r(E^1_w)\},\\
\tilde E^1&=A\sqcup B\sqcup C \sqcup D\text{ where }\\
A&=\{e\mid e\in E^1_{uw},r(e)\not\in r(E^1_w)\},\\ 
B&=\{e^{(1)},\dots,e^{(w(g^{r(e)}))}\mid e\in E^1_{uw}, r(e)\in r(E^1_w)\},\\
C&=\{e^{(1)}\mid e\in E^1_{w}\},\\
D&=\{e^{(2)},\dots,e^{(w(e))}\mid e\in E^1_{w}\},\\
\tilde s(e)&=s(e),~\tilde r(e)=r(e)\quad (e\in A),\\
\tilde s(e^{(i)})&=s(e),~\tilde r(e^{(i)})=r(e)^{(i)}\quad (e^{(i)}\in B),\\
\tilde s(e^{(1)})&=s(e),~\tilde r(e^{(1)})=r(e)^{(1)}\quad (e^{(1)}\in C),\\
\tilde s(e^{(i)})&=r(e)^{(i)},~\tilde r(e^{(i)})=s(e)\quad (e^{(i)}\in D).
\end{align*}
There is a $K$-algebra isomorphism $\phi:L(E,w)\rightarrow L(\tilde E)$ such that
\begin{align*}
\phi(v)&=\begin{cases}v,&\text{ if } v\not\in r(E^1_w),\\\sum_{i=1}^{w(g^v)}v^{(i)},&\text{ if }v\in r(E^1_w),\end{cases}\\
\phi(e_i)&=\begin{cases}e,&\text{ if } e\in E_{uw}^1,r(e)\not\in r(E^1_w),i=1,\\\sum_{j=1}^{w(g^{r(e)})}e^{(j)},&\text{ if }e\in E_{uw}^1,r(e)\in r(E^1_w),i=1,\\
e^{(1)},&\text{ if }e\in E_{w}^1,i=1,\\
(e^{(i)})^*,&\text{ if }e\in E_{w}^1,i>1,
\end{cases}\\
\phi(e_i^*)&=\begin{cases}e^*,&\text{ if } e\in E_{uw}^1,r(e)\not\in r(E^1_w),i=1,\\\sum_{j=1}^{w(g^{r(e)})}(e^{(j)})^*,&\text{ if }e\in E_{uw}^1,r(e)\in r(E^1_w),i=1,\\
(e^{(1)})^*,&\text{ if }e\in E_{w}^1,i=1,\\
e^{(i)},&\text{ if }e\in E_{w}^1,i>1
\end{cases}
\end{align*}
for any $v\in E^0$, $e\in E^1$ and $1\leq i\leq w(e)$. For details see \cite[Proof of Lemma 11]{Raimund4}. For example if $(E,w)$ is the weighted graph 
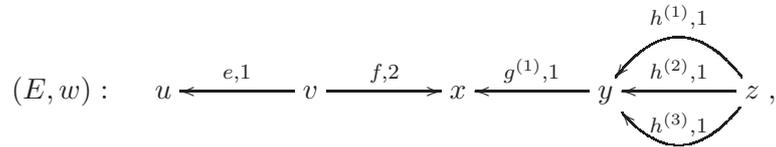
\[(E,w):\quad\xymatrix@C=1.5cm{u&v\ar[l]_{e,1}\ar[r]^{f,2} & x &y\ar[l]_{g^{(1)},1}& z\ar@/_1.7pc/[l]_{h^{(1)},1}\ar[l]_{h^{(2)},1}\ar@/^1.7pc/[l]_{h^{(3)},1}},\]
then the unweighted graph $\tilde E$ define above is given by
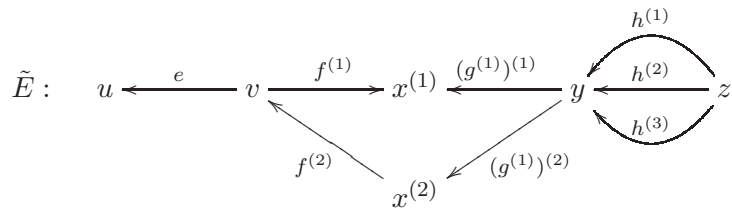
\[\tilde E:\quad\xymatrix@C=1.5cm{u&v\ar[l]_{e}\ar[r]^{f^{(1)}} & x^{(1)} &y\ar[l]_{(g^{(1)})^{(1)}}\ar[dl]^{(g^{(1)})^{(2)}}& z\ar@/_1.7pc/[l]_{h^{(1)}}\ar[l]_{h^{(2)}}\ar@/^1.7pc/[l]_{h^{(3)}}\\&&x^{(2)}\ar[ul]^{f^{(2)}}&&}\]
and we have $L(E,w)\cong L(\tilde E)$.\\
\\
The theorem follows from Steps 1 and 2 above.
\end{proofs}

\begin{example}\label{exlpa1}
Suppose $(E,w)$ is the weighted graph
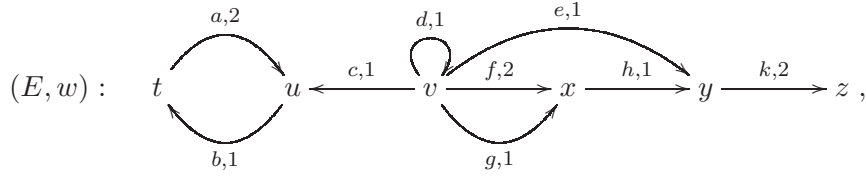
\[
(E,w):\quad\xymatrix@C+15pt{
t\ar@/^1.7pc/[r]^{a,2}&u\ar@/^1.7pc/[l]^{b,1}& v\ar[l]_{c,1}\ar@(ul,ur)^{d,1}\ar@/^1.9pc/[rr]^{e,1}\ar[r]^{f,2}\ar@/_1.7pc/[r]_{g,1}&x\ar[r]^{h,1}&y\ar[r]^{k,2}&z},
\]
which satisfies Condition (LPA). Let $(\tilde E,\tilde w)$ be the weighted graph
 \[
(\tilde E,\tilde w):\quad
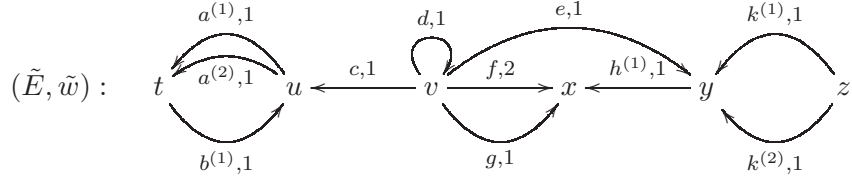
\xymatrix@C+15pt{
t\ar@/_1.7pc/[r]_{b^{(1)},1}&u\ar@/_1.7pc/[l]_{a^{(1)},1}\ar@/_1.0pc/[l]^{a^{(2)},1}& v\ar[l]_{c,1}\ar@(ul,ur)^{d,1}\ar@/^1.9pc/[rr]^{e,1}\ar[r]^{f,2}\ar@/_1.7pc/[r]_{g,1}&x&y\ar[l]_{h^{(1)},1}&z\ar@/_1.7pc/[l]_{k^{(1)},1}\ar@/^1.7pc/[l]^{k^{(2)},1}
}
\]
and $F$ the unweighted graph
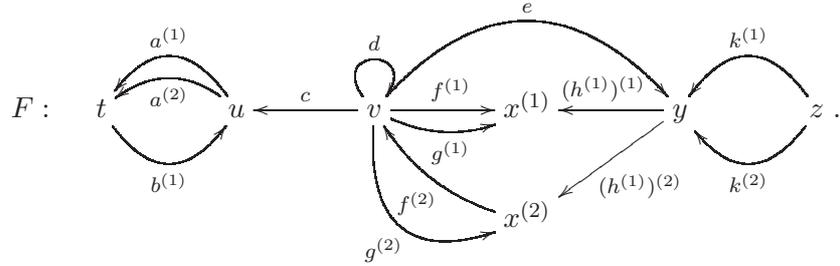
\[
F:\quad\xymatrix@C+15pt{
t\ar@/_1.7pc/[r]_{b^{(1)}}&u\ar@/_1.7pc/[l]_{a^{(1)}}\ar@/_1pc/[l]^{a^{(2)}}& v\ar[l]_{c}\ar@(ul,ur)^{d}\ar@/^2.8pc/[rr]^{e}\ar[r]^{f^{(1)}}\ar@/_0.7pc/[r]_{g^{(1)}}\ar@/_2.5pc/[dr]_{g^{(2)}}&x^{(1)}&y\ar[l]_{(h^{(1)})^{(1)}}\ar[dl]^{(h^{(1)})^{(2)}}&z\ar@/_1.7pc/[l]_{k^{(1)}}\ar@/^1.7pc/[l]^{k^{(2)}}\\
&&&x^{(2)}\ar@/^0.7pc/[ul]^{f^{(2)}}&&
}.
\]
Then $L(E,w)\cong L(\tilde E,\tilde w)$ and $L(\tilde E,\tilde w)\cong L(F)$ by Steps 1 and 2 in the proof of Theorem \ref{thmm}, respectively.
\end{example}

\subsection{Abscence of Condition (LPA)}
%\begin{definition}
%A {\it left-normal d-path} or {\it lenod}-path is a nod-path $p=x_1\dots x_n$ such that $x_1=e_i^*$ for some $e\in E^1,e\neq e^{s(e)}, 1\leq i\leq w(e)$ or $x_1=e_i$ for some $e\in E^1, 2\leq i\leq w(e)$. A {\it right-normal d-path} or {\it rinod}-path is a nod-path $p=x_1\dots x_n$ such that $x_n=e_i$ for some $e\in E^1,e\neq e^{s(e)}, 1\leq i\leq w(e)$ or $x_1=e_i^*$ for some $e\in E^1, 2\leq i\leq w(e)$. A {\it left-right-normal d-path} or {\it lerinod}-path is a nod-path $p$ that is left-normal and right-normal. 
%\end{definition}

%\begin{definition}
%A {\it nod$^2$-path} is a nod-path $p$ such that $p^2$ is again a nod-path.
%\end{definition}

The following lemma plays a crucial role in the proofs of Theorems \ref{thmm3} and \ref{thmm2}.
\begin{lemma}[{\cite[Lemma 17]{Raimund4}}]\label{lemimp}
Suppose that $(E,w)$ does not satisfy Condition (LPA). Then there is a nod-path whose first letter is $e_2$ and whose last letter is $e_2^*$ for some $e\in E^1_w$.
\end{lemma}
\begin{proofs}
We only consider the case that $(E,w)$ does not satisfy (LPA1). So suppose there is a vertex $v\in E^0$ which emits at least two weighted edges. Then we can choose a weighted edge $e\in s^{-1}(v)$ such that $e\neq e^v$. Clearly $e_2e_2^*$ is a nod-path since $e$ is not special.
\end{proofs}

\begin{theorem}[{\cite[Corollary 18 and Theorem 3]{Raimund4}}]\label{thmm3}
If $(E,w)$ does not satisfy Condition (LPA), then $L(E,w)$ is infinite-dimensional as a $K$-vector space, not simple, not graded simple with respect to its standard grading, not locally finite with respect to its standard grading, not Noetherian, not Artinian, not von Neumann regular, and has infinite Gelfand-Kirillov dimension.
\end{theorem}
\begin{proofs}
We only prove that $L(E,w)$ is not Noetherian. By Lemma \ref{lemimp}, we can choose a nod-path $p=x_1\dots x_n$ whose first letter is $e_2$ and whose last letter is $e_2^*$ for some $e\in E^1_w$. Let $q$ be the nod-path one obtains by replacing the first letter of $p$ by $e_1$. For any $n\in\N$ let $I_n$ be the left ideal generated by the nod-paths $p,pq,\dots,pq^n$. One directly checks that $I_n$ equals the linear span of all nod-paths $o$ such that one of the words $p,pq,\dots,pq^n$ is a suffix of $o$. It follows that 
$I_n\subsetneq I_{n+1}$ (clearly none of the words $p,pq,\dots,pq^n$ is a suffix of $pq^{n+1}$ since $p$ and $q$ have the same length but are distinct; hence $pq^{n+1}\not\in I_n$). 
\end{proofs}

\begin{lemma}[{\cite[Lemma 19]{Raimund4}}]\label{lemidem}
Let $p$ be a nod-path starting with $e_2$ and ending with $e_2^*$ for some $e\in E^1_w$. Then the ideal $I$ of $L(E,w)$ generated by $p$ contains no nonzero idempotent.
\end{lemma}

\begin{theorem}[{\cite[Theorem 2]{Raimund4}}]\label{thmm2}
If $(E,w)$ does not satisfy Condition (LPA), then there is no field $K'$ and graph $F$ (row-finite or not) such that $L_K(E,w)\cong L_{K'}(F)$ as rings.
\end{theorem}

\begin{proofs}
Assume there is a field $K'$, a graph $F$ and a ring isomorphism $\phi:L_K(E,w)\rightarrow L_{K'}(F)$. By Lemma \ref{lemimp}, there is a nod-path $p$ whose first letter is $e_2$ and whose last letter is $e_2^*$ for some $e\in E^1_w$. Let $q$ be the nod-path one obtains by replacing the last letter of $p$ by $e_1^*$. By Lemma \ref{lemidem}, the ideal $I$ of $L_K(E,w)$ generated by $p$ contains no nonzero idempotent. Similarly, for any $n\in \N$, the ideal $I_n$ of $L_K(E,w)$ generated by $qp^n$ contains no nonzero idempotent. It follows from \cite[Proposition 2.7.9]{abrams-ara-molina}, that $\phi(I),\phi(I_n)\subseteq I(P_c(F))~(n\in\N)$ where $I(P_c(F))$ is the ideal of $L_{K'}(F)$ generated by all vertices in $F^0$ which belong to a cycle without an exit. Hence $\phi(p),\phi(qp^n)\in I(P_c(F))~(n\in\N)$. By \cite[Theorem 2.7.3]{abrams-ara-molina} we have 
\begin{equation}
I(P_c(F))\cong \bigoplus_{i\in \Gamma}M_{\Lambda_i}(K'[x,x^{-1}])
\end{equation}
as a $K'$-algebra and hence also as a ring. The sets $\Gamma$ and $\Lambda_i~(i\in \Gamma)$ in (9) might be infinite if $F$ is not finite. \\
It follows from the previous paragraph that there is a subring $A$ of $L_K(E,w)$ such that $p,qp^n\in A~(n\in \N)$ and $A\cong \bigoplus_{i\in \Gamma}M_{\Lambda_i}(K'[x,x^{-1}])$. For any $n\in \N$ let $J_n$ be the left ideal of $A$ generated by $qp^2,\dots,qp^{n+1}$. Then $J_n$ is contained in the linear span of all nod-paths $o$ such that one of the words $qp^2,\dots,qp^{n+1}$ is a suffix of $o$. It follows that 
$J_n\subsetneq J_{n+1}$ (clearly none of the words $qp^2,\dots,qp^{n+1}$ is a suffix of $qp^{n+2}$ since $p$ and $q$ have the same length but are distinct). If the sets $\Gamma$ and $\Lambda_i~(i\in \Gamma)$ are finite, then we already have a contradiction since it is well-known that $\bigoplus_{i\in \Gamma}M_{\Lambda_i}(K'[x,x^{-1}])$ is Noetherian in this case. If one of the sets $\Gamma$ and $\Lambda_i~(i\in \Gamma)$ is infinite, then one can use \cite[Lemma 20]{Raimund4} to obtain a contradiction.
\end{proofs}

\section{Finite-dimensionality, Noetherianess and von Neumann regularity}

Throughout this section $(E,w)$ denotes a weighted graph. Suppose we would like to know if $L(E,w)$ is finite-dimensional (respectively Noetherian, von Neumann regular). Then one could first check if $(E,w)$ satisfies Condition (LPA). If not, then $L(E,w)$ is neither finite-dimensional nor Noetherian nor von Neumann regular by Theorem \ref{thmm3}. If $(E,w)$ does satisfy Condition (LPA), then one could apply Steps 1 and 2 in the proof of Theorem \ref{thmm} to obtain an unweighted graph $F$ such that $L(E,w)\cong L(F)$. Having done that one could use the well-known graph-theoretical criteria for finite-dimensionality, Noetherianess and von Neumann regularity of unweighted Leavitt path algebras, see \cite{abrams-ara-molina}. 

In this section we present graph-theoretical criteria which allow us to determine directly, without the detour via the unweighted graph $F$ mentioned above, if $L(E,w)$ is finite-dimensional (respectively Noetherian, von Neumann regular). Consider Conditions (W1) and (W2) below.

\begin{enumerate}[(W1)]
\item No cycle in $(E,w)$ is based at a vertex $v\in T(r(E^1_w))$.
\smallskip
\item There is no $n\geq 1$ and paths $p_1,\dots,p_n,q_1,\dots,q_n$ in $(E,w)$ such that $r(p_i)=r(q_i)~(1\leq i\leq n)$, $s(p_1)=s(q_n)$, $s(p_i)=s(q_{i-1})~(2\leq i \leq n)$ and for any $1\leq i \leq n$, the first letter of $p_i$ is a weighted edge, the first letter of $q_i$ is an unweighted edge and the last letters of $p_i$ and $q_i$ are distinct.
\end{enumerate}
Condition (W2) is based on an idea %which appeared in an unpublished work
by N. T. Phuc. Each of the Conditions (W1) and (W2) above ``forbids" a certain constellation in the weighted graph $(E,w)$. The pictures below illustrate these forbidden constellations.\\
\\
\begin{enumerate}[(W1)]
\item \[\vcenter
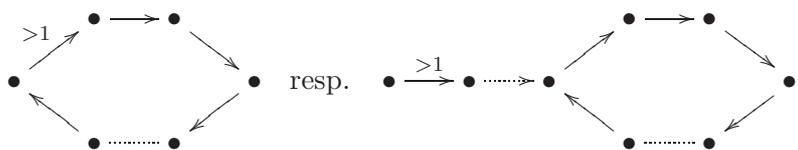
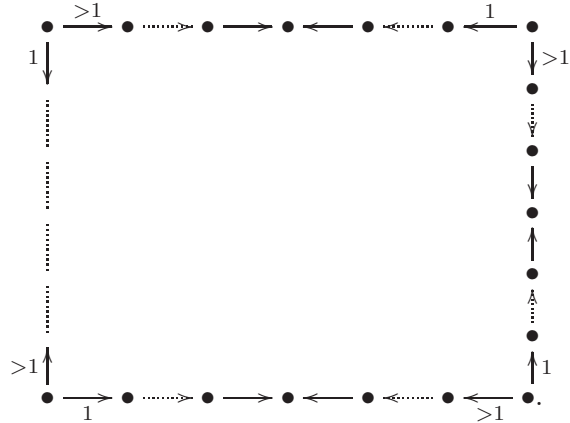
{\vbox{\xymatrix@R-1pc@C-0.5pc{
		&\bullet \ar[r]&\bullet \ar[rd]&\\
		\bullet \ar^{>1}[ru]&&&\bullet \ar[ld]\\
		&\bullet \ar[lu]&\bullet \ar@{.}[l]}}}
		~~~~\text{ resp. }~~~~
\vcenter{\vbox{\xymatrix@R-1pc@C-0.5pc{
		&&&\bullet \ar[r]&\bullet \ar[rd]&\\
		\bullet \ar^{>1}[r] & \bullet \ar@{.>}[r] &\bullet \ar[ru]&&&\bullet. \ar[ld]\\
		&&&\bullet \ar[lu]&\bullet \ar@{.}[l]}}}		
		\]
\medskip
\item \[\xymatrix@R-1pc@C-0.5pc{
\bullet\ar_{1}[d]\ar^{>1}[r]&\bullet\ar@{.>}[r]&\bullet\ar[r]&\bullet&\bullet\ar[l]&\bullet\ar@{.>}[l]&\bullet\ar_{1}[l]\ar^{>1}[d]\\
\ar@{.}[d]&&&&&&\bullet\ar@{.>}[d]\\
\ar@{.}[d]&&&&&&\bullet\ar[d]\\
&&&&&&\bullet\\
\ar@{.}[u]&&&&&&\bullet\ar[u]\\
\ar@{.}[u]&&&&&&\bullet\ar@{.>}[u]\\
\bullet\ar^{>1}[u]\ar_{1}[r]\ar@{.}[u]&\bullet\ar@{.>}[r]&\bullet\ar[r]&\bullet&\bullet\ar[l]&\bullet\ar@{.>}[l]&\bullet.\ar_{1}[u]\ar^{>1}[l]
}\]
\end{enumerate}
$~$\\
$(E,w)$ is called {\it well-behaved} if Conditions (LPA1), (LPA2), (LPA3), (W1) and (W2) are satisfied. Note that Condition (W1) is stronger than Condition (LPA4). Hence a well-behaved weighted graph satisfies Condition (LPA).

\subsection{Finite-dimensional weighted Leavitt path algebras}

\begin{lemma}[{\cite[Corollary 23]{preusser1}}]\label{cor3.0}
Suppose $(E,w)$ is well-behaved. Then $\{c,c^*\mid c \text{ is a cycle in }\hat E\}$ is the set of all quasicycles.  
\end{lemma}

$(E,w)$ is called {\it acyclic} if there is no cycle in $(E,w)$, and {\it aquasicyclic} if there is no quasicycle. 

\begin{theorem}[{\cite[Theorem 47]{preusser}}, {\cite[Theorem 25]{preusser1}}]\label{thmfd}
The following are equivalent:
\begin{enumerate}[(i)]
\item $L(E,w)$ is finite-dimensional.
\item $(E,w)$ is finite, acyclic and well-behaved.
\item $(E,w)$ is finite and aquasicyclic.
\item $(E,w)$ is finite and $\GKdim L(E,w)=0$.
\item $L(E,w)\cong \bigoplus_{i=1}^{m}\Mat_{n_i}(K)$ for some $m,n_1,\dots,n_m\in\N$.
\end{enumerate}
\end{theorem}
\begin{proof}
$(i)\Rightarrow (ii).$  Suppose that $L(E,w)$ is finite-dimensional. It follows from Theorem \ref{thmbasis} that $(E,w)$ is finite and acyclic. One checks easily that if one of the Conditions (LPA1), (LPA2), (LPA3), (W1) and (W2) was not satisfied, then there would be a nod$^2$-path (cf. Lemma \ref{lemimp}) and hence $L(E,w)$ would have infinite dimension, again by Theorem \ref{thmbasis}. Hence $(E,w)$ is well-behaved.

$(ii)\Rightarrow (iii).$ Suppose that $(E,w)$ is finite, acyclic and well-behaved. Then $\hat E$ is acyclic since a cycle in $\hat E$ would lift to a cycle in $(E,w)$. It follows from Lemma \ref{cor3.0} that there is no quasicycle. 

$(iii)\Rightarrow (iv).$ Suppose that $(E,w)$ is finite and aquasicyclic. Then $\GKdim L(E,w)=0$ by Theorem \ref{thmgk}.

$(iv)\Rightarrow (i).$ Follows from the fact that a finitely generated $K$-algebra $A$ is finite-dimensional as a $K$-vector space if and only if $\GKdim A=0$.

$(i)\Rightarrow (v).$ Suppose that $L(E,w)$ is finite-dimensional. Since $(i)\Rightarrow (ii)$, $(E,w)$ is well-behaved and therefore satisfies Condition (LPA). Hence, by Theorem \ref{thmm}, $L(E,w)$ is isomorphic to an unweighted Leavitt path algebra. By \cite[Theorem 2.6.17]{abrams-ara-molina} any finite-dimensional unweighted Leavitt path algebra is isomorphic to a finite direct sum of matrix rings over $K$.

$(v)\Rightarrow (i).$ This implication is obvious.
\end{proof}

\subsection{Noetherian weighted Leavitt path algebras}
A {\it left adhesive nod-path} or {\it lenod-path} is a nod-path $p$ such that the juxtaposition $op$ is a nod-path for any nontrivial nod-path $o$ such that $r(o)=s(p)$. A {\it right adhesive nod-path} or {\it rinod-path} is a nod-path $p$ such that the juxtaposition $po$ is a nod-path for any nontrivial nod-path $o$ such that $s(o)=r(p)$. A {\it left-right adhesive nod-path} or {\it lerinod-path} is a nod-path that is left adhesive and right adhesive. If $A, B$ are nonempty words over the same alphabet, then we write $A\sim B$ if $A$ is a suffix of $B$ or $B$ is a suffix of $A$, and $A\not\sim B$ otherwise. %If $X\subseteq L(E,w)$, then we denote by $LI(X)$ the left ideal of $L(E,w)$ generated by $X$.

\begin{lemma}[{\cite[Lemma 48]{preusser1}}]\label{lem4.1}
If $L(E,w)$ is Noetherian, then there is no lenod-path $p$ and nod$^2$-path $q$ such that $p\not\sim q$ and $pq$ is a nod-path.
\end{lemma}
\begin{proofs}
Assume there is a lenod-path $p$ and nod$^2$-path $q$ such that $p\not\sim q$ and $pq$ is a nod-path. For any $n\in \N$ let $I_n$ be the left ideal of $L(E,w)$ generated by $p, pq,\dots, pq^n$. Then one can show that $I_1 \subsetneq I_2\subsetneq I_3\subsetneq\dots$.
\end{proofs}

\begin{lemma}[{\cite[Lemma 50]{preusser1}}]\label{lem4.2}
If $L(E,w)$ is Noetherian, then there is no nod$^2$-path $p$ based at a vertex $v$ such that $pp^*$ is a nod-path and $p^*p=v$ in $L(E,w)$.
\end{lemma}
\begin{proofs}
Assume there is a nod$^2$-path $p$ based at a vertex $v$ such that $pp^*$ is a nod-path and $p^*p=v$ in $L(E,w)$. For any $n\in \N$ let $I_n$ be the left ideal of $L(E,w)$ generated by $v-p^n(p^*)^n$. Then one can show that $I_1 \subsetneq I_2\subsetneq I_3\subsetneq\dots$.
\end{proofs}

Recall that a graded $K$-algebra $A=\bigoplus_{g\in G} A_g$ is called {\it locally finite} if $\dim_K A_g < \infty$ for every $g\in G$.
\begin{theorem}[{\cite[Theorem 52]{preusser1}}]\label{thmnoet}
The following are equivalent:
\begin{enumerate}[(i)]
\item $L(E,w)$ is Noetherian.
\item $(E,w)$ is finite, well-behaved and no cycle has an exit.
\item $(E,w)$ is finite and $\GKdim L(E,w)\leq 1$.
\item $L(E,w)$ is locally finite with respect to its standard grading.
\item $L(E,w)\cong\big(\bigoplus_{i=1}^{l} \Mat_{m_i}(K)\big)\oplus\big(\bigoplus_{j=1}^{l'} \Mat_{n_j}(K[x,x^{-1}])\big)$ for some integers $l,l'\geq 0$ and $m_i,n_j\geq 1$.
%\item $L(E,w)$ is isomorphic to a locally finite Leavitt path algebra.
\end{enumerate}
\end{theorem}
\begin{proofs}
First we show that  $(ii)\Rightarrow (iii) \Rightarrow (iv)\Rightarrow (ii)$ (and hence $(ii)\Leftrightarrow (iii) \Leftrightarrow (iv)$).

$(ii)\Rightarrow (iii)$. Suppose that $(E,w)$ is finite, well-behaved and no cycle has an exit. Using Lemma \ref{cor3.0} one can show that no quasicycle is selfconnected. Moreover, if $q$ and $q'$ are quasicycles such that $q\not\approx q'$, then $q\Longrightarrow q'$ cannot hold, see \cite[Corollary 32]{preusser1}. Hence $\GKdim L(E,w)\leq 1$ by Theorem \ref{thmgk}.

$(iii)\Rightarrow (vi)$. Suppose that $(E,w)$ is finite and $\GKdim L(E,w)\leq 1$. By Lemma  \ref{lemnewgk} and Theorem \ref{thmgk} any nontrivial nod-path is of the form $o_1p^lqo_2$ where $o_1$ and $o_2$ are either the empty word or nod-paths in $P'$, $p$ is a quasicycle, $l\geq 0$ and $q\neq p$ is a prefix of $p$. By \cite[Remark 16(a)]{preusser} there are only finitely many quasicycles and, by \cite[Lemma 34]{preusser1}, none of them is of homogeneous degree $0$. Further we have $|P'|<\infty$ by \cite[Lemma 21]{preusser}. It follows that $L(E,w)$ is locally finite (for fixed $o_1$, $p$, $q$ and $o_2$ as above there can at most be one nod-path of the form $o_1p^lqo_2$ in each homogeneous component).

$(iv)\Rightarrow (ii)$. Suppose $L(E,w)$ is locally finite. Assume that $(E,w)$ is not finite. Then $|E^0|=\infty$. But $E^0\subseteq L(E,w)_0$ and $E^0$ is a linearly independent set by Theorem \ref{thmbasis}, which contradicts the assumption that $(E,w)$ is locally finite. Hence $(E,w)$ is finite. Assume now that one of the Conditions (LPA1), (LPA2), (LPA3), (W1) and (W2) is not satisfied. Then there is a nod$^2$-path $p$ starting with $e_2$ for some $e\in E_w^1$. Clearly for any $n\in\N$, $(p^*)^np^n$ is a nod-path in the homogeneous $0$-component of $L(E,w)$. But this contradicts the assumption that $(E,w)$ is locally finite. Hence $(E,w)$ is well-behaved.
It remains to show that no cycle has an exit. Assume there is a cycle $c=e^{(1)}\dots e^{(n)}$ in $(E,w)$ with an exit $f\in E^1$. Clearly $e^{(1)},\dots, e^{(n)}\in E^1_{uw}$ because Condition (W1) is satisfied. Without loss of generality assume that $s(f)=s(e^{(n)})$ and $f\neq e^{(n)}$. Since $w(e^{(n)})=1$, we can choose the special edge $e^{s(e^{(n)})}\neq e^{(n)}$. Set $\hat c:=e_1^{(1)}\dots e_1^{(n)}$. Then for any $n\in\N$, $\hat c^n(\hat c^*)^n$ is a nod-path in the homogeneous $0$-component of $L(E,w)$. But this contradicts the assumption that $(E,w)$ is locally finite. Hence no cycle in $(E,w)$ has an exit.

Next we show that $(i)\Rightarrow (ii)$ and $(iii)\Rightarrow (i)$ (and hence $(i)\Leftrightarrow(ii)\Leftrightarrow (iii) \Leftrightarrow (iv)$).

$(i)\Rightarrow (ii).$ Suppose that $L(E,w)$ is Noetherian. It follows from Theorem \ref{thmm3} that $(E,w)$ satisfies Condition (LPA). Moreover, $(E,w)$ is finite (otherwise taking a sequence $(v_n)_{n\geq 1}$ of pairwise distinct vertices one would have $I_1\subsetneq I_2\subsetneq\dots$ where for any $n\in\N$, $I_n$ is the left ideal of $L(E,w)$ generated the set $\{v_1,\dots,v_n\}$). It remains to show that $(E,w)$ satisfies (W1) and (W2), and that no cycle has an exit. Assume that Condition (W1) is not satisfied. Then, since $(E,w)$ satisfies Condition (LPA4), there is a cycle $e^{(1)}\dots e^{(n)}$ where $e^{(1)}\in E^1_w$. Set $p:=e^{(1)}_2e^{(2)}_1\dots e^{(n)}_1$ and $q:=e^{(1)}_1e^{(2)}_1\dots e^{(n)}_1$. Then $p$ is a lenod-path and $q$ is a nod$^2$-path such that $p\not\sim q$ and $pq$ is a nod-path. But this contradicts Lemma \ref{lem4.1}. Hence Condition (W1) is satisfied. Assume that Condition (W2) is not satisfied. Then there is a nod$^2$-path $p$ based at a vertex $v$ such that $p^*p$ is a nod-path and $pp^*=v$ in $L(E,w)$ (see the proof of \cite[Theorem 51]{preusser1}). But that contradicts Lemma \ref{lem4.2}. Hence Condition (W2) is satisfied. Assume now that there is a cycle $c=e^{(1)}\dots e^{(n)}$ in $(E,w)$ with an exit $f\in E^1$. Clearly $e^{(1)},\dots, e^{(n)}\in E^1_{uw}$ because Condition (W1) is satisfied. Without loss of generality assume that $s(f)=s(e^{(n)})$ and $f\neq e^{(n)}$. Clearly we can choose $e^{s(e^{(n)})}\neq e^{(n)}$ since $w(e^{(n)})=1$. Set $v:=s(c)$ and $\hat c:=e^{(1)}_1\dots e^{(n)}_1$. Then $\hat c$ is a nod$^2$-path based at $v$, $\hat c\hat c^*$ is a nod-path and $\hat c^*\hat c=v$. But that contradicts Lemma \ref{lem4.2}. Hence no cycle in $(E,w)$ has an exit. 

$(iii)\Rightarrow (i).$ Suppose that $(E,w)$ is finite and $\GKdim L(E,w)\leq 1$. Then $(E,w)$ is well-behaved and therefore satisfies Condition (LPA), since we already proved that $(iii)\Rightarrow (ii)$. Hence, by Theorem \ref{thmm} and its proof, there is a finite graph $F$ such that $L(E,w)\cong L(F)$ as $K$-algebras. It follows from \cite[Theorem 5]{zel12} that no cycle in $F$ has an exit (otherwise one would have $\GKdim L(F)\geq 2$). Thus, by \cite[Theorem 4.2.17]{abrams-ara-molina}, $L(E,w)$ is Noetherian (note that $L(E,w)$ is left Noetherian iff it is right Noetherian iff it is Noetherian since $L(E,w)$ has an involution).

Finally we show that $(i)\Leftrightarrow (v)$ (and thus $(i)\Leftrightarrow(ii)\Leftrightarrow (iii) \Leftrightarrow (iv)\Leftrightarrow (v)$).

$(i)\Rightarrow (v).$ Suppose that $L(E,w)$ is Noetherian. Then $(E,w)$ is well-behaved and therefore satisfies Condition (LPA), since we already proved that $(i)\Rightarrow (ii)$. Hence, by Theorem \ref{thmm}, there is a graph $F$ such that $L(E,w)\cong L(F)$ as $K$-algebras. Thus, by \cite[Theorem 4.2.17]{abrams-ara-molina}, $L(E,w)$ is isomorphic to $\big(\bigoplus_{i=1}^{l} \Mat_{m_i}(K)\big)\oplus\big(\bigoplus_{j=1}^{l'} \Mat_{n_j}(K[x,x^{-1}])\big)$ for some integers $l,l'\geq 0$ and $m_i,n_j\geq 1$. 

$(v)\Rightarrow (i).$
Follows from the fact that finite matrix rings over $K$ or $K[x,x^{-1}]$ are Noetherian.

\end{proofs}

\subsection{Von Neumann regular weighted Leavitt path algebras}
\begin{lemma}\label{lemNeu}
Suppose that $(E,w)$ satisfies Condition (LPA). If $(E,w)$ contains a cycle or is not well-behaved, then $L(E,w)\cong L(F)$ for some unweighted graph $F$ containing a cycle.
\end{lemma}
\begin{proof}
Assume $(E,w)$ contains a cycle $c$. Since $(E,w)$ satisfies Condition (LPA), we can apply Steps 1 and 2 in the proof of Theorem \ref{thmm} to obtain a graph $F$ such that $L(E,w)\cong L(F)$ as $K$-algebras. The weighted graph $(\tilde E, \tilde w)$ one obtains by applying Step 1 to $(E,w)$ clearly contains a cycle $\tilde c=\tilde e^{(1)}\dots\tilde e^{(n)}$ (since either all the edges of $c$ get reversed or none of them). Since the ranges of weighted edges in $(\tilde E, \tilde w)$ are sinks, we have $\tilde r(\tilde e^{(i)})\not\in \tilde r(\tilde E^1_w)$ for any $1\leq i\leq n$ (and hence $\tilde e^{(i)}\in \tilde E^1_{uw}$ for any $1\leq i\leq n$). Hence $\tilde c=\tilde e^{(1)}\dots\tilde e^{(n)}$ is also a cycle in the graph $F$ one gets by applying Step 2 to $(\tilde E,\tilde w)$. 

Assume now that $(E,w)$ is acyclic but not well-behaved. Then $(E,w)$ does not satisfy Condition (W2). Hence there are paths $p_i=e^{i,1}\dots e^{i,l_i},q_i=f^{i,1}\dots f^{i,m_i}~(1\leq i\leq n)$ such that $r(p_i)=r(q_i)~(1\leq i\leq n)$, $s(p_1)=s(q_n)$, $s(p_i)=s(q_{i-1})~(2\leq i \leq n)$ and for any $1\leq i \leq n$, $e^{i,1}$ is a weighted edge, $f^{i,1}$ is an unweighted edge and $e^{i,l_i}\neq f^{i,m_i}$. Since $(E,w)$ satisfies Condition (LPA), we can apply Steps 1 and 2 in the proof of Theorem \ref{thmm} to obtain a graph $F$ such that $L(E,w)\cong L(F)$ as $K$-algebras. Set $Z:=T(r(E^1_w))$. We need the claim below.

\underline{Claim} Let $1\leq i\leq n$, $1\leq j\leq l_i$ and $1\leq k\leq m_i$. Then $s(e^{i,j})\in Z\Leftrightarrow j>1$, and $s(f^{i,k})\not\in Z$.

\underline{Proof} Clearly $s(e^{i,j})\in Z$ for any $j>1$ because $e^{i,1}$ is a weighted edge. Since $s(e^{i,1})=s(f^{n,1})$ if $i=1$ respectively $s(e^{i,1})=s(f^{i-1,1})$ if $i>1$, it only remains to show that $s(f^{i,k})\not\in Z$. Assume $s(f^{i,k})\in Z$. We will show that this assumption leads to a contradiction. Since $Z$ is hereditary, it follows that $s(f^{i,m_i})\in Z$. If $l_i>1$, then $s(e^{i,l_i})\in Z$. But this contradicts \cite[Lemma 8]{Raimund4} (since $e^{i,l_i}$ and $f^{i,m_i}$ are distinct, $s(e^{i,l_i}),s(f^{i,m_i})\in Z$ and $r(e^{i,l_i})=r(f^{i,m_i})$). Hence $l_i=1$. Since $s(f^{i,m_i})\in Z$, there is an $h\in E_w^1$ and a path $p$ from $r(h)$ to $s(f^{i,m_i})$. Since $(E,w)$ satisfies Condition (LPA3), $h$ and $e^{i,1}$ are in line, i.e. $h=e^{i,1}$ or there is a path from $r(h)$ to $s(e^{i,1})$ or there is a path from $r(e^{i,1})$ to $s(h)$. Assume that $h=e^{i,1}$. Then $pf^{i,m_i}$ is a closed path at $r(h)=r(e^{i,1})=r(f^{i,m_i})$. But that is not possible since $(E,w)$ is acyclic. Assume now that there is a path from $r(h)$ to $s(e^{i,1})$. Then $s(e^{i,1})\in Z$, which contradicts \cite[Lemma 8]{Raimund4} (see above). Finally assume that there is a path $q$ from $r(e^{i,1})$ to $s(h)$. Then $qhpf^{i,m_i}$ is a closed path at $r(e^{i,1})=r(f^{i,m_i})$. But that is not possible since $(E,w)$ is acyclic. Hence the assumption $s(f^{i,k})\in Z$ leads to a contradiction. Thus $s(f^{i,k})\not\in Z$, which finishes the proof of the claim.\qed

It follows from the claim above that if one applies Step 1 in the proof of Theorem \ref{thmm} to $(E,w)$, then out of the edges $e^{i,j},f^{i,k}~(1\leq i\leq n,~1\leq j\leq l_i,~1\leq k\leq m_i)$ only the edges $e^{i,j}$ where $j>1$ get reversed (to be more precise, the edges $e^{i,j}$ where $j>1$ are replaced by one ore more unweighted edges with reversed orientation). Hence there are $\tilde e^{(1)},\dots,\tilde e^{(n)}\in \tilde E^1_w$ and paths $\tilde q_i=\tilde f^{i,1}\dots \tilde f^{i,t_i}~(1\leq i\leq n)$ in the resulting weighted graph $(\tilde E,\tilde w)$ such that $\tilde r(\tilde e^{(1)})=\tilde r(\tilde q_i)~(1\leq i\leq n)$, $\tilde s(\tilde e^{(1)})=\tilde s(\tilde q_n)$, $\tilde s(\tilde e^{(i)})=\tilde s(\tilde q_{i-1})~(2\leq i \leq n)$, $\tilde f^{i,1}\in\tilde E^1_{uw}~(1\leq i\leq n)$ and $\tilde e^{(i)}\neq \tilde f^{i,t_i}~(1\leq i\leq n)$. Since the ranges of weighted edges in $(\tilde E,\tilde w)$ are sinks, we have $\tilde f^{i,k}\in \tilde E_{uw}^1$ for any $1\leq i\leq n$ and $1\leq k<t_i$. Since no vertex in $(\tilde E,\tilde w)$ receives two distinct weighted edges, we also have $\tilde f^{i,t_i}\in \tilde E_{uw}^1$ for any $1\leq i\leq n$. Hence all the $\tilde f^{i,k}$ are unweighted. The picture below illustrates the situation (the weights are omitted).

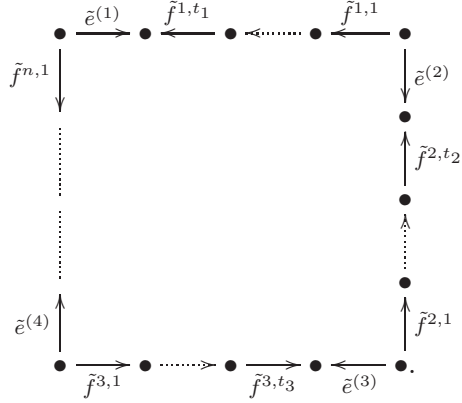
\[\xymatrix@R=0.7cm@C=0.7cm{
\bullet\ar_{\tilde f^{n,1}}[d]\ar^{\tilde e^{(1)}}[r]&\bullet&\bullet\ar[l]_{\tilde f^{1,t_1}}&\bullet\ar@{.>}[l]&\bullet\ar_{\tilde f^{1,1}}[l]\ar^{\tilde e^{(2)}}[d]\\
&&&&\bullet\\
\ar@{.}[u]&&&&\bullet\ar[u]_{\tilde f^{2,t_2}}\\
\ar@{.}[u]&&&&\bullet\ar@{.>}[u]\\
\bullet\ar^{\tilde e^{(4)}}[u]\ar_{\tilde f^{3,1}}[r]\ar@{.}[u]&\bullet\ar@{.>}[r]&\bullet\ar[r]_{\tilde f^{3,t_3}}&\bullet&\bullet.\ar_{\tilde f^{2,1}}[u]\ar^{\tilde e^{(3)}}[l]
}\]

Now we apply Step 2 to $(\tilde E,\tilde w)$. That has the following effect on the weighted subgraph $S$ of $(\tilde E,\tilde w)$ displayed above. Let $1\leq i\leq n$ and set $u_i:=s(\tilde e^{(i)})$, $v_i:=\tilde r(\tilde e^{(i)})$ and $n_i:=\tilde w(\tilde e^{(i)})$. Then Step 2 replaces $v_i$ by $n_i$ vertices $v_i^{(1)},\dots,v_n^{(n_i)}$. The edge $\tilde e^{(i)}$ is replaced by an unweighted edge $(\tilde e^{(i)})^{(1)}$ from $u_i$ to $v_i^{(1)}$ and unweighted edges $ (\tilde e^{(i)})^{(j)}~(2\leq j\leq n_i)$ from $v_i^{(j)}$ to $u_i$. Moreover, the edge $\tilde f^{i,t_i}$ is replaced by unweighted edges $(\tilde f^{i,t_i})^{(j)}~(1\leq j\leq n_i)$ from $\tilde s(f^{i,t_i})$ to $v_i^{(j)}$. The other vertices and edges of $S$ are not changed (since ranges of weighted edges in $(\tilde E,\tilde w)$ are sinks). The picture below illustrates the situation.

\[\xymatrix@R=2cm@C=2cm{&\bullet\ar_{(\tilde e^{(1)})^{(n_1)}}[dl]\ar@{.}[d]&&&&\\
\bullet\ar_{\tilde f^{n,1}}[d]\ar_(.6){(\tilde e^{(1)})^{(1)}}[dr]&\bullet\ar_{(\tilde e^{(1)})^{(2)}}[l]&\bullet\ar[ld]^(.6){(\tilde f^{1,t_1})^{(1)}}\ar[l]_{(\tilde f^{1,t_1})^{(2)}}\ar[lu]_{(\tilde f^{1,t_1})^{(n_1)}}&\bullet\ar@{.>}[l]&\bullet\ar_{\tilde f^{1,1}}[l]\ar_{(\tilde e^{(2)})^{(1)}}[dl]&\\
&\bullet&&\bullet &\bullet \ar_(.4){(\tilde e^{(2)})^{(2)}}[u]&\bullet\ar_{(\tilde e^{(2)})^{(n_2)}}[ul]\ar@{.}[l]\\
\ar@{.}[u]&&&&\bullet\ar[ul]^{(\tilde f^{2,t_2})^{(1)}}\ar[u]_(.6){(\tilde f^{2,t_2})^{(2)}}\ar[ur]_(.5){(\tilde f^{2,t_2})^{(n_2)}}&\\
\ar@{.}[u]&&&\bullet&\bullet\ar@{.>}[u]&\\
\bullet\ar@{.}[u]\ar_{\tilde f^{3,1}}[r]\ar@{.}[u]&\bullet\ar@{.>}[r]&\bullet\ar[ur]^(.5){(\tilde f^{3,t_3})^{(1)}}\ar[r]_(.5){(\tilde f^{3,t_3})^{(2)}}\ar[dr]_(.5){(\tilde f^{3,t_3})^{(n_3)}}&\bullet\ar_{(\tilde e^{(3)})^{(2)}}[r]&\bullet.\ar_{\tilde f^{2,1}}[u]\ar_{(\tilde e^{(3)})^{(1)}}[ul]&\\
&&&\bullet\ar@{.}[u]\ar_{(\tilde e^{(3)})^{(n_3)}}[ur]&&
}\]
One checks easily that
\[\Big(\tilde f^{n,1}\dots\tilde f^{n,t_n-1}(\tilde f^{n,t_n})^{(2)}(\tilde e^{(n)})^{(2)}\Big)\Big(\tilde f^{n-1,1}\dots(\tilde f^{n-1,t_{n-1}})^{(2)}(\tilde e^{(n-1)})^{(2)}\Big)\dots \Big(\tilde f^{1,1}\dots(\tilde f^{1,t_1})^{(2)}(\tilde e^{(1)})^{(2)}\Big)\]
is a closed path in the graph $F$ one gets by applying Step 2 to $(\tilde E,\tilde w)$. Thus $F$ contains a cycle.
\end{proof}

Recall that a ring $R$ is called {\it von Neumann regular} if for any $x\in R$ there is a $y\in R$ such that $xyx=x$. A $K$-algebra is called {\it matricial} if it is isomorphic to a finite direct sum of full finite-dimensional matrix algebras over $K$, and {\it locally matricial} if it is a direct limit of matricial $K$-algebras.

\begin{theorem}\label{thmNeu}
The following are equivalent:
\begin{enumerate}[(i)]
\item $L(E,w)$ is von Neumann regular.
\item $(E,w)$ is acyclic and well-behaved.
\item $L(E,w)$ is locally matricial.
\end{enumerate}
\end{theorem}
\begin{proof}
$(i)\Rightarrow (ii)$. Suppose $L(E,w)$ is von Neumann regular. Then $(E,w)$ satisfies Condition (LPA) by Theorem \ref{thmm3} (if $(E,w)$ did not satisfy Condition (LPA), then there would be a lerinod-path $p$ by Lemma \ref{lemimp}; but for such a $p$ there is no $a\in L(E,w)$ such that $pap=p$). Assume that $(E,w)$ contains a cycle or is not well-behaved. Then, by Lemma \ref{lemNeu}, $L(E,w)\cong L(F)$ for some unweighted graph $F$ containing a cycle. But this contradicts \cite[Theorem 3.4.1]{abrams-ara-molina}. Thus $(E,w)$ is acyclic and well-behaved.

$(ii)\Rightarrow (iii)$. Suppose that $(E,w)$ is acyclic and well-behaved. By Lemma \ref{lemdirlim} we have $L(E,w)=\varinjlim_i L(E_i,w_i)$ where $\{(E_i,w_i)\}$ is the direct system of all finite complete weighted subgraphs of $(E,w)$. It is easy to see that the subgraphs $(E_i,w_i)$ are also acyclic and well-behaved. By Theorem \ref{thmfd} each $L(E_i,w_i)$ is matricial. Thus $L(E,w)$ is locally matricial.

$(iii)\Rightarrow (i)$. Suppose that $L(E,w)$ is locally matricial. It is well known that every matricial $K$-algebra is von Neumann regular and hence so is any direct union of such algebras. Therefore $L(E,w)$ is von Neumann regular.
\end{proof}

 \section{Realisation as generalised corner skew Laurent polynomial rings}

\subsection{Generalised corner skew Laurent polynomial rings}\label{gecorner}

Let $R$ be a unital ring and $p$ an idempotent in $R$. Let $\phi:R\rightarrow pRp$ be a \emph{corner isomorphism}, \index{corner isomorphism} i.e. a ring isomorphism with $\phi(1)=p$. A \emph{corner skew Laurent polynomial ring}  \index{corner skew Laurent polynomial ring} with coefficients in $R$, denoted by $R[t_{+},t_{-},\phi]$, is a unital ring which is constructed as follows:  The elements of $R[t_{+},t_{-},\phi]$ are the formal expressions
\[t^j_{-}r_{-j} +t^{j-1}_{-}r_{-j+1}+\dots+t_{-}r_{-1}+r_0 +r_1t_{+}+\dots +r_it^i_{+},\]
where $r_{-n} \in p_n R$ and $r_n \in R p_n$ for any $n\geq 0$, where $p_0 =1$ and $p_n =\phi^n(p_0)$. The addition is component-wise, and the multiplication is determined by the distributive law and the following rules:
\begin{equation*}
t_{-}t_{+} =1, \qquad t_{+}t_{-} =p, \qquad rt_{-} =t_{-}\phi(r),\qquad  t_{+}r=\phi(r)t_{+}.
\end{equation*}

The corner skew Laurent polynomial rings are studied in~\cite{arabrucom}, where their $K_1$-groups are calculated. This construction is a special case of the so-called fractional skew monoid rings constructed in~\cite{arafrac}. Assigning $1$ to $t_{+}$ and $-1$ to $t_{-}$ makes $A:=R[t_{+},t_{-},\phi]$ a $\mathbb Z$-graded ring with $A=\bigoplus_{i\in \mathbb Z}A_i$, where 
\begin{align*}
A_i& = Rp_it^i_{+}, \text{ for  } i>0,\\
A_i&=t^{-i}_{-}p_{-i}R, \text{ for } i<0,\\
A_0& =R,
\end{align*}
see~\cite[Proposition~1.6]{arafrac}. Clearly, when $p=1$ and $\phi$ is the identity map, then $R[t_{+},t_{-},\phi]$ reduces to the familiar ring $R[t,t^{-1}]$. 
As it was shown in ~\cite[Example~2.5]{arafrac} (see also~\cite[Example~1.6.14]{hazrat16}), Leavitt path algebras of finite graphs with no sources are examples of corner skew Laurent polynomial rings.  This was used to prove that Leavitt path algebras are graded von Neumann regular rings~\cite[Corollary~1.6.17]{hazrat16}. 

In this section we introduce the notion of a generalised corner skew Laurent polynomial ring. A corner skew Laurent polynomial ring is characterised as a $\mathbb Z$-graded unital ring $A=\bigoplus_{i\in \mathbb Z}A_i$ with an $s\in A_{-1}$ and a $t\in A_1$ such that $st=1$. A generalised corner skew Laurent polynomial ring is characterised as a $\mathbb Z^n$-graded unital ring $A=\bigoplus_{i\in \mathbb Z^n}A_i$ with $s_1,\dots,s_n \in \Mat_{1\times m}(A_{\alpha_i})$ and $t_1,\dots,t_n\in \Mat_{m\times 1}(A_{-\alpha_i})$, where $m\in \mathbb N$, such that 
$s_1t_1=\dots =s_nt_n=1$. Here $\{\alpha_i \mid 1\leq i \leq n\}$ is the standard basis of $\mathbb Z^n$.

\begin{definition}\label{def15}
Let $G$ be a group, written multiplicatively, and $A=\{A_g\ | \ g\in G\}$ a family of abelian groups, written additively. Furthermore, let for any $g,h\in G$,
\begin{align*}\phi_{g,h}:A_g\times A_h&\rightarrow A_{gh}
\end{align*}
be a biadditive map. If $a\in A_g$ and $b\in A_h$, we sometimes write $a.b$ instead of $\phi_{g,h}(a,b)$. Suppose that 
\begin{enumerate}[(i)]
\item there is an element $1\in A_0$ such that $a.1=a=1.a$ for any $g\in G$ and $a\in A_g$ and 
\medskip 
\item $(a.b).c=a.(b.c)$ for any $g,h,k\in G$ and $a\in A_g,b\in A_h,c\in A_k$.
\end{enumerate}
Define the map
\begin{align*}
\cdot:  \bigoplus_{g\in G}A_g\times \bigoplus_{g\in G}A_g&\longrightarrow \bigoplus_{g\in G}A_g\\
\big((a_g)_{g\in G},(b_g)_{g\in G}\big)&\longmapsto (c_g)_{g\in G}
\end{align*}
where 
\[c_g=\sum_{hk=g}a_h.b_k.\]
Then $(\bigoplus_{g\in G}A_g,+,\cdot)$ is a ring which we denote by $A[G]$ and call a {\it generalised group ring}. 
\end{definition}

\begin{remark}
\begin{enumerate}[(a)]
\item If each $A_g=R$ where $R$ is a unital ring and each $\phi_{g,h}$ is the multiplication in $R$, then $A[G]$ is the group ring $R[G]$.
\medskip 
\item $A[G]$ is a $G$-graded ring such that its $g$-component equals $A_g$ for any $g\in G$ (we identify $A_g$ with its image in $\bigoplus_{g\in G}A_g$).
\end{enumerate}
\end{remark}

Until the end of this subsection $n$ denotes a fixed positive integer. We define a map $~\hat{}:\Z^n\to \N_0^n$ by $\hat g=(\max(g_1,0),\dots, \max(g_n,0))$. Recall that $\{\alpha_i\mid 1\leq i \leq n\}$ denotes the standard basis of $\mathbb Z^n$.

\begin{definition}\label{def17}
Let $A[\Z^n]=\bigoplus_{g\in \mathbb Z^n}A_g$ be a generalised group ring such that the conditions (i)-(iv) below are satisfied. 
\begin{enumerate}[(i)]
\item There is an $m\in\N$, a unital ring $R$ and idempotents $p_g\in \Mat_{m^{g_1+\dots+g_n}\times m^{g_1+\dots+g_n}}(R)\hspace{0.1cm}(g\in \N_0^n)$
where $p_0=1$ such that $A_g=p_{\widehat g} \Mat_{m^{\widehat g_1+\dots+\widehat g_n}\times m^{(\widehat{-g})_1+\dots+(\widehat{-g})_n}}(R)p_{\widehat{-g}}$ for any $g\in \Z^n$.
\medskip 
\item $\sigma.\tau=\sigma\tau$ for any $1\leq i\leq n$, $\sigma\in A_{-\alpha_{i}}$ and $\tau\in A_{\alpha_i}$, where $\sigma\tau$ is the usual matrix product of $\sigma$ and $\tau$.
\vspace{-0.25cm}
\item
$1\in A_{\alpha_i}.A_{-\alpha_{i}}$ for any $1\leq i\leq n$. Here $A_{\alpha_i}.A_{-\alpha_{i}}=\{\sum_{j=1}^ka_j.b_j\mid k\in\N, a_j\in A_{\alpha_i}, b_j\in A_{-\alpha_i}\}$.
\medskip 
\item
$\sigma.\tau.\omega=\sigma\tau\omega$ for any $g\in \Z^n$, $\sigma,\omega\in A_g$ and $\tau\in A_{-g}$, where $\sigma\tau\omega$ is the usual matrix product of $\sigma$, $\tau$ and $\omega$.
%\item if $p_g=e_g$ for any $g\in (\N_0)^n$ where $e_g$ denotes the identity matrix in  $M_{m^{g_1+\dots+g_n}}$ $_{\times m^{g_1+\dots+g_n}}(R)$, then each $A_g$ contains an (in $\R(\Z^n)$) invertible element, 
\end{enumerate}
Then $A[\Z^n]$ is called a {\it generalised corner skew Laurent polynomial ring}. 
\end{definition}

\begin{remark}\begin{enumerate}[(a)]
\item We usually denote the ring $A[\Z^n]$ by $R[t_1,\dots t_n,t_{-n},\dots,t_{-1},\phi]$ and write an element $(\sigma_g)_{g\in \Z^n}\in A[\Z^n]$ in the form $\sum_{g\in \Z^n}t_1^{\widehat{g}_1}\dots t_n^{\widehat{g}_n}\sigma_gt_{-n}^{(\widehat{-g})_{n}}\dots t_{-1}^{(\widehat{-g})_{1}}$.
\medskip
\item  Let $R$ be a unital ring and $p$ an idempotent of $R$. Let $\phi: R\rightarrow pRp$ be a corner isomorphism, i.e. a ring isomorphism with $\phi(1) = p$. Set $p_g:=\phi^g(1)$ for any $g\in \N_0$, $A_g:=p_{\widehat g}Rp_{\widehat{-g}}$ for any $g\in \Z$, and $A:=\{A_g\ | \ g\in \Z\}$. Define for any $g,h\in\Z$ the map
$\phi_{g,h}:A_g\times A_h\longrightarrow A_{g+h}$ by 
\begin{displaymath}
\phi_{g,h}(a,b) 
=\left\{
\begin{array}{ll}
\phi^h(a)b             & \text{ if }g,h\geq 0\\
\phi^h(ab)             & \text{ if }g\geq 0,h\leq 0,g+h\geq 0\\
\phi^{-g}(ab)         &\text{ if }g\geq 0,h\leq 0,g+h\leq 0\\
\phi^{g+h}(a)b      &\text{ if }g\leq 0,h\geq 0,g+h\geq 0\\
a \phi^{-(g+h)}(b) &\text{ if }g\leq 0,h\geq 0,g+h\leq 0\\
a\phi^{-g}(b)        &\text{ if }g,h\leq 0.
\end{array}\right.
\end{displaymath}

One checks easily that the conditions (i) and (ii) in Definition~\ref{def15} are satisfied. Furthermore, the generalised group ring $A[\Z]$ satisfies the conditions (i)-(iv) in Definition~\ref{def17}. Hence $A[\Z]=R[t_1,t_{-1},\phi]$ is a generalised corner skew Laurent polynomial ring. Let $R[t_+,t_-,\phi]$ be the corner skew Laurent polynomial ring defined by $R$ and $\phi$ (see ~\cite[\S 1.6.2]{hazrat16}). Define the map 
\begin{align*}
\psi:R[t_1,t_{-1},\phi]&\longrightarrow R[t_+,t_-,\phi]\\
\sum_{g\in \Z}t_1^{\widehat{g}_1}\sigma_gt_{-1}^{(\widehat{-g})_{1}}&\longmapsto \sum_{g\in \Z}t_-^{\widehat{g}_1}\sigma_gt_+^{(\widehat{-g})_{1}}.
\end{align*}
It is easy to show, that $\psi$ is a ring isomorphism. Defining a grading on $R[t_+,t_-,\phi]$ by assigning $1$ to $t_-$ and $-1$ to $t_+$, we get $\psi((R[t_1,t_{-1},\phi])_g)=(R[t_+,t_-,\phi])_g$ for any $g\in\Z$, i.e. $\psi$ is a graded ring isomorphism. Hence corner skew Laurent polynomial rings are special cases of generalised corner skew Laurent polynomial rings. 
\end{enumerate}
\end{remark} 

Recall that a $G$-graded ring $A=\bigoplus_{g\in G} A_g$ is called {\it strongly graded} if $A_gA_h=A_{gh}$ for any $g, h\in G$. Here $A_gA_h=\{\sum_{j=1}^k a_jb_j\mid k\in\N, a_1,\dots,a_k\in A_g, b_1,\dots,b_k\in A_h\}$.

\begin{proposition} \label{prop17}
Let $A=R [t_1,\dots t_n,t_{-n},\dots,t_{-1},\phi]$ denote a generalised corner skew Laurent polynomial ring. Then $A$ is strongly graded if and only if $\Mat_{1\times m}(R)  p_{\alpha_i}\Mat_{m\times 1}(R)=R$ for any $1\leq i\leq n$ (where $\Mat_{1\times m}(R) p_{\alpha_i}\Mat_{m\times 1}(R)=\big \{\sum_{j=1}^k u_jp_{\alpha_i}v_j\ | \ k\in\N, u_1,\dots,u_k\in \Mat_{1\times m}(R), v_1,\dots,v_k\in \Mat_{m\times 1}(R)\big \}$).
\end{proposition}

\begin{proof} 
($\Rightarrow$) Suppose $A$ is strongly graded. Let $1\leq i\leq n$. Then $1\in A_{-\alpha_{i}}A_{\alpha_{i}}$. Hence there is a $k\in \N$, $u_1,\dots,u_k\in \Mat_{1\times m}(R)$ and $v_1,\dots,v_k\in \Mat_{m\times 1}(R)$ such that $\sum_{j=1}^k(u_jp_{\alpha_i}t_{-i})\cdot(t_ip_{\alpha_i}v_j)=1$. But $(u_jp_{\alpha_i}t_{-i})\cdot(t_ip_{\alpha_i}v_j)=u_jp_{\alpha_i}v_j$ by Definition~\ref{def17}(ii). Hence $\sum_{j=1}^k u_jp_{\alpha_i}v_j=1$ and therefore $\Mat_{1\times m}(R)p_{\alpha_i}\Mat_{m\times 1}(R)=R$.

($\Leftarrow$)
Suppose $\Mat_{1\times m}(R)p_{\alpha_i}\Mat_{m\times 1}(R)=R$ for any $1\leq i\leq n$. Since $\Z^n$ is generated by the $\alpha_i$'s, it suffices to show that $1\in A_{\alpha_i}A_{-\alpha_{i}}$ and $1\in A_{-\alpha_{i}}A_{\alpha_i}$ for any $1\leq i\leq n$ in order to prove that $A$ is strongly graded (see~\cite[\S 1.1.3]{hazrat16}). Let $1\leq i\leq n$. By Definition~\ref{def17}(iii), $1\in A_{\alpha_i}A_{-\alpha_{i}}$. Let $k\in\N$, $u_1,\dots,u_k\in \Mat_{1\times m}(R)$ and $v_1,\dots,v_k\in \Mat_{m\times 1}(R)$ such that $\sum_{j=1}^k u_jp_{\alpha_i}v_j=1$. Then $\sum_{j=1}^k(u_jp_{\alpha_i}t_{-i})\cdot(t_ip_{\alpha_i}v_j)=1$ by Definition~\ref{def17}(ii) and thus $1\in A_{-\alpha_{i}}A_{\alpha_i}$.
\end{proof}

Recall that a ring $A$ is called {\it von Neumann regular} if for
any $a\in A$ there is a $b\in A$ such that $a = aba$. A graded ring $A$ is called {\it graded von Neumann regular} if for
any homogeneous $a\in A$ there is a $b\in A$ such that $a = aba$. Proposition~\ref{prop22} determines when a generalised corner skew Laurent polynomial ring is graded von Neumann regular. We need the following Definition~\ref{def20} and Lemma~\ref{lem21} for the proof of Proposition~\ref{prop22}.

\begin{definition}\label{def20}
Let $C$ be a category. Then $C$ is called {\it von Neumann regular} if for any morphism $f\in \Hom(X,Y)$ there is a $g\in \Hom(Y,X)$ such that $fgf=f$.
\end{definition}

\begin{lemma} \label{lem21}
Let $R$ denote a unital ring and $C$ the category of finitely generated free right $R$-modules. Then $C$ is von Neumann regular if and only if $R$ is von Neumann regular.
\end{lemma}

\begin{proof}
($\Rightarrow$) 
This direction is obvious since $(\Hom(R,R),\circ)\cong (R,\cdot)$ as a monoid.

($\Leftarrow$) 
This follows from~\cite[Theorem 24, p.114]{kaplansky72} and the following easy to prove remark: If $A$ is a von Neumann regular ring and $e,e'\in A$ are idempotents, then for any $f\in eAe'$ there is a $g\in e'Ae$ such that $fgf=f$.
\end{proof}

\begin{proposition} \label{prop22}
Let $A=R[t_1,\dots t_n,t_{-n},\dots,t_{-1},\phi]$ denote a generalised corner skew Laurent polynomial ring. Then $A$ is graded von Neumann regular if and only if $R$ is von Neumann regular.
\end{proposition}

\begin{proof}
($\Rightarrow$)
If a graded ring is graded von Neumann regular, then it is easy to see
that its zero component ring is von Neumann regular.

($\Leftarrow$)
Suppose $R$ is von Neumann regular. Let $g\in \Z^n$ and $\sigma\in A_g\subseteq \Mat_{m^{\widehat g_1+\dots+\widehat g_n}\times}$ $_{m^{(\widehat{-g})_1+\dots+(\widehat{-g})_n}}(R)$. By Lemma~\ref{lem21}, there is a $\tau'\in \Mat_{m^{(\widehat{-g})_1+\dots+(\widehat{-g})_n}\times m^{\widehat g_1+\dots+\widehat g_n}}(R)$ such that $\sigma\tau'\sigma=\sigma$. Set $\tau:=p_{\widehat{-g}}\tau'p_{\widehat g}\in A_{-g}$. Then clearly $\sigma\tau\sigma=\sigma\tau'\sigma=\sigma$ since $p_{\widehat g}$ and $p_{\widehat{-g}}$ are idempotents. By Definition~\ref{def17}(iv), 
\begin{align*}&(t_1^{\widehat{g}_1}\dots t_n^{\widehat{g}_n}\sigma t_{-n}^{(\widehat{-g})_{n}}\dots t_{-1}^{(\widehat{-g})_{1}})(t_1^{(\widehat{-g})_1}\dots t_n^{(\widehat{-g})_n}\tau t_{-n}^{\widehat{g}_{n}}\dots t_{-1}^{\widehat{g}_{1}})(t_1^{\widehat{g}_1}\dots t_n^{\widehat{g}_n}\sigma t_{-n}^{(\widehat{-g})_{n}}\dots t_{-1}^{(\widehat{-g})_{1}})\\
=&~t_1^{\widehat{g}_1}\dots t_n^{\widehat{g}_n}\sigma\tau\sigma t_{-n}^{(\widehat{-g})_{n}}\dots t_{-1}^{(\widehat{-g})_{1}}\\
=&~t_1^{\widehat{g}_1}\dots t_n^{\widehat{g}_n}\sigma t_{-n}^{(\widehat{-g})_{n}}\dots t_{-1}^{(\widehat{-g})_{1}}.
\qedhere\end{align*}
\end{proof}

Let $R$ be a unital ring and $m\in\mathbb{N}$. Below we define a multiplication on the set $\bigcup_{i,j\in \N_0} \Mat_{m^i\times m^j}(R)$ which extends the usual matrix multiplication and makes $\bigcup_{i,j\in \N_0} \Mat_{m^i\times m^j}(R)$ a monoid. This monoid structure will be used in the proof of Theorem~\ref{thm24}.

\begin{definition}\label{def23}
Let $R$ be a unital ring and $m\in\mathbb{N}$. Set $\Mf=\Mf(m,R):=\bigcup_{i,j\in \N_0} \Mat_{m^i\times m^j}(R)$. We define a multiplication $\cdot$ on $\Mf$ as follows. Let $A,B\in \Mf$. Then there are $i,j,k,l\in\N_0$ such that $A\in \Mat_{m^i\times m^j}(R)$ and $B\in \Mat_{m^k\times m^l}(R)$.
%\\
%\underline{case 1} Assume that $j=k$.\\
%Then $A\cdot B$ by definition is just the usual matrix product of $A$ and $B$. Hence $A%\cdot B\in \Mat_{n^i\times n^l}(R)$ in this case.\\

\begin{enumerate}

\item[Case 1] Assume that $j\leq k$. Let 
\[A':=\begin{pmatrix} A&&\\&\ddots&\\&&A\end{pmatrix}\in \Mat_{m^{i+k-j}\times m^k}(R)\]
be the matrix with $m^{k-j}$ copies of $A$ on the diagonal and zeros elsewhere. We define $A\cdot B:=A'B\in \Mat_{m^{i+k-j}\times m^l}(R)$ where $A'B$ is the usual matrix product of $A'$ and $B$.

\medskip 

\item[Case 2]  Assume that $j\geq k$.
Let 
\[B':=\begin{pmatrix} B&&\\&\ddots&\\&&B\end{pmatrix}\in \Mat_{m^{j}\times m^{l+j-k}}(R)\]
be the matrix with $m^{j-k}$ copies of $B$ on the diagonal and zeros elsewhere. We define $A\cdot B=AB'\in \Mat_{m^{i}\times m^{l+j-k}}(R)$ where $AB'$ is the usual matrix product of $A$ and $B'$.

\end{enumerate}
\end{definition}

One checks easily that $(\Mf,\cdot)$ is a monoid whose identity element is $(1)\in \Mat_{1\times 1}(R)$. Further $A(B+C)=AB+AC$ and $(B+C)A=BA+CA$ whenever $B+C$ is defined, i.e. whenever $B$ and $C$ have the same size. However, one can show that there is no binary operation $+$ on $\Mf$ such that $(\Mf,+,\cdot)$ is a ring.

The theorem below will be used in the next section in order to identify weighted Leavitt path algebras with generalised corner skew Laurent polynomial rings.

\begin{theorem} \label{thm24}
Let $A$ be a $\Z^n$-graded unital ring. Assume that there is an $m\in \N$ and elements $T_1,\dots, T_n\in \Mat_{1\times m}(A_{\alpha_i})$ and $T_{-1},\dots, T_{-n}\in \Mat_{m\times 1}(A_{-\alpha_i})$ such that $T_iT_{-i}=1$ for any $1\leq i \leq n$. Then $A$ is graded isomorphic to a generalised corner skew Laurent polynomial ring.
\end{theorem}
\begin{proof} 
Set $R:=A_0$ and for any $g\in \N_0^n$,
\[p_g:=T^{g_n}_{-n}\dots T^{g_1}_{-1}T^{g_1}_{1}\dots T^{g_n}_{n}\in \Mat_{m^{g_1+\dots+g_n}\times m^{g_1+\dots+g_n}}(R)\]
(where the multiplication is taken in $\Mf(m,A)$). Moreover, set for any $g\in \Z^n$,
\[\tilde A_g:=p_{\widehat g} \Mat_{m^{\widehat g_1+\dots+\widehat g_n}\times m^{(\widehat{-g})_1+\dots+(\widehat{-g})_n}}(R)p_{\widehat{-g}}\]
and define for any $g,h\in \Z^n$ the map
\begin{align*}
\phi_{g,h}:\tilde A_g\times \tilde A_h&\rightarrow \tilde A_{g+h}\\
(\sigma,\tau)&\mapsto T^{(\widehat{g+h})_n}_{-n}\dots T^{(\widehat{g+h})_1}_{-1}T^{\widehat g_1}_{1}\dots T^{\widehat g_n}_{n}\sigma T^{(\widehat{-g})_n}_{-n}\dots T^{(\widehat{-g})_1}_{-1}\cdot\\
&~~~\cdot T^{\widehat h_1}_{1}\dots T^{\widehat h_n}_{n}\tau T^{(\widehat{-h})_n}_{-n}\dots T^{(\widehat{-h})_1}_{-1}T^{(\widehat{-(g+h)})_1}_{1}\dots T^{(\widehat{-(g+h)})_n}_{n}.
\end{align*}
One checks easily that $\mathcal A:=\{\tilde A_g\ | \ g\in \Z^n\}$ is a family of abelian groups and the $\phi_{g,h}$'s are biadditive maps which satisfy the conditions (i) and (ii) in Definition~\ref{def15}. Furthermore, the generalised group ring $\mathcal A[\Z^n]$ satisfies the conditions (i)-(iv) in Definition~\ref{def17}. Hence $\mathcal A[\Z^n]=R[t_1,\dots t_n,t_{-n},\dots,t_{-1},\phi]$ is a generalised corner skew Laurent polynomial ring. We will show that $A$ is graded isomorphic to $\mathcal A[\Z^n]$. Define the map
\begin{align*}
\psi:\mathcal A[\Z^n] &\longrightarrow A\\
\sum_{g\in \Z^n}t_1^{\widehat{g}_1}\dots t_n^{\widehat{g}_n}\sigma_g t_{-n}^{(\widehat{-g})_{n}}\dots  t_{-1}^{(\widehat{-g})_{1}}&\longmapsto \sum_{g\in \Z^n}T_1^{\widehat{g}_1}\dots T_n^{\widehat{g}_n}\sigma_gT_{-n}^{(\widehat{-g})_{n}}\dots T_{-1}^{(\widehat{-g})_{1}}.
\end{align*}
One checks easily that $\psi$ is a graded ring homomorphism. In order to show that $\psi$  is an isomorphism, it suffices to show that $\psi|_{\mathcal A[\Z^n]_g}:\mathcal A[\Z^n]_g\rightarrow A_g$ is a bijection for any $g\in \Z^n$. Let $g\in \Z^n$. Suppose 
\[T_1^{\widehat{g}_1}\dots T_n^{\widehat{g}_n}\sigma T_{-n}^{(\widehat{-g})_{n}}\dots T_{-1}^{(\widehat{-g})_{1}}=\psi(t_1^{\widehat{g}_1}\dots t_n^{\widehat{g}_n}\sigma t_{-n}^{(\widehat{-g})_{n}}\dots t_{-1}^{(\widehat{-g})_{1}})=0\]
for some $\sigma\in \mathcal A[\Z^n]_g$. By multiplying $T^{\widehat g_n}_{-n}\dots T^{\widehat g_1}_{-1}$ from the left and $T^{(\widehat{-g})_1}_{1}\dots T^{(\widehat{-g})_n}_{n}$ from the right, it follows that $p_{\widehat g}\sigma p_{\widehat{-g}}=0$. But clearly $p_{\widehat g}\sigma p_{\widehat{-g}}=\sigma$, since $\sigma\in \mathcal A[\Z^n]_g$ and $p_{\widehat g}$ and $p_{\widehat{-g}}$ are idempotents. Hence $\sigma=0$ and therefore 
$t_1^{\widehat{g}_1}\dots t_n^{\widehat{g}_n}\sigma t_{-n}^{(\widehat{-g})_{n}}\dots t_{-1}^{(\widehat{-g})_{1}}=0$. Thus we have shown that $\psi|_{\mathcal A[\Z^n]_g}:\mathcal A[\Z^n]_g\rightarrow A_g$ is injective. Suppose now that $a\in A_g$. Then
\[T_{-n}^{\widehat{g}_{n}}\dots T_{-1}^{\widehat{g}_{1}} a T_1^{(\widehat{-g})_1}\dots T_n^{(\widehat{-g})_n} \in p_{\widehat g} \Mat_{m^{\widehat g_1+\dots+\widehat g_n}\times m^{(\widehat{-g})_1+\dots+(\widehat{-g})_n}}(R)p_{\widehat{-g}}.\]
Clearly 
\begin{align*}
&\psi(\underbrace{t_1^{\widehat{g}_1}\dots  t_n^{\widehat{g}_n}(T_{-n}^{\widehat{g}_{n}}\dots T_{-1}^{\widehat{g}_{1}} a T_1^{(\widehat{-g})_1}\dots T_n^{(\widehat{-g})_n})t_{-n}^{(\widehat{-g})_{n}}\dots t_{-1}^{(\widehat{-g})_{1}}}_{\in \mathcal A[\Z^n]_g})\\
=&T_1^{\widehat{g}_1}\dots  T_n^{\widehat{g}_n}T_{-n}^{\widehat{g}_{n}}\dots T_{-1}^{\widehat{g}_{1}} a T_1^{(\widehat{-g})_1}\dots T_n^{(\widehat{-g})_n}T_{-n}^{(\widehat{-g})_{n}}\dots T_{-1}^{(\widehat{-g})_{1}}\\
=&a
\end{align*}
since $T_iT_{-i}=1$ for any $i\in \{1,\dots,n\}$. Thus we have shown that $\psi|_{\mathcal A[\Z^n]_g}:\mathcal A[\Z^n]_g\rightarrow A_g$ is surjective as well. This finishes the proof. 
\end{proof}

\subsection{Weighted Leavitt path algebras as generalised corner skew Laurent polynomial rings}\label{sec4}

Let $(E,w)$ be a finite weighted graph and $n$ the maximal weight of an edge in $E$. Recall that the standard grading of $L(E,w)$ is a $\Z^n$-grading such that $\deg(v)=0$, $\deg(e_i)=\alpha_i$ and $\deg(e_i^*)=-\alpha_i$ for any $v\in E^0$, $e \in E^{1}$ and $1\leq i\leq w(e)$ (we continue to denote the standard basis of $\mathbb Z^n$ by $\{\alpha_i \mid 1\leq i \leq n\}$).

\begin{theorem}\label{thm25}
Let $(E,w)$ be a finite weighted graph without sinks. Then $L(E,w)$ is graded isomorphic (with respect to its standard grading) to a generalised corner skew Laurent polynomial
ring.
\end{theorem}
\begin{proof}
Let $n$ be the maximal weight of an edge in $E$. Write $E^0=\{v_1,\dots,v_k\}$ and for any $1\leq j\leq k$, $s^{-1}(v_j)=\{e^{j,1},\dots, e^{j,n_j}\}$. For any $1\leq i\leq n$ and $1\leq j\leq k$ set 
\[T_i^j:=\begin{pmatrix}e_i^{j,1}&\dots&e_i^{j,n_j}\end{pmatrix}\in \Mat_{1\times n_j}(L(E,w)_{\alpha_i})\]
and
\[T_{-i}^j:=\begin{pmatrix}(e_i^{j,1})^*\\\vdots\\(e_i^{j,n_j})^*\end{pmatrix}\in \Mat_{n_j\times 1}(L(E,w)_{-\alpha_i})\]
where $e_i^{j,l}=(e_i^{j,l})^*=0$ if $i>w(e^{j,l})$. Moreover, for any $1\leq i\leq n$ set
\[T_i:=\begin{pmatrix}T_i^1&\dots&T_i^{k}\end{pmatrix}\in \Mat_{1\times m}(L(E,w)_{\alpha_i})\]
and
\[T_{-i}:=\begin{pmatrix}T_{-i}^1\\\vdots\\T_{-i}^k\end{pmatrix}\in \Mat_{m\times 1}(L(E,w)_{-\alpha_i})\]
where $m=n_1+\dots+n_k=|E^1|$. It follows from Definition~\ref{def3}(iii) that $T_iT_{-i}=v_1+\dots+v_k=1$ for any $i\in \{1,\dots,n\}$. Thus, by Theorem~\ref{thm24}, $L(E,w)$ is graded isomorphic to a generalised corner skew Laurent polynomial ring.\end{proof}

Theorems \ref{thm26} and \ref{thm28} below follow directly from Theorem~\ref{thm25} and Propositions~\ref{prop17} and~\ref{prop22}.

\begin{theorem} \label{thm26}
Let $(E,w)$ be a finite weighted graph without sinks. Then $L(E,w)$ is strongly graded with respect to its standard grading if and only if
\[\Mat_{1\times m}(L(E,w)_0)T_{-i}T_i\Mat_{m\times 1}(L(E,w)_0)=L(E,w)_0\]
for any $1\leq i \leq n$ where $m$, $n$, $T_i$ and $T_{-i}$ are defined as in the proof of Theorem~\ref{thm25}.
\end{theorem}

\begin{corollary}\label{cor27}
Let $(E,w)$ be a finite weighted graph without sinks and sources such that $w\equiv 1$. Then $L(E,w)$ is strongly graded with respect to its standard grading.
\end{corollary}
\begin{proof}
We use the same notation as in the proof of Theorem \ref{thm25}. Clearly 
\[T_{-1}T_1=\diag\big(r(e^{1,1}),\dots,r(e^{1,n_1}),r(e^{2,1}),\dots,r(e^{2,n_2}),\dots,r(e^{k,1}),\dots,r(e^{k,n_k})\big)\]
by Definition~\ref{def3}(iv). Since there are no sources in $E$, any vertex of $E$ appears as a diagonal entry in $T_{-1}T_1$. It follows that $\Mat_{1\times m}(L(E,w)_0)T_{-1}T_1 \Mat_{m\times 1}(L(E,w)_0)=L(E,w)_0$. Thus $L(E,w)$ is strongly graded with respect to its standard grading by Theorem \ref{thm26}. 
\end{proof}

\begin{theorem}\label{thm28}
Let $(E,w)$ be a finite weighted graph without sinks. Then $L(E,w)$ is graded von Neumann regular with respect to its standard grading if and only if $L(E,w)_0$ is von Neumann regular.
\end{theorem}

It has yet to be worked out whether there is a good description for the zero component  $L(E,w)_0$ (cf. Section 12). When $w\equiv 1$, this ring is well understood. 

\begin{corollary} \label{cor29}
Let $(E,w)$ be a finite weighted graph without sinks such that $w\equiv 1$. Then $L(E,w)$ is graded von Neumann regular with respect to its standard grading.
\end{corollary}
\begin{proof} By~\cite[pp. 168-169]{hazrat16} (or \cite[Corollary~2.1.16]{abrams-ara-molina}), $L(E,w)_0$, is an ultramatricial algebra, i.e. it is isomorphic to the union of an increasing chain of finite products of matrix algebras over a field. Hence $L(E,w)_0$ is von Neumann regular and thus, by Theorem~\ref{thm28}, $L(E,w)$ is graded von Neumann regular.
\end{proof}

\section{Local valuations}
\subsection{General results}
In this subsection we recall some notions and results from \cite[\S 7.1]{Raimund2}.

\begin{definition}\label{defval}
Let $R$ be a ring. A {\it valuation} on $R$ is a map $\nu:R\rightarrow \N_0\cup\{-\infty\}$ such that 
\begin{enumerate}[(i)]
\item $\nu(x)=-\infty\Leftrightarrow x=0$,
\item $\nu(x-y)\leq \max\{\nu(x),\nu(y)\}$ for any $x,y\in R$ and
\item $\nu(xy)=\nu(x)+\nu(y)$ for any $x,y\in R$.
\end{enumerate}
We use the conventions $-\infty<n$ for any $n\in\N_0$ and $(-\infty)+n=n+(-\infty)=-\infty$ for any $n\in\N_0\cup\{-\infty\}$.
\end{definition}

\begin{comment}
\begin{remark}
One easily checks that condition (ii) in Definition \ref{defval} is satisfied iff the conditions (iia) and (iib) below are satisfied.
\begin{enumerate}[({ii}a)]
\item $\nu(x+y)\leq \max\{\nu(x),\nu(y)\}$ for any $x,y\in R$.
\item $\nu(x)=\nu(-x)$ for any $x\in R$.
\end{enumerate}
\end{remark}

Recall that a {\it domain} is a nonzero ring without zero divisors. 
\begin{lemma}\label{lemval}
Let $R$ be a nonzero ring that has a valuation. Then $R$ is a domain.
\end{lemma}
\begin{proof}
Let $\nu$ be a valuation on $R$. Let $x,y\in R$ such that $xy=0$. Then $-\infty=\nu(0)=\nu(xy)=\nu(x)+\nu(y)$ and hence $\nu(x)=-\infty$ or $\nu(y)=-\infty$. Thus $x=0$ or $y=0$.
\end{proof}
\end{comment}

\begin{definition}
A ring with enough idempotents is a pair $(R,E)$ where $R$ is a ring and $E$ is a set
of nonzero orthogonal idempotents in $R$ such that the set of finite sums of
distinct elements of $E$ is a set of local units for $R$. Note that if $(R,E)$ is a ring with enough idempotents, then $R=\bigoplus_{e\in E}eR=\bigoplus_{f\in E}Rf=\bigoplus_{e,f\in E}eRf$. A ring with enough idempotents $(R,E)$ is called {\it connected} if $eRf\neq\{0\}$ for any $e,f\in E$.
\end{definition}

\begin{definition}\label{deflocval}
Let $(R,E)$ be a ring with enough idempotents. A {\it local valuation} on $(R,E)$ is a map $\nu:R\rightarrow \N_0\cup\{-\infty\}$ such that 
\begin{enumerate}[(i)]
\item $\nu(x)=-\infty\Leftrightarrow x=0$,
\item $\nu(x-y)\leq \max\{\nu(x),\nu(y)\}$ for any $x,y\in R$ and
\item $\nu(xy)=\nu(x)+\nu(y)$ for any $e\in E$, $x\in Re$ and $y\in eR$.
\end{enumerate}
A local valuation $\nu$ on $(R,E)$ is called {\it trivial} if $\nu(x)=0$ for any $x\in R\setminus\{0\}$ and {\it nontrivial} otherwise.
\end{definition}

%While any free algebra is an example of a ring with a valuation (define $\nu(x)$ as the degree of $x$ in the free generators), any path algebra is an example of a ring with enough idempotents a local valuation (let $E$ be the set of vertices and define $\nu(x)$ as the maximal length of a path appearing with nonzero coefficient in an expression of $x$ as linear combination of pairwise distinct paths).

Let $R$ be a ring. Recall that a left ideal $I$ of $R$ is called {\it essential} if $I \cap J=\{0\}\Rightarrow J=\{0\}$, for any left ideal $J$ of $R$. If $I$ is an essential left ideal of $R$, then we write $I\subseteq_e R$. For any $x\in R$ define the left ideal $\ann(x):=\{y\in R \mid yx=0\}$. The ring $R$ is said to be {\it left nonsingular} if for any $x\in R$, $\ann(x)\subseteq_e R~\Leftrightarrow x=0$. A right nonsingular ring is defined similarly. $R$ is called {\it nonsingular} if it is left and right nonsingular.

\begin{proposition}[{\cite[Proposition 37]{Raimund2}}]\label{proplocval1}
Let $(R,E)$ be a ring with enough idempotents that has a local valuation. Then $R$ is nonsingular. 
\end{proposition}
\begin{proof}
We show only left singularity of $R$ and leave the right singularity to the reader. Let $\nu$ be a local valuation on $(R,E)$ and $x\in R\setminus\{0\}$. Choose an $e\in E$ such that $ex\neq 0$. Suppose that $ye\in \ann(x)$ for some $y\in R$. Then 
\[\nu(ye)+\nu(ex)=\nu(yex)=\nu(0)=-\infty\]
and hence $ye=0$. This shows that $\ann(x)\cap Re=\{0\}$. But $Re\neq \{0\}$ since $e\in Re$. Hence $\ann(x)$ is not essential.
\end{proof}

Recall that a nonzero ring $R$ is called a {\it prime ring} if $IJ=\{0\}\Rightarrow (I=\{0\} \lor J=\{0\})$, for any ideals $I$ and $J$ of $R$. If $R$ has local units, then it is prime iff $xRy=\{0\}\Rightarrow(x=0 \lor y=0)$, for any $x,y\in R$.

\begin{proposition}[{\cite[Proposition 38]{Raimund2}}]\label{proplocval2}
Let $(R,E)$ be a nonzero, connected ring with enough idempotents that has a local valuation. Then $R$ is a prime ring. 
\end{proposition}
\begin{proof}
Let $\nu$ be a local valuation on $(R,E)$ and $x,y\in R\setminus\{0\}$. Clearly there are $e,f\in E$ such that $xe, fy\neq 0$. Since $(R,E)$ is connected, we can choose a $z\in eRf\setminus\{0\}$. Clearly 
\[\nu(xzy)=\nu(xezfy)=\nu(xe)+\nu(z)+\nu(fy)\geq 0.\]
Hence $xzy\neq 0$ and thus $xRy\neq \{0\}$.
\end{proof}

\begin{comment}
Recall that a ring $R$ is called {\it von Neumann regular} if for any $x\in R$ there is a $y\in R$ such that $xyx=x$. 

\begin{proposition}[{\cite[Proposition 39]{Raimund2}}]\label{proplocval3}
Let $(R,E)$ be a ring with enough idempotents that has a nontrivial local valuation. Then $R$ is not von Neumann regular. 
\end{proposition}
\begin{proof}
Let $\nu$ be a nontrivial local valuation on $(R,E)$. Choose an $x\in R$ such that $\nu(x)>0$. We may assume that $x\in eRf$ for some $e,f\in E$ because $(R,E)$ is a ring with enough idempotents. Let $y\in R$. Clearly 
\[\nu(xyx)=\nu(xfyex)=\nu(x)+\nu(fye)+\nu(x).\]
It follows that either $\nu(xyx)=-\infty$ (if $\nu(fye)=-\infty$) or $\nu(xyx)\geq 2\nu(x)$ (if $\nu(fye)\geq 0$). Thus $xyx\neq x$.
\end{proof}
\end{comment}

%We call a ring with enough idempotents $(A,E)$ a {\it $K$-algebra with enough idempotents} if $A$ is a $K$-algebra. 
Recall that a ring is called {\it semiprimitive} if its Jacobson radical is the zero ideal. 
\begin{proposition}[{\cite[Proposition 40]{Raimund2}}]\label{proplocval4}
Let $(R,E)$ be a connected ring with enough idempotents. Suppose $R$ is a $K$-algebra and there is a local valuation $\nu$ on $(R,E)$ such that $\nu(x)=0$ iff $x\in \Span(E)\setminus\{0\}$ where $\Span(E)$ denotes the linear subspace of $R$ spanned by $E$. Then $R$ is semiprimitive.
\end{proposition}
\begin{proof}
Let $\nu$ be the local valuation on $(R,E)$ such that $\nu(x)=0$ iff $x\in \Span(E)\setminus\{0\}$. Assume that the Jacobson radical $J$ of $R$ is not zero. Since $R=\bigoplus_{v,w\in E}eRf$, there are $e,f\in E$ and an $x'\in J\cap eRf\setminus\{0\}$. Since $R$ is connected, we can choose an element $z\in fRe\setminus\{0\}$. Then $x:=x'z\in J\cap eRe\setminus\{0\}$ since $\nu(x)=\nu(x'z)=\nu(x')+\nu(z)\geq 0$. Since $J$ does not contain any nonzero idempotents, it follows that $\nu(x)>0$ (if $\nu(x)=0$, then $x=ke$ for some $k\in K$ and hence $J$ contains the nonzero idempotent $e$). Since $x\in J$, we have that $x$ is left quasi-regular, i.e. there is a $y\in R$ such that $x+y=yx$. By multiplying $e$ from the right and from the left one gets $x+eye=eyex$. Hence we may assume that $y\in eRe$. It follows that 
\[\max\{\nu(x),\nu(y)\}\geq \nu(x+y)=\nu(yx)=\nu(x)+\nu(y).\]
This implies that $\nu(y)=0$ and hence $y=ke$ for some $k\in K$. It follows that $y=yx-x=kex-x=kx-x=(k-1)x$. But this yields a contradiction since $\nu(y)=0$ but either $\nu((k-1)x)=-\infty$, if $k=1$, or $\nu((k-1)x)=\nu(x)>0$, if $k\neq 1$ (note that $\nu((k-1)x)=\nu((k-1)ex)=\nu((k-1)e)+\nu(x)$). Thus the Jacobson radical of $R$ is zero.
\end{proof}

\subsection{Applications to weighted Leavitt path algebras}
Throughout this subsection $(E,w)$ denotes a weighted graph. Note that $(L(E,w),E^0)$ is a ring with enough idempotents. We say that $(E,w)$ {\it satisfies Condition (LV)} if $w(e)\geq 2$ for any  $e\in E^{1}$ and $\#\{e\in s^{-1}(v)\mid w(e)=w(v)\}\geq 2$ for any  $v\in E^0_{\reg}$. 

As in \S 4.1 let $K\X$ denote the free $K$-algebra on the set $X=\{v,e_i,e_i^*\mid v\in E^0,e\in E^1,1\leq i\leq w(e)\}$ and $K\X_{\nod}$ the linear subspace of $K\X$ generated by the nod-paths. Let $\NF:L(E,w)\rightarrow K\X_{\nod}$ be the isomorphism of $K$-vector spaces defined in the proof of Theorem \ref{thmbasis}. For an $a\in L(E,w)$ we define its {\it support} $\supp(a)$ as the set of all nod-paths which appear in $\NF(a)$ with nonzero coefficient. Recall that the length of a path $p$ is denoted by $|p|$. 

\begin{theorem}[{\cite[Proposition 40]{hazrat-preusser}}]\label{thmlocval}
If $(E,w)$ satisfies Condition (LV), then the map
\begin{align*}
\nu:L(E,w)&\rightarrow \N_0\cup\{-\infty\}\\
a&\mapsto \max\{|p|\ | \ p\in \supp(a)\}.
\end{align*}
is a local valuation on $(L(E,w),E^0)
$. Here we use the convention $\max(\emptyset)=-\infty$.
\end{theorem}

\begin{corollary}[{\cite[Theorem 47]{hazrat-preusser}}]\label{corlocval1}
If $(E,w)$ satisfies Condition (LV), then $L(E,w)$ is nonsingular. 
\end{corollary}
\begin{proof}
Follows from Proposition \ref{proplocval1} and Theorem \ref{thmlocval}.
\end{proof}

\begin{corollary}[{\cite[Theorem 46]{hazrat-preusser}}]\label{corlocval2}
If $(E,w)$ satisfies Condition (LV), then $L(E,w)$ is a prime ring. 
\end{corollary}
\begin{proof}
Clearly $L(E,w)$ is not the zero ring since $E$ is nonempty. Let $u,v\in E^0$. Since $E$ is connected, there is a d-path $p$ from $u$ to $v$. Let $\nu$ be the local valuation on $(L(E,w),E^0)$ defined in Theorem \ref{thmlocval}. Clearly $\nu(p)=|p|\geq 0$ since $\nu(x)=0$ for any $x\in E^0$ and $\nu(e_i)=\nu(e_{i}^*)=1$ for any $e\in E^1$ and $1\leq i\leq w(e)$. Hence $uL(E,w)v\neq \{0\}$ for any $u, v\in E^0$ and therefore $(L(E,w),E^0)$ is a nonzero, connected ring with enough idempotents. It follows from Proposition \ref{proplocval2} that $L(E,w)$ is a prime ring.
\end{proof}

\begin{comment}
\begin{corollary}\label{corlocval3}
If $(E,w)$ satisfies Condition (LV), then $L(H)$ is not von Neumann regular. 
\end{corollary}
\begin{proof}
Let $\nu$ be the local valuation on $(L(H),H^0)$ defined in Theorem \ref{thmlocval}. Choose an $h\in H^1$, an $i\in I_h$ and a $j\in J_h$. Then $\nu(h_{ij})=1$ and therefore $\nu$ is nontrivial. It follows from Proposition \ref{proplocval3} that $L(H)$ is not von Neumann regular.
\end{proof}
\end{comment}

\begin{corollary}[{\cite[Theorem 50]{hazrat-preusser}}]\label{corlocval4}
If $(E,w)$ satisfies Condition (LV), then $L(E,w)$ is semiprimitive. 
\end{corollary}
\begin{proof}
Let $\nu$ be the local valuation on $(L(E,w),E^0)$ defined in Theorem \ref{thmlocval}. Then clearly $\nu(x)=0$ iff $x\in \Span(E^0)\setminus\{0\}$ where $\Span(E^0)$ denotes the linear subspace of $L(E,w)$ spanned by $E^0$. Moreover, $(L(E,w),E^0)$ is connected since $E$ is connected (see the proof of Corollary \ref{corlocval2}). It follows from Proposition \ref{proplocval4} that $L(E,w)$ is semiprimitive.
\end{proof}

%\begin{lemma}\label{freg}
%Suppose $H$ satisfies Condition (LV), $H$ is connected and $H^1\neq\emptyset$. If $I$ is a nonzero ideal of $L(H)$, then for any $n\in\N$ and $v,w\in H^0$ there is an $a\in I\cap v L(H)w$ such that $\nu(a)> n$.
%\end{lemma}
%\begin{proof}
%Let $I$ be a nonzero ideal of $L(H)$, $n\in\N$ and $v,w\in H^0$. Choose a nonzero element $a'\in I$. Then there are $z_1,z_2\in H^0$ such that $z_1a'z_2\neq 0$. Since $H$ is connected there are d-paths $p$ and $q$ that $s(p)=v$, $r(p)=z_1$, $s(q)=z_2$ and $r(q)=w$. Clearly we assume that $|p|, |q|>n$ (choose d-paths $p'$ and $q'$ that $s(p')=z_1$, $r(p')=v$, $s(q')=w$ and $r(q')=z_2$ and replace $p$ and $q$ by $(pp')^kp$ resp. $(qq')^lq$ for sufficiently large $k$ and $l$). It follows that $\nu(p)=|p|,\nu(q)=|q|>n$. Set $a:=pa'q\in I\cap vL(H)w$. Then $\nu(a)=\nu(pa'q)=\nu(p)+\nu(z_1a'z_2)+\nu(q)>n$.
%\end{proof}

\subsection{Classification of the weighted graphs $(E,w)$ such that $L(E,w)$ is a domain}

Recall that a {\it domain} is a nonzero ring without zero divisors. A weighted graph $(E,w)$ satisfying Condition (LV) is called an {\it LV-rose} if $|E^0|=1$. 

\begin{theorem}[{\cite[Theorem 41]{hazrat-preusser}}]\label{thmdomain}
Let $(E,w)$ be a weighted graph. Then $L(E,w)$ is a domain if and only if $(E,w)$ is an LV-rose or the weighted graph $\xymatrix{\bullet\ar@(ur,rd)^1}$.
\end{theorem}
\begin{proof}
One checks easily that $L(E,w)$ has zero divisors if $(E,w)$ is neither an LV-rose nor the weighted graph $\xymatrix{\bullet\ar@(ur,rd)^1}$. On the other hand, if $(E,w)$ is an LV-rose, then the local valuation $\nu$ on $(L(E,w), E^0)$ defined in Theorem \ref{thmlocval} is a valuation and hence $L(E,w)$ is a domain. Moreover, $L(\xymatrix{\bullet\ar@(ur,rd)^1})\cong K[x,x^{-1}]$ is a domain.
\end{proof}

\section{The $\V$-monoid and $K_0$}
Recall from Section 4 that $\G$ denotes the category whose objects are the weighted graphs and whose morphisms are the complete weighted graph homomorphisms, that $\A$ denotes the category of $K$-algebras and that $L:\G\to \A$ is a functor which commutes with direct limits. In Subsection 9.1 we define functors $\V:\A\rightarrow \M$ and $\cM:\G\rightarrow \M$ where $\M$ denotes the category of abelian monoids. In Subsection 9.2 we recall some universal ring constructions by G. Bergman which are used in the proof of Theorem \ref{thmmon}. In Subsection 9.3 we sketch the proof of Theorem \ref{thmmon} which states that $\V\circ L\cong \cM$ and that $L(E,w)$ is left and right hereditary provided $(E,w)$ is finite. In Subsection 9.4 we compute the Grothendieck group of a weighted Leavitt path algebra.

\subsection{The functors $\V$ and $\cM$}
\begin{definition}
Let $A$ be a $K$-algebra. Let $\Mat_\infty(A)$ be the directed union of the rings $\Mat_n(A)~(n\in\N)$, where the transition maps $\Mat_n(A)\rightarrow \Mat_{n+1}(A)$ are given by $x\mapsto \begin{pmatrix}x&0\\0&0\end{pmatrix}$. Let $I(\Mat_\infty(A))$ denote the set of all idempotent elements of $\Mat_\infty(A)$. If $e, f\in I(\Mat_\infty(A))$, write $e\sim f$ iff there are $x,y\in \Mat_\infty(A)$ such that
$e = xy$ and $f = yx$. Then $\sim$ is an equivalence relation on $I(\Mat_\infty(A))$. Let $\V(A)$ be the set of all $\sim$-equivalence classes, which becomes an abelian monoid by defining
\[[e]+[f]=\left[\begin{pmatrix}e&0\\0&f\end{pmatrix}\right]\]
for any $[e],[f]\in \V(A)$. If $\phi:A\rightarrow B$ is a morphism in $\A$, let $\V(\phi):\V(A)\rightarrow \V(B)$ be the canonical monoid homomorphism induced by $\phi$. One checks easily that $\V:\A\rightarrow \M$ is a functor that commutes with direct limits.
% Choose an associative, unital ring $R$ containing $A$ as a two-sided ideal (if $A$ is unital, one can choose $R=A$). Denote by $FP(A,R)$ the class of finitely generated projective right $R$-modules $P$ for which $P = PA$. Then $\V(A)$ is defined as the monoid of isomorphism classes of the elements of $FP(A,R)$ where the binary operation is given by the direct sum of modules. This definition of $\V(A)$ does not depend on the ring $R$, cf. \cite[Definition 3.2.1]{abrams-ara-molina}. If $f:A\rightarrow B$ is a morphism in $\A$, then $V$ defines a functor from the category of associative rings to the category of abelian monoids, which commutes with direct limits.
\end{definition}

\begin{remark}\label{rempro}
For a $K$-algebra $A$ with local units there is the following alternative description of the monoid $\V(A)$. Recall that a left $A$-module $M$ is called {\it unital} if $AM=M$. Let $\V'(A)$ denote the set of isomorphism classes of finitely generated projective unital left $A$-modules, which becomes an abelian monoid by defining $[P]+[Q]:=[P\oplus Q]$ for any $[P],[Q]\in \V'(A)$. Then $\V'(A)\cong \V(A)$ as abelian monoids, see \cite[Subsection 4A]{ara-hazrat-li-sims}.
\end{remark}

\begin{definition}\label{defM}
Let $(E,w)$ be a weighted graph. For any $v\in E^0$ write $w(s^{-1}(v))=\{w_1(v),\dots,$ $w_{k_v}(v)\}$ where $k_v\geq 0$ and $w_1(v)<\dots<w_{k_v}(v)$ (hence $k_v$ is the number of different weights of edges in $s^{-1}(v)$). Let $\cM(E,w)$ be the abelian monoid presented by the generating set $\{v,q_1^v,\dots,q^v_{k_v-1}\mid v\in E^0\}$ and the relations
\begin{equation}
q^v_{i-1}+(w_i(v)-w_{i-1}(v))v=q_i^v+\sum_{\substack{e\in s^{-1}(v),\\w(e)=w_i(v)}}r(e)\quad\quad(v\in E^0,1\leq i\leq k_v)
\end{equation}
where $q^v_0=q^v_{k_v}=w_0(v)=0$. If $\phi:(E,w)\rightarrow (E',w')$ is a morphism in $\G$, then there is a unique monoid homomorphism $\cM(\phi):\cM(E,w)\rightarrow \cM(E',w')$ such that $\cM(\phi)([v])=[\phi^0(v)]$ and $\cM(\phi)([q_i^v])=[q_i^{\phi^0(v)}]$ for any $v\in E^0$ and $1\leq i \leq k_v-1$. One checks easily that $\cM:\G\rightarrow \M$ is a functor that commutes with direct limits.
\end{definition}

%\begin{remark}
%If $k_v\leq 1$ for any $v\in E^0$, then $M(E,w)$ is the abelian monoid $M_E$ defined in \cite[Theorem 5.21]{hazrat13}.
%\end{remark}

\subsection{Some universal ring constructions by G. Bergman}
In this subsection all rings are assumed to be unital. Let $k$ be a commutative ring and $R$ a $k$-algebra (i.e. $R$ is a ring given with a homomorphism of $k$ into its center). By an $R$-module we mean a left $R$-module. An {\it $R$-ring$_k$} is a $k$-algebra $S$ given with a $k$-algebra homomorphism $R\rightarrow S$. In \cite{bergman74}, G. Bergman described the following two key constructions: 

\begin{itemize}
\item ADJOINING MAPS Let $M$ be any $R$-module and $P$ a finitely generated projective $R$-module. Then there exists an $R$-ring$_k$ $S$, having a universal module homomorphism $f: M\otimes S \rightarrow P\otimes S$, see \cite[Theorem 3.1]{bergman74}. $S$ can be obtained by adjoining to $R$ a family of generators subject to certain relations, see \cite[Proof of Theorem 3.1]{bergman74}.\\
\item IMPOSING RELATIONS Let $M$ be any $R$-module, $P$ a projective $R$-module and $f: M\rightarrow P$ any module homomorphism. Then there exists an $R$-ring$_k$ $S$ such that $f\otimes S = 0$ and universal for that property. $S$ can be chosen to be a quotient of $R$, see \cite[Proof of Theorem 3.2]{bergman74}.
\end{itemize}

Using the key constructions above Bergman described more complicated constructions. Two of them are used in this paper:

\begin{itemize}
\item ADJOINING ISOMORPHISMS Given two finitely generated projective $R$-modules $P$
and $Q$, one can adjoin a universal isomorphism between $P\otimes $ and $Q\otimes$ by first
freely adjoining maps $i:P\otimes \rightarrow Q \otimes$ and $i^{-1}:Q\otimes \rightarrow P \otimes$ (via ADJOINING MAPS) and then setting $1_{Q\otimes}-ii^{-1}$ and $1_{P\otimes}-i^{-1}i$ equal to $0$ (via IMPOSING RELATIONS), see \cite[p. 38]{bergman74}. Bergman denoted the resulting $R$-ring$_k$ by $R\langle i, i^{-1} : \overline{P}\cong\overline{Q}\rangle$.\\
\item ADJOINING IDEMPOTENT ENDOMORPHISMS Given a finitely generated projective $R$-module $P$, one can adjoin a universal idempotent endomorphism of $P\otimes $ by first
freely adjoining a map $e:P\otimes \rightarrow P \otimes$ (via ADJOINING MAPS) and then setting $e-e^2$ equal to $0$ (via IMPOSING RELATIONS), see \cite[p. 39]{bergman74}. Bergman denoted the resulting $R$-ring$_k$ by $R\langle e: \overline{P}\rightarrow\overline{P};e^2=e\rangle$. Note that the adjoined idempotent endomorphism $e$ yields a universal direct sum decomposition $P\otimes = \ker(e)\oplus \im(e)$.
\end{itemize}

Set $S:=R\langle i, i^{-1} : \overline{P}\cong\overline{Q}\rangle$ and $T:=R\langle e: \overline{P}\rightarrow\overline{P};e^2=e\rangle$. Bergman proved the following (for these results he required that $k$ is a field and that $P$ and $Q$ are nonzero): The abelian monoid $\V'(S)$ (see Remark \ref{rempro}) may be obtained from $\V'(R)$ by imposing one relation $[P] = [Q]$. The abelian monoid $\V'(T)$ may be obtained from $\V'(R)$ by adjoining two new generators $[P_1]$ and $[P_2]$ and one relation $[P_1] + [P_2] =[P]$. Furthermore, the left global dimension of $S$ (resp. $T$) equals the left global dimension of $R$, unless the left global dimension of $R$ is $0$, in which case the left global dimension of $S$ (resp. $T$) is $\leq 1$. See \cite[Theorems 5.1, 5.2 and the last paragraph of p. 48]{bergman74}.

\subsection{$\V\circ L$ and $\cM$ are isomorphic}
\begin{comment}
The lemma below will be used in the proof of Theorem \ref{thmmon}.
\begin{lemma}\label{lemmon}
Let $G$ be an abelian group (resp. an abelian monoid) presented by a generating set $X$ and relations 
\[l_i=r_i~ (i\in I)\text{ and }y=\sum_{x\in X\setminus\{y\}}n_xx\]
where for any $i\in I$, $l_i$ and $r_i$ are elements of the free abelian group (resp. the free abelian monoid) $G\X$ generated by $X$, $y$ is an element of $X$, the $n_x$ are integers (resp. nonnegative integers) and only finitely of them are nonzero. Let $G\langle X\setminus\{y\}\rangle$ be the free abelian group (resp. the free abelian monoid) generated by $X\setminus \{y\}$ and $f:G\X\rightarrow G\langle X\setminus\{y\}\rangle$ the homomorphism which maps each $x\in X\setminus\{y\}$ to $x$ and $y$ to $\sum_{x\in X\setminus\{y\}}n_xx$. Then $G$ is also presented by the generating set $X\setminus\{y\}$ and the relations $f(l_i)=f(r_i)~(i\in I)$.
\end{lemma}
\begin{proof}
Straightforward.
\end{proof}
\end{comment}

\begin{theorem}[{\cite[Theorem 14]{preusser-1}}]\label{thmmon}
$\V\circ L\cong \cM$. Moreover, if $(E,w)$ is a finite weighted graph, then $L(E,w)$ is left and right hereditary.
\end{theorem}
\begin{proofs}
In Part I below we define a natural transformation $\theta:\cM\rightarrow \V\circ L$. In Part II we explain why $\theta$ is a natural isomorphism and why $L(E,w)$ is left and right hereditary provided that $(E,w)$ is finite.\\
\\  
{\bf Part I}
Let $(E,w)$ be a weighted graph and $v\in E_{\reg}^0$. Write $s^{-1}(v)=\{e^{v,1}, \dots, e^{v,n(v)}\}$ where $w(e^{v,1})\leq\dots \leq w(e^{v,n(v)})$. Let $X_v$ be defined as in Remark \ref{remwlpa2}. Then $X_v$ has the upper triangular block form
\[X_v=\begin{pmatrix}
X_v^{1,1}&X_v^{1,2}&X_v^{1,3}&\dots&X_v^{1,k_v}\\0&X_v^{2,2}&X_v^{2,3}&\dots&X_v^{2,k_v}\\0&0&X_v^{3,3}&\dots&X_v^{3,k_v}\\\vdots&\vdots&\ddots&\ddots&\vdots\\0&0&\dots&0&X_v^{k_v,k_v}
\end{pmatrix}\]
where $k_v$ is the number of different weights of edges in $s^{-1}(v)$ and none of the matrices $X_v^{i,j}$ has a zero entry. For any $1\leq l \leq k_v-1$ define the matrix
\begin{comment}
\[X_{v,l}=\begin{pmatrix}
X_v^{l+1,l+1}&X_v^{l+1,l+2}&X_v^{l+1,l+3}&\dots&X_v^{l+1,k_v}\\0&X_v^{l+2,l+2}&X_v^{l+2,l+3}&\dots&X_v^{l+2,k_v}\\0&0&X_v^{l+3,l+3}&\dots&X_v^{l+3,k_v}\\\vdots&\vdots&\ddots&\ddots&\vdots\\0&0&\dots&0&X_v^{k_v,k_v}
\end{pmatrix}\]
\end{comment}
\[X_{v,l}=\begin{pmatrix}X_v^{1,l+1}&\dots&X_v^{l,k_v}\\\vdots&\ddots&\vdots\\X_v^{l,l+1}&\dots&X_v^{l,k_v}
\end{pmatrix}\]
and set $\epsilon_{v,l}:=X_{v,l}X_{v,l}^*$. Here $X_{v,l}^*$ is the matrix one obtains by transposing $X_{v,l}$ and applying the involution $*$ to each entry.
One can show that $\epsilon_{v,l}$ is an idempotent matrix for any $1\leq l\leq k_v-1$. There is a unique monoid homomorphism $\theta_{(E,w)}:\cM(E,w)\rightarrow \V(L(E,w))$ such that $\theta_{(E,w)}(v)=[(v)]$ and $\theta_{(E,w)}(q_l^v)=[\epsilon_{v,l}]$ for any $v\in E^0$ and $1\leq l\leq k_v-1$. It is an easy exercise to show that $\theta:\cM\rightarrow \V\circ L$ is a natural transformation.\\
\\
{\bf Part II} It remains to show that the natural transformation $\theta:\cM\rightarrow \V\circ L$ defined in Part I is a natural isomorphism, i.e. that $\theta_{(E,w)}:\cM(E,w)\rightarrow \V(L(E,w))$ is an isomorphism for any weighted graph $(E,w)$. By Lemma \ref{lemdirlim} any weighted graph is a direct limit of a direct system of finite weighted graphs. Hence is suffices to show that $\theta_{(E,w)}$ is an isomorphism for any finite weighted graph $(E,w)$ (note that $\cM$, $\V$ and $L$ commute with direct limits). So let $(E,w)$ be a finite weighted graph. Set $B_0:=K^{E^0}$. We denote by $\alpha_v$ the element of $B_0$ whose $v$-component is $1$ and whose other components are $0$. Let $v_1, \dots, v_m$ be the distinct elements of $E_{\reg}^0$. Let $1\leq t \leq m$ and assume that $B_{t-1}$ has already been defined. We define a $K$-algebra $B_t$ as follows. Set $C_{t,0}:=B_{t-1}$ and let $\beta^{t,0}:C_{t,0}\rightarrow C_{t,0}$ be the map sending any element to $0$. For $1\leq l \leq k_{v_t}-1$ define inductively $C_{t,l}:=C_{t,l-1}\langle \beta^{t,l}: \overline{O_{t,l}}\rightarrow \overline{O_{t,l}};(\beta^{t,l})^2=\beta^{t,l}\rangle $ (see Subsection 9.2 or \cite[p. 39]{bergman74}) where 
\[O_{t,l}=\im (\beta^{t,l-1})\oplus \bigoplus_{h=w_{l-1}(v_t)+1}^{w_l(v_t)}  \alpha_{v_t}C_{t,l-1}.\]
Set $D_{t,0}:=C_{t,k_{v_t}-1}$. For $1\leq l \leq k_{v_t}-1$ define inductively $D_{t,l}:=D_{t,l-1}\langle \gamma^{t,l},(\gamma^{t,l})^{-1}:\overline{P_{t,l}}\cong\overline{Q_{t,l}}\rangle $ (see Subsection 9.2 or \cite[p. 38]{bergman74}) where 
\[P_{t,l}=\bigoplus_{h=n_{l-1}(v_t)+1}^{n_l(v_t)}\alpha_{r(e^{h,v_t})}D_{t,l-1}\text{ and }Q_{t,l}=\ker(\beta^{t,l}).\]
Finally define $B_t:=D_{t,k_{v_t}-1}\langle \gamma^{t,k_{v_t}},(\gamma^{t,k_{v_t}})^{-1}:\overline{P_{t,k_{v_t}}}\cong\overline{Q_{t,k_{v_t}}}\rangle $ where
\begin{align*}
&P_{t,l}=\bigoplus_{h=n_{k_{v_t}-1}(v_t)+1}^{n_{k_{v_t}}(v_t)}\alpha_{r(e^{h,v_t})}D_{t,k_{v_t}-1}\text{ and }\\
&Q_{t,k_{v_t}}=\im(\beta^{t,k_{v_t}-1})\oplus\bigoplus_{h=w_{k_{v_t}-1}(v_t)+1}^{w_{k_{v_t}}(v_t)}  \alpha_{v_t}D_{t,k_{v_t}-1}.
\end{align*}
Then one can show that $L(E,w)\cong B_m$. It follows from \cite[Theorems 5.1, 5.2]{bergman74} that $\cM(E,w)\cong \V'(B_m)\cong \V'(L(E,w))\cong \V(L(E,w))$. One checks easily that the monoid isomorphism $\cM(E,w)\rightarrow \V(L(E,w))$ one gets in this way is precisely $\theta_{(E,w)}$. Furthermore, the left global dimension of $B_m\cong L(E,w)$ is $\leq 1$ by \cite[Theorems 5.1, 5.2]{bergman74}, i.e. $L(E,w)$ is left hereditary. Since $L(E,w)$ is a ring with involution, we have $L(E,w)\cong L(E,w)^{op}$. Thus $L(E,w)$ is also right hereditary.
\end{proofs}

\begin{corollary}[{\cite[Corollary 15]{preusser-1}}]\label{cornum}
Let $(E,w)$ be a weighted graph. If there is a vertex $v\in E^0$ that emits edges of different weights, then $|\V(L(E,w))|=\infty$.
\end{corollary}
\begin{proof}
Let $F$ denote the free abelian monoid generated by the set $\{v,q_1^v,\dots,q^v_{k_v-1}\mid v\in E^0\}$ and $\sim$ the congruence on $F$ defined by the relations (10) in Definition \ref{defM}. Let $v\in E^0$ be a vertex such that $k_v>1$. For any $n\in \N_0$ let $[nq_1^v]$ denote the $\sim$-congruence class of $nq^v_1$. In $F$ one cannot write $nq_1^v$ as $x+y$ where $x\in F$ and $y$ is the left or right hand side of one of the relations (10) (note that in the left hand side as well as in the right hand side of each of the relations (10) a nonempty sum of vertices appears). Hence $[nq_1^v]=\{nq_1^v\}$ for any $n\in \N_0$ (i.e. each $nq_1^v$ is only congruent to itself). Therefore the elements $[nq_1^v]~(n\in\N_0)$ are pairwise distinct in $\cM(E,w)$. It follows from Theorem \ref{thmmon} that $|\V(L(E,w))|=\infty$.
\end{proof}

\begin{corollary}[{\cite[Corollary 16]{preusser-1}}]
Let $(E,w)$ be an LV-rose such that there are edges of different weights. Then $L(E,w)$ is a domain that is neither $K$-algebra isomorphic to an unweighted Leavitt path algebra nor to a Leavitt algebra.
\end{corollary}
\begin{proof}
By Theorem \ref{thmdomain}, $L(E,w)$ is a domain. Clearly $(E,w)$ does not satisfy Condition (LPA1) and hence $L(E,w)$ is not isomorphic to an unweighted Leavitt path algebra by Theorem \ref{thmm2}. It remains to show that $L(E,w)$ is not isomorphic to a Leavitt algebra $L(m,n)$ where $1\leq m<n$. By Example \ref{wlpapp} and Theorem \ref{thmmon} we have $\V(L(m,n))\cong \N_0/\langle m=n\rangle$ and hence $|\V(L(m,n))|=n<\infty$. %If $k=0$, then $\V(L(n,n+k))\cong \N_0$ and therefore $L(E,w)$ is not isomorphic to $L(n,n+k)$ by the previous paragraph. Suppose now that $k\geq 1$. 
But by Corollary \ref{cornum} we have $|\V(L(E,w))|=\infty$. Thus $L(E,w)$ is not isomorphic to a Leavitt algebra $L(m,n)$.
\end{proof}

\subsection{The Grothendieck group $K_0(L(E,w))$}
One can use the adjacency matrix and the weighted identity matrix of a weighted graph $(E,w)$ to describe $K_0(L(E,w))$. We define those matrices below.
\begin{definition}
Let $(E,w)$ be a weighted graph. The {\it adjacency matrix of $(E,w)$} is the matrix $N\in \N_0^{E^0\oplus E^0}$ whose entry at position $(u,v)$ is the number of edges from $u$ to $v$. The {\it weighted identity matrix of $(E,w)$} is the matrix $I_w\in \N_0^{E^0\oplus E^0}$ whose entry at position $(u,v)$ is $w(v)$ if $u=v$ and $0$ otherwise.
\end{definition}

Let $(E,w)$ be a weighted graph. Denote the transpose of its adjacency matrix $N$ by $N^t$. Multiplying the matrix $N^t-I_w$ from the left defines a group homomorphism $\Z^{E^0}\rightarrow \Z^{E^0}$ where $\Z^{E^0}$ is the direct sum of copies of $\Z$ indexed by $E^0$. The theorem below shows that the cokernel of this map is the Grothendieck group of $L(E,w)$.

\begin{theorem}[{\cite[Theorem 18]{preusser-1}}]\label{thmkzero}
Let $(E,w)$ be a weighted graph. Then
\[K_0(L(E,w))\cong \coker(N^t-I_w: \Z^{E^0}\rightarrow \Z^{E^0}).\]
\end{theorem}
\begin{proof}
Since $L(E,w)$ is a ring with local units, $K_0(L(E,w))$ is the group completion $(\V(L(E,w)))^+$ of the abelian monoid $\V(L(E,w))$, see \cite[p. 77]{abrams-ara-molina}. By Theorem \ref{thmmon}, $(\V(L(E,w)))^+\cong (\cM(E,w))^+$. It follows from \cite[Equation (45)]{hazrat13} that the abelian group $(\cM(E,w))^+$ is presented by the generating set $\{v,q_1^v,\dots,q^v_{k_v-1}\mid v\in E^0\}$ and the relations
\begin{equation*}
q^v_{i-1}+(w_i(v)-w_{i-1}(v))v=q_i^v+\sum_{\substack{e\in s^{-1}(v),\\w(e)=w_i(v)}}r(e)\quad\quad(v\in E^0,1\leq i\leq k_v)
\end{equation*}
where $q^v_0=q^v_{k_v}=0$. We can rewrite the relations above in the form 
\begin{align*}
&q_i^v=q^v_{i-1}+(w_i(v)-w_{i-1}(v))v-\sum_{\substack{e\in s^{-1}(v),\\w(e)=w_i(v)}}r(e)\quad\quad(v\in E^0,1\leq i\leq k_v).
\end{align*}
By successively applying Tietze transformations we get that $(\cM(E,w))^+$ is presented by the generating set $E^0$ and the relations 
\[w(v)v=\sum_{e\in s^{-1}(v)}r(e)\quad\quad(v\in E^0).\] 
Hence $(\cM(E,w))^+\cong \Z^{E^0}/H$ where $H$ is the subgroup of $\Z^{E^0}$ generated by the set \[\{\sum_{e\in s^{-1}(v)}\alpha_{r(e)}-w(v)\alpha_{v}\mid v\in E^0\}\]
(where for a vertex $v$, $\alpha_v$ denotes the element of $\Z^{E^0}$ whose $v$-component is $1$ and whose other components are $0$). One checks easily that $H$ is the image of the homomorphism $N^t-I_w: \Z^{E^0}\rightarrow \Z^{E^0}$. Thus 
\[K_0(L(E,w))\cong (\V(L(E,w)))^+\cong (\cM(E,w))^+\cong \Z^{E^0}/H=\coker(N^t-I_w: \Z^{E^0}\rightarrow \Z^{E^0}).\]
\end{proof}

\begin{example}\label{exkzero}
Consider the weighted graphs\\
\[
(E,w):\xymatrix@C+15pt{ u& v\ar[l]_{e,1}\ar[r]^{f,2}& x}\quad\text{ and }\quad(E,w'):\xymatrix@C+15pt{ u& v\ar[l]_{e,2}\ar[r]^{f,2}& x}.
\]\\
Note that by Section 5, $L(E,w)\cong L(F)\cong \Mat_3(K)\oplus \Mat_3(K)$ where $F$ is the unweighted graph
\[\xymatrix@C+15pt{ &\bullet\ar[d]&\\F:\quad\bullet& \bullet\ar[l]\ar[r]& \bullet,}\]
while $L(E,w')$ is not isomorphic to an unweighted Leavitt path algebra.
By Theorem \ref{thmmon} we have
\[\V(L(E,w))\cong\langle u,v,x,q\mid v=q+u,q+v=x\rangle\cong \N_0^2 ~\text{ and }~\V(L(E,w'))\cong\langle u,v,x\mid 2v=u+x\rangle.\]
It follows that $\V(L(E,w))\not\cong \V(L(E,w'))$ since $\V(L(E,w))$ is a refinement monoid but $\V(L(E,w'))$ is not. On the other hand we have
\[K_0(L(E,w))\cong K_0(L(E,w'))\cong\langle u,v,x\mid 2v=u+x\rangle\cong\Z^2\]
by Theorem \ref{thmkzero}.
\end{example}

\section{The graded $\V$-monoid and $K_0^{\gr}$}
Throughout this section $\Gamma$ denotes a group with identity $\epsilon$.

\subsection{Basic definitions and results}
Let $R$ be a ring. Recall that a left $R$-module $M$ is called {\it unital} if $RM=M$. We denote by $R$-$\Mod$ the category of unital left $R$-modules. Furthermore, we denote by $R$-$\Mod_{\proj}$ the full subcategory of $R$-$\Mod$ whose objects are the projective objects of $R$-$\Mod$ that are finitely generated as a left $R$-module. If $R$ has local units, we define
\[\V(R)=\{[P]\mid P\in R\text{-}\Mod_{\proj}\}\]
where $[P]$ denotes the isomorphism class of $P$ as a left $R$-module. $\V(R)$ becomes an abelian monoid by defining $[P]+[Q]=[P\oplus Q]$. The Grothendieck group $K_0(R)$ is the group completion of $\V(R)$ (cf. \cite[Subsection 4A]{ara-hazrat-li-sims}).

Let $R$ now be a $\Gamma$-graded ring. Recall that a left $R$-module $M$ is called {\it $\Gamma$-graded} if there is a decomposition $M=\bigoplus_{\gamma\in\Gamma}M_\gamma$ such that $R_\alpha M_\gamma\subseteq M_{\alpha\gamma}$ for any $\alpha, \gamma\in \Gamma$. We denote by $R$-$\Gr$ the category of $\Gamma$-graded unital left $R$-modules with morphisms the $R$-module homomorphisms that preserve grading. Moreover, we denote by $R$-$\Gr_{\proj}$ the full subcategory of $R$-$\Gr$ whose objects are the projective objects of $R$-$\Gr$ that are finitely generated as a left $R$-module. If $R$ has graded local units, we define
\[\V^{\gr}(R)=\{[P] \mid P \in R\text{-}\Gr_{\proj}\}\]
where $[P]$ denotes the isomorphism class of $P$ as a graded left $R$-module. $\V^{\gr}(R)$ becomes an abelian monoid by defining $[P]+[Q]=[P\oplus Q]$. The graded Grothendieck group $K^{\gr}_0(R)$ is the group completion of $\V^{\gr}(R)$ (cf. \cite[Subsection 4A]{ara-hazrat-li-sims}).

Let $R$ be a $\Gamma$-graded ring. The {\it smash product} ring $R \# \Gamma$ is defined as the set of all formal sums $\sum_{\gamma\in\Gamma}r^{(\gamma)}p_{\gamma}$ where $r^{(\gamma)}\in R$ for any $\gamma\in\Gamma$, the $p_{\gamma}$'s are symbols, and all but finitely many coefficients $r^{(\gamma)}$ are zero. Addition is defined component-wise and multiplication is defined by linear extension of the rule $(rp_{\alpha})(sp_{\beta})=rs_{\alpha\beta^{-1}}p_{\beta}$ where $r,s\in R$ and $\alpha,\beta\in \Gamma$. We will use the proposition below to compute the graded $\V$-monoid of a weighted Leavitt path algebra.

\begin{proposition}[{\cite[Proposition 66]{Raimund2}}]\label{propsmash}
Let $R$ be a $\Gamma$-graded ring with graded local units. Then $R$-$\Gr_{\proj}\cong R\#\Gamma$-$\Mod_{\proj}$ by an isomorphism that commutes with direct sums. It follows that $\V^{\gr}(R)\cong \V(R\#\Gamma)$.
\end{proposition}

\subsection{Admissible weight maps}
In this subsection $(E,w)$ denotes a weighted graph. For any $v\in E_{\reg}^0$ we choose an edge $e^v\in s^{-1}(v)$ such that $w(e^v)=w(v)$. Recall that the unweighted graph associated to $(E,w)$ is denoted by $\hat E$.

\begin{definition}\label{defweight}
An {\it admissible weight map} for $(E,w)$ is a map $W:\hat E^1\rightarrow \Gamma$ which has the property that \[W(e_i)W(e_j)^{-1}=W(f_i)W(f_j)^{-1}\quad\text{ and }\quad W(e_i)^{-1}W(f_i)=W(e_j)^{-1}W(f_j)\]
for any $v\in E^0_{\reg}$, $e,f\in s^{-1}(v)$ and $1\leq i,j\leq \min\{w(e),w(f)\}$.
\end{definition}

The lemma below is straightforward to check.
\begin{lemma}\label{lemweight1}
A map $W:\hat E^1\rightarrow \Gamma$ is an admissible weight map for $(E,w)$ if and only if $W(e_i)=W(e_i^{s(e)})W(e_1^{s(e)})^{-1}W(e_1)$ for any $e\in E^1$ and $1\leq i\leq w(e)$.
\end{lemma}

Set $S:=\{e_1,e^v_i\mid e\in E, v\in E^0_{\reg}, 1\leq i\leq w(v)\}\subseteq \hat E^1$. It follows from Lemma \ref{lemweight1} that there is a $1$-$1$ correspondence between the set of all maps $S\rightarrow \Gamma$ and the set of all admissible weight maps $\hat E^1\rightarrow \Gamma$ for $(E,w)$.

\begin{lemma}\label{lemweight2}
Let $W$ be an admissible weight map for $(E,w)$. Then $W$ induces a $\Gamma$-grading on $L(E,w)$ such that $\deg(v)=\epsilon$, $\deg(e_i)=W(e_i)$ and $\deg(e_i^*)=W(e_i)^{-1}$ for any $v\in E^0$, $e\in E^1$ and $1\leq i\leq w(e)$.  
\end{lemma}
\begin{proof}
Let $K\X$ denote the free $K$-algebra on the set $X=\{v,e_i,e_i^*\mid v\in E^0, e\in E^1, 1\leq i\leq w(e)\}$. There is a $\Gamma$-grading on $K\X$ defined by $\deg(v)=\epsilon$, $\deg(e_i)=W(e_i)$ and $\deg(e_i^*)=W(e_i)^{-1}$ for any $v\in E^0$, $e\in E^1$ and $1\leq i\leq w(e)$.  Clearly the relations (i)-(iv) in Definition \ref{def3} are homogeneous with respect to this grading. Hence the $\Gamma$-grading on $K\X$ induces a $\Gamma$-grading on $L(E,w)$. 
\end{proof}

\begin{example}\label{exweight1}
Let $\lambda$ be defined as in the last paragraph of Section 3. Define a map $W:\hat E^1\rightarrow \Z^\lambda$ by $W(e_i)=\alpha_i$ where $\alpha_i$ denotes the element of $\mathbb Z^\lambda$ whose $i$-th component is $1$ and whose other components are $0$. Then $W$ is an admissible weight map for $(E,w)$. It induces the standard grading on $L(E,w)$.
\end{example}

\begin{example}\label{exweight2}
Let $\lambda$ be defined as in the previous example. For any $v\in E^0_{\reg}$ write $s^{-1}(v)=\{e^{v,1},\dots,e^{v,n_v}\}$. Set $\mu:=\sup\{n_v\mid v\in E^0_{\reg}\}$ if this supremum is finite and otherwise $\mu:=\omega$ where $\omega$ is the smallest infinite ordinal. Define a map $W:\hat E^1\rightarrow \Z^\lambda\oplus \Z^\mu$ by $W(e^{v,j}_i)=(\alpha_i,\beta_j)$ for any $v\in E^0_{\reg}$, $1\leq j\leq n_v$ and $1\leq i\leq w(e^{v,j})$ ($\alpha_i$ is defined as in the previous example and $\beta_j$ is defined analogously). Then $W$ is an admissible weight map for $(E,w)$ and therefore it induces a $\Z^\lambda\oplus \Z^\mu$-grading on $L(E,w)$. Obviously this grading is finer than the one defined in the previous example.
\end{example}

\subsection{Leavitt path algebras of bi-separated graphs}

In this subsection we recall some definitions and results from \cite{mohan-suhas}.

\begin{definition}\label{defbisep}
A {\it bi-separated graph} is a triple $\dot{E}=(E,C,D)$ such that  
\begin{enumerate}[(i)]
\item $E=(E^0,E^1,s,r)$ is a graph,
\smallskip
\item $C=\bigsqcup_{v\in E^0}C_v$ where $C_v$ is a partition of $s^{-1}(v)$ for any $v\in E^0$,
\smallskip
\item $D=\bigsqcup_{v\in E^0}D_v$ where $D_v$ is a partition of $r^{-1}(v)$ for any  $v\in E^0$,
\smallskip
\item $|X\cap Y|\leq 1$ for any $X\in C$ and $Y\in D$.
\end{enumerate}
\end{definition}

For a bi-separated graph $\dot{E}=(E,C,D)$ we set $C_{\fin}:=\{X\in C\mid |X|<\infty\}$ and $D_{\fin}:=\{Y\in D\mid |Y|<\infty\}$. In the following we assume that $C=C_{\fin}$ and $D=D_{\fin}$. Recall from Section 2 that if $E$ is a graph, then $E_d$ denotes the double graph of $E$ and $P(E)$ denotes the path algebra of $E$. In \cite{mohan-suhas} the double graph of $E$ was denoted by $\hat E$ and the path algebra of $E$ by $K(E)$.

Let $\dot{E}=(E,C,D)$ be a bi-separated graph. For $X\in C$ we denote by $s(X)$ the common source of the edges in $X$. For $Y\in D$ we denote by $r(Y)$ the common range of the edges in $Y$. Moreover, for $X\in C$ and $Y\in D$ we define
\[XY=YX=\begin{cases}
e,&\text{if }X\cap Y=\{e\},\\
0,&\text{if }X\cap Y=\emptyset.
\end{cases}\]

\begin{definition}\label{defbslpa}
Let $\dot{E}=(E,C,D)$ be a bi-separated graph. The {\it Leavitt path algebra of $\dot E$}
with coefficients in $K$, denoted by $L(\dot E)$, is the quotient of $P(E_d)$ obtained by imposing the following relations:
\begin{enumerate}[(L1)]
\item $\sum_{Y\in D}(XY)(YX')^*=\delta_{XX'}s(X)$ for any $X,X'\in C$,
\smallskip
\item $\sum_{X\in C}(YX)^*(XY')=\delta_{YY'}r(Y)$ for any $Y,Y'\in D$.
\end{enumerate}
\end{definition}

\begin{example}
Let $(E,w)$ be a weighted graph. For any $v\in E^0$ and $1\leq i\leq w(v)$ define $X^i_v:=\{e_i\mid e\in s^{-1}(v),w(e)\geq i\}$. For any $e\in E^1$ define $Y^e:=\{e_i\mid 1\leq i\leq w(e)\}$. Moreover, define $C_v:=\{X^i_v\mid 1\leq i\leq w(v)\}$ and $D_v:=\{Y^e \mid e\in r^{-1}(v)\}$ for any $v\in E^0$. Then $(\hat E, C, D)$ is a bi-separated graph and $L(E,w)\cong L(\hat E, C, D)$.
\end{example}

Let $\dot{E}=(E,C,D)$ be a bi-separated graph. We define a equivalence relation $\sim_D$ on $C$ as follows: For $X,X'\in C$ define $X\sim_D X'$ if $X=X'$ or there exists a finite sequence $X_0, Y_1, X_1, Y_2, X_2,\dots, Y_n, X_n$ such that $X_i\in C~(0\leq i\leq n)$, $Y_i\in D~(1\leq i\leq n)$, $X_0=X$, $X_n=X'$, $X_{i-1}\cap Y_i\neq \emptyset~(1\leq i\leq n)$ and $Y_i\cap X_i\neq \emptyset~(1\leq i\leq n)$. Let $C=\bigsqcup_{\lambda\in\Lambda}\cX_\lambda$ be the partition of $C$ induced by $\sim_D$. Similarly we define an equivalence relation $\sim_C$ on $D$ and let $D=\bigsqcup_{\lambda'\in\Lambda'}\cY_{\lambda'}$ be the partition of $D$ induced by $\sim_C$. There is a canonical bijective map $\Lambda\to\Lambda'$, see \cite[\S 4.2]{mohan-suhas}. Therefore we will denote the indexing sets of the partitions of both $C$ and $D$ by $\Lambda$. %Moreover, we may assume that if $X\in C$ and $Y\in D$ such that $X\cap Y\neq \emptyset$, then there is a $\lambda\in\Lambda$ such that $X\in\cX_\lambda$ and $Y\in\cY_\lambda$.
The bi-separated graph $\dot{E}$ is called {\it tame}, if $|\cX_\lambda|,|\cY_\lambda|<\infty$ for any $\lambda\in\Lambda$. $\dot{E}$ is called {\it docile} if it is tame and for any $\lambda\in\Lambda$ there are distinguished elements $X_\lambda\in\cX_\lambda$ and $Y_\lambda\in\cY_\lambda$ such that $X_\lambda\cap Y\neq \emptyset$ for any $Y\in\cY_\lambda$ and $X\cap Y_\lambda\neq \emptyset$ for any $X\in\cX_\lambda$.

Let $\dot{E}=(E,C,D)$ be a docile bi-separated graph. The words $(XY_\lambda)(Y_\lambda X')^*~(\lambda\in\Lambda,~X,X'\in \cX_\lambda)$ and $(YX_\lambda)^*(X_\lambda Y')~(\lambda\in\Lambda,~Y,Y'\in \cY_\lambda)$ are called {\it forbidden}. We call a path in the double graph $E_d$ a {\it d-path}. A {\it normal d-path} or {\it nod-path} is a d-path such that none of its subwords is forbidden.

\begin{theorem}[{\cite[Theorem 5.5]{mohan-suhas}}] \label{thmbsbasis}
Let $\dot{E}$ be a docile bi-separated graph. Then the nod-paths form a linear basis for $L(\dot E)$.
\end{theorem} 

\subsection{The covering bi-separated graph of a weighted graph defined by an admissible weight map}

\begin{definition} 
Let $(E,w)$ be a weighted graph and $W$ an admissible weight map for $(E,w)$. Define a graph $F$ by 
\begin{align*}
F^0&=\{v^{(\gamma)}\mid v\in E^0,\gamma\in\Gamma\},\\
F^1&=\{e_i^{(\gamma)}\mid e\in E^1,1\leq i\leq w(e),\gamma\in\Gamma\},\\
s_F(e_i^{(\gamma)})&=s(e)^{(W(e_i^{s(e)})W(e_1^{s(e)})^{-1}\gamma)},\\
r_F(e_i^{(\gamma)})&=r(e)^{(W(e_1)^{-1}\gamma)}.
\end{align*}
For any $\gamma\in\Gamma$, $v\in E^0$ and $1\leq i\leq w(v)$ set $X_{\gamma,v,i}=\{e_i^{(\gamma)}\mid e\in s^{-1}(v),w(e)\geq i\}$. For any $\gamma\in\Gamma$ and $e\in E^1$ set $Y_{\gamma,e}=\{e_i^{(\gamma)}\mid 1\leq i\leq w(e)\}$. For any $v^{(\gamma)}\in F^0$ define
\begin{align*}
C_{v^{(\gamma)}}&=\{X_{W(e_1^v)W(e_i^v)^{-1}\gamma,v,i}\mid 1\leq i\leq w(v)\},\\
D_{v^{(\gamma)}}&=\{Y_{W(e_1)\gamma,e}\mid e\in r^{-1}(v)\}.
\end{align*}
The bi-separated graph $\dot F=(F,C,D)$ is called the {\it covering bi-separated graph} of $(E,w)$ defined by $W$.
\end{definition}

Until the end of this subsection $(E,w)$ denotes a weighted graph, $W$ an admissible weight map for $(E,w)$ and $\dot F=(F,C,D)$ the covering bi-separated graph of $(E,w)$ defined by $W$. Recall that $W$ induces a $\Gamma$-grading on $L(E,w)$. We will show that $L(\dot F)\cong L(E,w)\#\Gamma$. 

The equivalence relations $\sim_D$ on $C$ and $\sim_C$ on $D$ defined in the previous subsection are given by $X_{\beta,u,i}\sim_D X_{\gamma,v,j} \Leftrightarrow (\beta=\gamma\land u=v)$ and $Y_{\beta,e}\sim_C Y_{\gamma,f} \Leftrightarrow (\beta=\gamma\land s(e)=s(f))$, respectively. Hence the partition of $C$ induced by $\sim_D$ is 
\[C=\bigsqcup_{(\gamma,v)\in \Gamma\times E^0_{\reg}}\cX_{(\gamma,v)}\text{ where }\cX_{(\gamma,v)}=\{X_{\gamma,v,i}\mid 1\leq i\leq w(v)\}.\]
Similarly the partition of $D$ induced by $\sim_C$ is 
\[D=\bigsqcup_{(\gamma,v)\in \Gamma\times E^0_{\reg}}\cY_{(\gamma,v)}\text{ where }\cY_{(\gamma,v)}=\{Y_{\gamma,e}\mid e\in s^{-1}(v)\}.\]
Moreover, $\dot F$ is docile. As distinguished elements we choose $X_{(\gamma,v)}:=X_{\gamma,v,1}$ and $Y_{(\gamma,v)}:=Y_{\gamma,e^v}$ for any $(\gamma,v)\in \Gamma\times E^0_{\reg}$.

\begin{lemma}\label{lemcov1}
There is a surjective $K$-algebra homomorphism $\psi:L(\dot F)\to L(E,w)\#\Gamma$ such that 
\begin{equation}
\psi(v^{(\gamma)})=vp_{\gamma}, \quad\psi(e_i^{(\gamma)})=e_ip_{W(e_1)^{-1}\gamma},\quad\psi((e_i^{(\gamma)})^*)=e_i^*p_{W(e_i^{s(e)})W(e_1^{s(e)})^{-1}\gamma}
\end{equation}
for any $v^{(\gamma)}\in F^0$ and $e_i^{(\gamma)}\in F^1$.
\end{lemma}
\begin{proof}
In order to show that there is a $K$-algebra homomorphism $\psi:L(\dot F)\to L(E,w)\#\Gamma$ such that (11) holds, it suffices to prove the relations (i)-(iv) below, where \[A_{v^{(\gamma)}}=vp_{\gamma},\quad A_{e_i^{(\gamma)}}=e_ip_{W(e_1)^{-1}\gamma}, \quad B_{e_i^{(\gamma)}}=e_i^*p_{W(e_i^{s(e)})W(e_1^{s(e)})^{-1}\gamma}\]
for any $v^{(\gamma)}\in F^0$ and $e_i^{(\gamma)}\in F^1$, and moreover $A_0=0$ and $B_0=0$ (cf. \cite[Proposition 3.18]{mohan-suhas}).
\begin{enumerate}[(i)]
\item $A_{u}A_{v}=\delta_{uv}A_{v}$ for any $u,v\in F^0$.
\smallskip
\item $A_{s(f)}A_{f}=A_{f}A_{r(f)}=A_{f}$ and $A_{r(f)}B_{f}=B_{f}A_{s(f)}=B_{f}$ for any $f\in F^1$.
\smallskip
\item $\sum_{Y\in D}A_{XY}B_{YX'}=\delta_{XX'}A_{s(X)}$ for any $X,X'\in C$.
\smallskip
\item $\sum_{X\in C}B_{YX}A_{XY'}=\delta_{YY'}A_{r(Y)}$ for any $Y,Y'\in D$.
\end{enumerate}
We leave (i) and (ii) to the reader and show only (iii) and (iv). 

First we prove (iii). Let $X,X'\in C$. Clearly we may assume that $X,X'\in \cX_{(\gamma,v)}$ for some $(\gamma,v)\in \Gamma\times E^0_{\reg}$ (otherwise (iii) is trivially satisfied). It follows that $X=X_{\gamma,v,i}$ and $X'=X_{\gamma,v,j}$ for some $1\leq i,j\leq w(v)$. Clearly
\begin{align*}
&\sum_{Y\in D}A_{XY}B_{YX'}\\
=&\sum_{e\in s^{-1}(v),w(e)\geq i,j}A_{e^{(\gamma)}_i}B_{e^{(\gamma)}_j}\\
=&\sum_{e\in s^{-1}(v),w(e)\geq i,j} e_ip_{W(e_1)^{-1}\gamma}e_j^*p_{W(e_j^{v})W(e_1^{v})^{-1}\gamma}\\
=&\sum_{e\in s^{-1}(v),w(e)\geq i,j} e_i(e_j^*)_{W(e_1)^{-1}W(e_1^{v})W(e_j^{v})^{-1}}p_{W(e_j^{v})W(e_1^{v})^{-1}\gamma}\\
=&\sum_{e\in s^{-1}(v),w(e)\geq i,j} e_ie_j^*p_{W(e_j^{v})W(e_1^{v})^{-1}\gamma}\\
=&\delta_{ij}vp_{W(e_j^{v})W(e_1^{v})^{-1}\gamma}\\
=&\delta_{XX'}A_{s(X)}
\end{align*}
by Lemma \ref{lemweight1} and hence (iii) holds. 

Next we prove (iv). Let $Y,Y'\in D$. Clearly we may assume that $Y,Y'\in \cY_{(\gamma,v)}$ for some $(\gamma,v)\in \Gamma\times E^0_{\reg}$ (otherwise (iv) is trivially satisfied). It follows that $Y=Y_{\gamma,e}$ and $Y'=Y_{\gamma,f}$ for some $e,f\in s^{-1}(v)$. Clearly
\begin{align*}
&\sum_{X\in C}B_{YX}A_{XY'}\hspace{6.4cm}\\
=&\sum_{1\leq i\leq w(e),w(f)}B_{e^{(\gamma)}_i}A_{f^{(\gamma)}_i}
\end{align*}
\begin{align*}
=&\sum_{1\leq i\leq w(e),w(f)}e_i^*p_{W(e_i^{v})W(e_1^{v})^{-1}\gamma}f_ip_{W(f_1)^{-1}\gamma}\hspace{2cm}\\
=&\sum_{1\leq i\leq w(e),w(f)} e_i^*(f_i)_{W(e_i^{v})W(e_1^{v})^{-1}W(f_1)}p_{W(f_1)^{-1}\gamma}\\
=&\sum_{1\leq i\leq w(e),w(f)} e_i^*f_ip_{W(f_1)^{-1}\gamma}\\
=&\delta_{ef}r(e)p_{W(f_1)^{-1}\gamma}\\
=&\delta_{YY'}A_{r(Y)}
\end{align*}
by Lemma \ref{lemweight1} and hence (iv) holds.

Thus there is a $K$-algebra homomorphism $\psi:L(\dot F)\to L(E,w)\#\Gamma$ such that (11) holds. Clearly the image of $\psi$ contains the set $S:=\{v p_\gamma, e_ip_\gamma, e_i^*p_\gamma\mid v\in E^0, e\in E^1, 1\leq i\leq w(e), \gamma\in \Gamma\}$. But $S$ generates $L(E,w)\# \Gamma$ as a $K$-algebra and therefore $\psi$ is surjective.
\end{proof}

Our next goal is to show that the homomorphism $\psi:L(\dot F)\to L(E,w)\#\Gamma$ from Lemma \ref{lemcov1} is injective. Let $A$ denote the set of all nod-paths for $\dot F$ and $B$ the set of all nod-paths for $(E,w)$. Set $Y:=\{v^{(\gamma)},e_{i}^{(\gamma)},(e^{(\gamma)}_{i})^*\mid \gamma\in \Gamma, v\in E^0, e\in E^1, 1\leq i\leq w(e)\}$ and $Z:=\{v,e_i,e_i^*\mid v\in E^0, e\in E^1, 1\leq i\leq w(e)\}$. Let $\langle Y\rangle$ be the set of all finite nonempty words over $Y$ and $\langle Z \rangle$ the set of all finite nonempty words over $Z$. With juxtaposition of words, $\langle Y\rangle$ and $\langle Z\rangle$ become semigroups. Let $\xi:\langle Y\rangle\rightarrow \langle Z\rangle$ be the unique semigroup homomorphism such that $\xi(v^{(\gamma)})=v$, $\xi(e_i^{(\gamma)})=e_i$ and $\xi((e_i^{(\gamma)})^*)=e_i^*$. For any $a\in A$ set $b(a):=\xi(a)$.

\begin{lemma}\label{lemcov3}
For any $a\in A$ there is a uniquely determined $\gamma(a)\in \Gamma$ such that $\psi(a)=b(a)p_{\gamma(a)}$. 
\end{lemma}
\begin{proof}
The assertion of the lemma clearly holds if $|a|\leq 1$. Suppose now that $a=y_1\dots y_n$ where $n\geq 2$ and $y_1,\dots, y_n\in Y\setminus F^0$. Then for any $1\leq i\leq n$ there is a $\gamma_i$ such that $\psi(y_i)=b(y_i)p_{\gamma_i}$. 

First we show that 
\begin{equation}
\deg{b(y_k)}=\gamma_{k-1}\gamma_k^{-1}~(2\leq k\leq n)
\end{equation}
where $\deg{b(y_k)}$ denotes the homogeneous degree of $b(y_k)$. We consider only the case that $y_{k-1}=e_i^{(\beta)}$ and $y_k=f_j^{(\gamma)}$ for some $e_i^{(\beta)},f_j^{(\gamma)}\in F^1$ and leave the other cases to the reader. Clearly $\gamma_{k-1}=W(e_1)^{-1}\beta$ and $\gamma_{k}=W(f_1)^{-1}\gamma$. Since $r(e)^{(W(e_1)^{-1}\beta)}=r_F(e_i^{(\beta)})=s_F(f_j^{(\gamma)})=s(f)^{(W(e^{s(f)}_j)W(e^{s(f)}_1)^{-1}\gamma)}$, we have 
\begin{equation*}
\gamma_{k-1}\gamma_k^{-1}=W(e_1)^{-1}\beta\gamma^{-1}W(f_1)=W(e^{s(f)}_j)W(e^{s(f)}_1)^{-1}W(f_1)=W(f_j)=\deg b(y_k)
\end{equation*}
in view of Lemma \ref{lemweight1}. Thus (12) holds.

Now we show by induction on $t$ that 
\begin{equation}
\psi(y_1\dots y_t)=b(y_1\dots y_t)p_{\gamma_t}~(2\leq t\leq n).
\end{equation}
\\
\underline{$t=2$:} Clearly $
\psi(y_1y_2)=b(y_1)p_{\gamma_1}b(y_2)p_{\gamma_2}=b(y_1)(b(y_2))_{\gamma_1\gamma_2^{-1}}p_{\gamma_2}=b(y_1y_2)p_{\gamma_2}$
by (12).\\
\\
\underline{$t-1\to t$:} 
Clearly 
\begin{align*}
&\psi(y_1\dots y_{t-1}y_{t})\\
=&\psi(y_1\dots y_{t-1})\psi(y_{t})\\
=&b(y_1\dots y_{t-1})p_{\gamma_{t-1}}b(y_{t})p_{\gamma_{t}}\\
=&b(y_1\dots y_{t-1})(b(y_{t}))_{\gamma_{t-1}\gamma_{t}^{-1}}p_{\gamma_{t}}=b(y_1\dots y_t)p_{\gamma_t}
\end{align*}
by the induction assumption and (12).

Hence (13) holds and thus in particular $\psi(a)=b(a)p_{\gamma(a)}$ where $\gamma(a)=\gamma_n$. The uniqueness of $\gamma(a)$ follows from the definition of $L(E,w)\# \Gamma$.
\end{proof}

\begin{lemma}\label{lemcov2}
$b(a)\in B$ for any $a\in A$.
\end{lemma}
\begin{proof}
We leave it to the reader to check that $b(a)$ is a d-path for $(E,w)$. We show now that $b(a)$ contains no forbidden subword and hence lies in $B$. That is clear if $a$ is a vertex, so let us assume that $a=y_1\dots y_n$ where $y_1, \dots, y_n\in Y\setminus F^0$. It is straightforward to check that the forbidden words for $\dot F$ are
\[(e^v)_i^{(\gamma)}((e^v)_j^{(\gamma)})^*~(\gamma\in\Gamma, v\in E^0_{\reg}, 1\leq i,j\leq w(v))\quad\text{ and }(e_1^{(\gamma)})^*f_1^{(\gamma)}~(\gamma\in\Gamma, v\in E^0_{\reg}, e,f\in s^{-1}(v)).\]
Assume that there is a $t\in \{1,\dots,n-1\}$ such that $\xi(y_ty_{t+1})=e^v_i(e^v_j)^*$ for some $v\in E^0$ and $1\leq i,j\leq w(v)$. Then $y_ty_{t+1}=(e^v)_i^{(\beta)}((e^v)_j^{(\gamma)})^*$ for some $\beta,\gamma\in \Gamma$. Since $a$ is a d-path, it follows that $\beta=\gamma$ (by the definition of $r_F$). Hence we get the contradiction that $a$ contains a forbidden subword. Analogously one can show that $b(a)$ contains no forbidden subword of the form $e_1^*f_1$ where $e,f\in s^{-1}(v)$ for some $v\in E^0$. Thus $b(a)\in B$ for any $a\in A$.
\end{proof}

\begin{lemma}\label{lemcov4}
The map $(b,\gamma):A\rightarrow B\times \Gamma$, $a\mapsto (b(a), \gamma(a))$ is injective.
\end{lemma}
\begin{proof}
Let $a=y_n\dots y_1, a'=y'_{n'}\dots y'_{1}\in A$ such that $b(a)=b(a')$ and $\gamma(a)=\gamma(a')$. Then clearly $n=n'$. We show by induction on $t$ that $y_t=y'_t$ for any $1\leq t\leq n$ (and hence $a=a'$).\\
\\
\underline{$t=1$:} Suppose that $y_1=e_i^{(\gamma)}$ for some $e_i^{(\gamma)}\in F^1$. Then $y'_1=e_i^{(\gamma')}$ for some $\gamma'\in \Gamma$ since $b(a)=b'(a)$. Clearly $W(e_1)^{-1}\gamma=\gamma(a)=\gamma(a')=W(e_1)^{-1}\gamma'$. Hence $\gamma=\gamma'$ and therefore $y_1=y'_1$. The cases $y_1=(e_i^{(\gamma)})^*$ and $y_1=v^{(\gamma)}$ are similar.\\
\\
\underline{$t-1\to t$:} Assume that $y_{t-1}=y'_{t-1}$ for some $2\leq t\leq n$. Suppose that $y_{t}=e_i^{(\gamma)}$ for some $e_i^{(\gamma)}\in F^1$. Then $y'_{t}=e_i^{(\gamma')}$ for some $\gamma'\in \Gamma$ since $b(a)=b'(a)$. Since $y_{t-1}=y'_{t-1}$ and $a,b$ are d-paths, we have $r_F(e_i^{(\gamma)})=r_F(e_i^{(\gamma')})$. It follows from the definition of $r_F$ that $\gamma=\gamma'$ and therefore $y_{t}=y'_{t}$. The case $y_{t}=(e_i^{(\gamma)})^*$ is similar.
\end{proof}

\begin{proposition}\label{propcov}
Let $(E,w)$ be a weighted graph, $W$ an admissible weight map for $(E,w)$ and $\dot F$ the covering bi-separated graph of $(E,w)$ defined by $W$. Then $L(\dot F)\cong L(E,w)\# \Gamma$.
\end{proposition}
\begin{proof}
Let $\psi:L(\dot F)\to L(E,w)\#\Gamma$ be the surjective $K$-algebra homomorphism from Lemma \ref{lemcov1}. Let $A, B$ and $b(a)~(a\in A)$ be defined as in the first paragraph after the proof of Lemma \ref{lemcov1}. By Lemma \ref{lemcov3} there is for any $a\in A$  a uniquely determined $\gamma(a)$ such that $\psi(a)=b(a)p_{\gamma(a)}$. By Lemma \ref{lemcov2}, $b(a)\in B$ for any $a\in A$. By Lemma \ref{lemcov4} the map $(b,\gamma):A\to B\times \Gamma$ is injective. 

In order to prove the proposition it suffices to show that $\psi$ is injective. Let $x\in L(\dot F)$. By Theorem \ref{thmbsbasis} we can write $x=\sum_{a\in A}k_{a}a$ where almost all coefficients $k_a\in K$ are zero. Clearly 
\[\psi(x)=\sum_{a\in A} k_{a}b(a) p_{\gamma(a)}=\sum_{\beta\in \Gamma}(\sum_{a\in A,\gamma(a)=\beta} k_{a}b(a))p_{\beta}.\] 
Assume now that $\psi(x)=0$. Then $\sum_{a\in A,\gamma(a)=\beta} k_{a}b(a)=0$ in $L(E,w)$ for any $\beta\in \Gamma$. But since the map $(b,\gamma):A\to B\times \Gamma$ is injective, we have $b(a)\neq b(a')$ for any $a\neq a'\in A$ such that $\gamma(a)=\gamma(a')$. It follows from Theorem \ref{thmbasis} that $k_a=0$ for any $a\in A$ and hence $x=0$. Thus $\psi$ is injective.
\end{proof}

\subsection{The monoid $\V^{\gr}(L(E,w))$}
In this subsection we fix a weighted graph $(E,w)$ and an admissible weight map $W:\hat E^1\rightarrow \Gamma$ for $(E,w)$. Recall that $W$ induces a $\Gamma$-grading on $L(E,w)$.

\begin{definition}\label{defgrmonoid}
For any $v\in E^0$ write $s^{-1}(v)=\{e^{v,1},\dots,e^{v,n(v)}\}$ where $w(e^{v,1})\leq\dots \leq w(e^{v,n(v)})$. Moreover, write $w(s^{-1}(v))=\{w_1(v),\dots,w_{k_v}(v)\}$ where $w_1(v)<\dots <w_{k_v}(v)$. Set $n_l(v):=\#\{e\in s^{-1}(v)\mid w(e)\leq w_i(v)\}$ for any $1\leq l\leq k_v$. We define $\cM^{\gr}(E,w)$ as the abelian monoid presented by the generating set $\{v^{(\gamma)}, I_1^{v,\gamma},\dots,I_{k_v-1}^{v,\gamma}\mid v\in E^0, \gamma\in \Gamma\}$ and the relations 
\[I^{v,\gamma}_{l-1}+\sum_{i=w_{l-1}(v)+1}^{w_l(v)}v^{(W(e_i^v)W(e_1^v)^{-1}\gamma)}=I^{v,\gamma}_{l}+\sum_{j=n_{l-1}(v)+1}^{n_l(v)}r(e^{v,j})^{(W(e^{v,j}_1)^{-1}\gamma)}\quad (v\in E^0, \gamma\in\Gamma,1\leq l\leq k_v).\]
Here we use the convention $I^{v,\gamma}_0=I^{v,\gamma}_{k_v}=w_0(v)=n_0(v)=0$ for any $v\in E^0$ and $\gamma\in\Gamma$. 
\end{definition}

\begin{theorem}\label{thmgrmonoid}
$\V^{\gr}(L(E,w))\cong \cM^{\gr}(E,w)$.
\end{theorem}
\begin{proof}
Let $\dot F=(F,C,D)$ be the covering bi-separated graph of $(E,w)$ defined by $W$. By Propositions \ref{propsmash} and \ref{propcov} we have $\V^{\gr}(L(E,w))\cong \V(L(E,w)\#\Gamma)\cong \V(L(\dot F))$. Hence it suffices to show that $\V(L(\dot F))\cong \cM^{\gr}(E,w)$. We use the notation of Definition \ref{defgrmonoid}. Recall that 
\[C=\bigsqcup_{(\gamma,v)\in \Gamma\times E^0_{\reg}}\cX_{(\gamma,v)}\text{ where }\cX_{(\gamma,v)}=\{X_{\gamma,v,i}\mid 1\leq i\leq w(v)\}\]
and
\[D=\bigsqcup_{(\gamma,v)\in \Gamma\times E^0_{\reg}}\cY_{(\gamma,v)}\text{ where }\cY_{(\gamma,v)}=\{Y_{\gamma,e}\mid e\in s^{-1}(v)\}\]
where $X_{\gamma,v,i}=\{e_i^{(\gamma)}\mid e\in s^{-1}(v),w(e)\geq i\}$ and $Y_{\gamma,e}=\{e_i^{(\gamma)}\mid 1\leq i\leq w(e)\}$. Moreover, $\dot F$ is docile. For any $(\gamma,v)\in \Gamma\times E^0_{\reg}$ let $A^{\gamma,v}$ be the matrix $w(v)\times n(v)$-matrix whose entry at position $(i,j)$ is 
\[X_{\gamma,v,i}\cap Y_{\gamma,e^{v,j}}=\begin{cases}(e^{v,j})_i^{(\gamma)},&\text{ if }w(e^{v,j})\geq i,\\
0,&\text{ otherwise.}\end{cases}\]
Then each $A^{\gamma,v}$ has the upper triangular block form
\[A^{\gamma,v}=\begin{pmatrix}
A^{\gamma,v}_{11}&A^{\gamma,v}_{12}&A^{\gamma,v}_{13}&\dots&A^{\gamma,v}_{1k_v}\\0&A^{\gamma,v}_{22}&A^{\gamma,v}_{23}&\dots&A^{\gamma,v}_{2k_v}\\0&0&A^{\gamma,v}_{33}&\dots&A^{\gamma,v}_{3k_v}\\\vdots&\vdots&\ddots&\ddots&\vdots\\0&0&\dots&0&A^{\gamma,v}_{k_vk_v}
\end{pmatrix}\]
where each $A^{\gamma,v}_{xy}$ is a matrix having $w_x(v)-w_{x-1}(v)$ rows, $n_y(v)-n_{y-1}(v)$ columns and no zero entry. It follows from \cite[Theorem 6.1 and preceding paragraph]{mohan-suhas} that $\V(L(\dot F))$ is presented by the generating set $\{v^{(\gamma)}, I_1^{v,\gamma},\dots,I_{k_v-1}^{v,\gamma}\mid v\in E^0, \gamma\in \Gamma\}$ and the relations 
\[I^{v,\gamma}_{l-1}+\sum_{i=w_{l-1}(v)+1}^{w_l(v)}s_F(X_{\gamma,v,i})=I^{v,\gamma}_{l}+\sum_{j=n_{l-1}(v)+1}^{n_l(v)}r_F(Y_{\gamma,e^{v,j}})\quad (v\in E^0, \gamma\in\Gamma,1\leq l\leq k_v).\]
It follows from the definition of $s_F$ and $r_F$ that $\V(L(\dot F))\cong \cM^{\gr}(E,w)$.
\end{proof}

\subsection{The graded Grothendieck group $K_0^{\gr}(L(E,w))$}
\begin{theorem}\label{thmkzerogr}
Let $(E,w)$ be a weighted graph and $W$ an admissible weight map for $(E,w)$. Then the abelian group $K_0^{\gr}(L(E,w))$ is presented by the generating set $\{v^{(\gamma)}\mid v\in E^0, \gamma\in \Gamma\}$ and the relations 
\[\sum_{i=1}^{w(v)}v^{(W(e_i^v)W(e_1^v)^{-1}\gamma)}=\sum_{e\in s^{-1}(v)}r(e)^{(W(e_1)^{-1}\gamma)}\quad (v\in E^0, \gamma\in\Gamma).\]
\end{theorem}
\begin{proof}
The theorem follows from Theorem \ref{thmgrmonoid} and the fact that $K_0^{gr}(L(E,w))$ is the group completion of the abelian monoid $\V^{\gr}(L(E,w)$ (cf. the proof of Theorem \ref{thmkzero}).
\end{proof}

\begin{example}
Consider again the weighted graphs\\
\[
(E,w):\xymatrix@C+15pt{ u& v\ar[l]_{e,1}\ar[r]^{f,2}& x}\quad\text{ and }\quad(E,w'):\xymatrix@C+15pt{ u& v\ar[l]_{e,2}\ar[r]^{f,2}& x}
\]
from Example \ref{exkzero}. Let $W$ and $W'$ be the admissible weight maps for $(E,w)$ and $(E,w')$ that induce the standard $\Z^2$-gradings on $L(E,w)$ and $L(E,w')$, respectively (see Example \ref{exweight1}).
By Theorem \ref{thmgrmonoid}, $\V(L(E,w))$ is presented by the generating set $\{u^{(m,n)},v^{(m,n)},x^{(m,n)},I^{(m,n)}\mid (m,n)\in\Z^2\}$ and the relations
\[v^{(m,n)}=I^{(m,n)}+u^{(m-1,n)},\quad I^{(m,n)}+v^{(m-1,n+1)}=x^{(m-1,n)}\quad((m,n)\in \Z^2).\]
Hence $\V(L(E,w))$ is the free abelian monoid generated by a countably infinite set. On the other hand, Theorem \ref{thmgrmonoid} implies that $\V(L(E,w'))$ is presented by the generating set $\{u^{(m,n)},v^{(m,n)},x^{(m,n)}\mid (m,n)\in\Z^2\}$ and the relations
\[v^{(m,n)}+v^{(m-1,n+1)}=u^{(m-1,n)}+x^{(m-1,n)}\quad((m,n)\in \Z^2).\]
%It follows that $\V(L(E,w))\not\cong \V(L(E,w'))$ since $\V(L(E,w))$ is a refinement monoid but $\V(L(E,w'))$ is not. 
Theorem \ref{thmkzerogr} implies that both graded Grothendieck groups $K^{\gr}_0(L(E,w))$ and $K^{\gr}_0(L(E,w'))$ are presented by presented by the generating set $\{u^{(m,n)},v^{(m,n)},x^{(m,n)}\mid (m,n)\in\Z^2\}$ and the relations
\[v^{(m,n)}+v^{(m-1,n+1)}=u^{(m-1,n)}+x^{(m-1,n)}\quad((m,n)\in \Z^2).\]
Hence $K^{\gr}_0(L(E,w))\cong K^{\gr}_0(L(E,w'))$ is the free abelian group generated by a countably infinite set.
\end{example}

%\section{Ideal structure}
\section{Representations}
In this section $(E,w)$ denotes a fixed weighted graph. Recall that $\hat E$ denotes the unweighted graph associated to $(E,w)$. If $e_i\in \hat E^1$, then we call $\tg(e_i):=i$ the {\it tag} of $e_i$ and $\mot(e_i):=e$ is the {\it structure edge} of $e_i$. If $F$ is a graph, then we denote by $\Path(F)$ the set of all paths in $F$. Moreover, if $u,v\in F^0$, then we denote by $_u\!\Path(F)$ the set of all paths starting in $u$, by $\Path_v(F)$ the set of all paths ending in $v$ and by $_u\!\Path_v(F)$ the intersection of $_u\!\Path(F)$ and $\Path_v(F)$. 

\subsection{Representation graphs for weighted graphs}

\subsubsection{Representation graphs}

\begin{definition}\label{defwp}
A {\it representation graph for $(E,w)$} is a pair $(F,\phi)$ where $F=(F^0,F^1,s_F,r_F)$ is a graph and $\phi=(\phi^0,\phi^1):F\rightarrow \hat E$ a graph homomorphism such that (i) and (ii) below are satisfied.
\begin{enumerate}[(i)]
\item For any $v\in F^0$ and $1\leq i\leq w(\phi^0(v))$ there is precisely one $f\in s_F^{-1}(v)$ such that $\tg(\phi^1(f))=~i$.

\smallskip

\item For any $v\in F^0$ and $e\in r^{-1}(\phi^0(v))$ there is precisely one $f\in r_F^{-1}(v)$ such that $\mot(\phi^1(f))=e$.
\end{enumerate}
\end{definition}

Usually we visualise a representation graph $(F,\phi)$ by drawing the graph $F$ and labelling each vertex $v\in F^0$ with $\phi^0(v)$ and each edge $f\in F^1$ with $\phi^1(f)$, see the four examples below.

\begin{example}\label{exrepgr1}
Suppose that $(E,w)$ is the weighted graph $\xymatrix{ v\ar@(dr,ur)_{f,1}\ar@(dl,ul)^{e,1}}$. Then a representation graph $(F,\phi)$ for $(E,w)$ is given by
\begin{equation*}
\xymatrix@C=0.3cm{&\ar@{.>}[rrrr]^{e_1}&&&&v\ar[rrrr]^{e_1}&&&&v\ar[rrrr]^{e_1}&&&&v\ar@{.>}[rrrr]^{e_1}&&&&\\
(F,\phi):\,\,&&&&&v\ar[u]^{f_1}&&&&v\ar[u]^{f_1}&&&&v.\ar[u]^{f_1}&&&&\\
&&&&\ar@{.>}[ur]^{e_1}&&\ar@{.>}[ul]_{f_1}&&\ar@{.>}[ur]^{e_1}&&\ar@{.>}[ul]_{f_1}&&\ar@{.>}[ur]^{e_1}&&\ar@{.>}[ul]_{f_1}&&&}
\end{equation*}
\end{example}

\begin{example}\label{exrepgr2}
Suppose again that $(E,w)$ is the weighted graph $\xymatrix{ v\ar@(dr,ur)_{f,1}\ar@(dl,ul)^{e,1}}$. Then a representation graph $(F,\phi)$ for $(E,w)$ is given by
\begin{equation*}
\xymatrix@C=0.5cm{&&v\ar@(ul,ur)^{e_1}&\\
(F,\phi):\,\,&&v.\ar[u]^{f_1}&\\
&\ar@{.>}[ur]^{e_1}&&\ar@{.>}[ul]_{f_1}}
\end{equation*}
\end{example}

\begin{example}\label{exrepgr3}
Suppose that $(E,w)$ is the weighted graph $\xymatrix{ v\ar@(dr,ur)_{f,2}\ar@(dl,ul)^{e,2}}$. Then a representation graph $(F,\phi)$ for $(E,w)$ is given by
\begin{equation*}
\xymatrix{(F,\phi):\,\,&v \ar@(dl,ul)^{e_1} \ar@(dr,ur)_{f_2}&.}
\end{equation*}
\end{example}

\begin{example}\label{exrepgr4}
Suppose that $(E,w)$ is the weighted graph $\xymatrix{ v\ar@(dr,ur)_{f,3}\ar@(dl,ul)^{e,3}}$. Then a representation graph $(F,\phi)$ for $(E,w)$ is given by
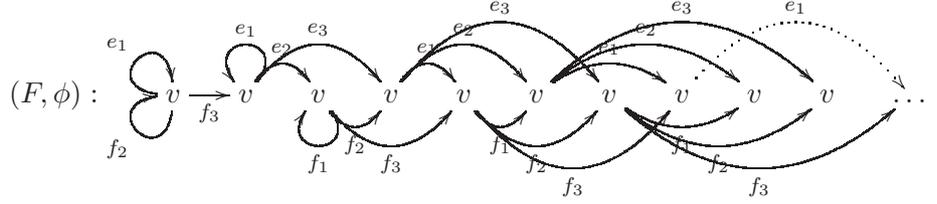
\begin{equation*}
  \xymatrix@=15pt{ (F,\phi): \,\,&v \ar@(l,u)^{e_1}  \ar@(d,l)^{f_2}  \ar[r]_{f_3}  & v \ar@(ur,ul)_{e_1}    \ar@/^1pc/[r]^{e_2}  \ar@/^1.6pc/[rr]^{e_3} & v\ar@(dr,dl)^{f_1}    \ar@/_1pc/[r]_{f_2}  
 \ar@/_1.6pc/[rr]_{f_3} & v  \ar@/^1pc/[r]^{e_1}  \ar@/^1.6pc/[rr]^{e_2}  \ar@/^2.3pc/[rrr]^{e_3} & 
 v  \ar@/_1pc/[r]_{f_1}  \ar@/_1.6pc/[rr]_{f_2}  \ar@/_2.3pc/[rrr]_{f_3} & v \ar@/^1pc/[rr]^{e_1}  \ar@/^1.6pc/[rrr]^{e_2}   \ar@/^2.3pc/[rrrr]^{e_3}
 & v   \ar@/_1pc/[rr]_{f_1}  \ar@/_1.6pc/[rrr]_{f_2}  \ar@/_2.3pc/[rrrr]_{f_3} & v  \ar@{.>}@/^2.3pc/[rrr]^{e_1} & v& v & \dots}.
\end{equation*}
\end{example}

Let $(F,\phi)$ be a representation graph for $(E,w)$. Let $\hat E_d$ and $F_d$ be the double graphs of $\hat E$ and $F$, respectively. Clearly the homomorphism $\phi :F\rightarrow \hat E$ induces a map $\Path(F_d)\rightarrow \Path(\hat E_d)$, which we also denote by $\phi$. The lemma below is easy to check.
 
\begin{lemma}[{\cite[Lemma 3]{hazrat-preusser-shchegolev}}]\label{lemwell} Let $(F,\phi)$ be a representation graph for $(E,w)$. Let $q, q'\in \Path(F_d)$ such that $\phi(q)=\phi(q')$. If $s(q)=s(q')$ or $r(q)=r(q')$, then $q=q'$.
\end{lemma} 

In Subsection 11.2 we will associate to any representation graph for $(E,w)$ a module for $L(E,w)$. The irreducible representation graphs defined below are precisely those representation graphs that yield a simple module.

\begin{definition}\label{defrepirres}
A representation graph $(F,\phi)$ for $(E,w)$ is called {\it irreducible} if 
$\phi({}_u\!\Path(F_d))\neq \phi({}_v\!\Path(F_d))$ for any $u\neq v\in F^0$.
\end{definition}

We leave it to the reader to check that the representation graph in Example \ref{exrepgr1} is not irreducible, while the representation graphs in Examples \ref{exrepgr2}, \ref{exrepgr3} and \ref{exrepgr4} are irreducible.

\subsubsection{The category of representation graphs}\label{catrepa} 

We denote by $\RG(E,w)$ the category whose objects are the representation graphs for $(E,w)$. A morphism $\alpha:(F,\phi)\to (G,\psi)$ in $\RG(E,w)$ is a graph homomorphism $\alpha:F\to G$ such that $\psi\alpha=\phi$. One checks easily that a morphism $\alpha:(F,\phi)\to (G,\psi)$ is an isomorphism if and only if $\alpha^0$ and $\alpha^1$ are bijective.

The lemma below and the two propositions thereafter are straightforward to check. 

\begin{lemma}[{\cite[Lemma 8]{hazrat-preusser-shchegolev}}]\label{lempath}
Let $(F,\phi)$ and $(G,\psi)$ be objects in $\RG(E,w)$. Let $u\in F^0$ and $v\in G^0$. If $\phi({}_u\!\Path(F_d))\subseteq \psi({}_v\!\Path(G_d))$, then $\phi({}_u\!\Path(F_d))=\psi({}_v\!\Path(G_d))$.
\end{lemma}
 
\begin{proposition}[{\cite[Proposition 9]{hazrat-preusser-shchegolev}}]\label{propcruc}
Let $\alpha:(F,\phi)\to (G,\psi)$ be a morphism in $\RG(E,w)$. If $u\in F^0$, then $\phi({}_u\!\Path(F_d))=\psi({}_{\alpha^0(u)}\!\Path(G_d))$.
\end{proposition}

Recall that a graph homomorphism $\alpha=(\alpha^0,\alpha^1):F\to G$ is called a {\it covering} if $\alpha^0$ and $\alpha^1$ are surjective, and for any $v\in F^0$ the maps $\alpha^1|_{r^{-1}(v)}: r^{-1}(v)\rightarrow r^{-1}(\alpha^0(v))$ and $\alpha^1|_{s^{-1}(v)}: s^{-1}(v)\rightarrow s^{-1}(\alpha^0(v))$ are bijective (i.e. $\alpha$ is surjective and ``locally'' bijective). 

\begin{proposition}[{\cite[Proposition 10]{hazrat-preusser-shchegolev}}]\label{propsurj}
Let $\alpha:(F,\phi)\to (G,\psi)$ be a morphism in $\RG(E,w)$. Then $\alpha$ is a covering. 
\end{proposition}

\subsubsection{Quotients of representation graphs}
For any object $(F,\phi)$ in $\RG(E,w)$ we define an equivalence relation $\sim$ on $F^0$ by $u\sim v$ if $\phi({}_u\!\Path(F_d))=\phi({}_v\!\Path(F_d))$. Recall that if $\sim$ and $\approx$ are equivalence relations on a set $X$, then one writes $\approx~\leq~ \sim $ and calls $\approx$ {\it finer} than $\sim$ (and $\sim$ {\it coarser} than $\approx$) if $x\approx y$ implies that $x\sim y$, for any $x,y\in X$. 

\begin{definition}\label{defadm}
Let $(F,\phi)$ be an object in $\RG(E,w)$. An equivalence relation $\approx$ on $F^0$ is called {\it admissible} if the following hold:
\begin{enumerate}[(i)]
\item $\approx~\leq~\sim$.
\item If $u\approx v$, $p\in{} _u\!\Path_x(F_d)$, $q\in{} _v\!\Path_y(F_d)$ and $\phi(p)=\phi(q)$, then $x\approx y$.
\end{enumerate}
\end{definition}

The admissible equivalence relations on $F^0$ (with partial order $\leq$) form a bounded lattice whose maximal element is $\sim$ and whose minimal element is the equality relation $=$.

Let $(F,\phi)$ be an object in $\RG(E,w)$ and $\approx$ an admissible equivalence relation on $F^0$. If $f, g\in F^1$ we write $f\approx g$ if $s(f)\approx s(g)$ and $\phi(f)=\phi(g)$. This defines an equivalence relation on $F^1$. Define an object $(F_\approx,\phi_\approx)$ in $\RG(E,w)$ by 
\begin{align*}
F_\approx^0&=F^0/\approx,\\
F_\approx^1&=F^1/\approx,\\
s([f])&=[s(f)],\\
r([f])&=[r(f)],\\
\phi_\approx^0([v])&=\phi^0(v),\\
\phi_\approx^1([f])&=\phi^1(f).
\end{align*}
We call $(F_\approx,\phi_\approx)$ a {\it quotient} of $(F,\phi)$. 

\begin{theorem}[{\cite[Theorem 12]{hazrat-preusser-shchegolev}}]\label{thmadm}
Let $(F,\phi)$ and $(G,\psi)$ be objects in $\RG(E,w)$. Then there is a morphism $\alpha:(F,\phi)\to (G,\psi)$ if and only if $(G,\psi)$ is isomorphic to a quotient of $(F,\phi)$.
\end{theorem}
\begin{proof}
$(\Rightarrow)$ Suppose there is a morphism $\alpha:(F,\phi)\to (G,\psi)$. If $u,v\in F^0$, we write $u\approx v$ if $\alpha^0(u)=\alpha^0(v)$. Clearly $\approx$ defines an equivalence relation on $F^0$. Below we check that $\approx$ is admissible. 
\begin{enumerate}[(i)]
\item Suppose $u\approx v$. Then $\phi({}_u\!\Path(F_d))=\psi({}_{\alpha^0(u)}\!\Path(G_d))=\psi({}_{\alpha^0(v)}\!\Path(G_d))=\phi({}_v\!\Path(F_d))$ by Proposition \ref{propcruc}. Hence $u\sim v$.
\item Suppose $u\approx v$, $p\in{} _u\!\Path_x(F_d)$, $q\in{} _v\!\Path_y(F_d)$ and $\phi(p)=\phi(q)$. Clearly $\alpha(p)\in{} _{\alpha^0(u)}\!\Path_{\alpha^0(x)}$ $(G_d)$ and $\alpha(q)\in{} _{\alpha^0(v)}\!\Path_{\alpha^0(y)}(G_d)$. Moreover, $\psi(\alpha(p))=\phi(p)=\phi(q)=\psi(\alpha(q))$. Since $\alpha^0(u)=\alpha^0(v)$, it follows from Lemma \ref{lemwell} that $\alpha(p)=\alpha(q)$. Hence $\alpha^0(x)=r(\alpha(p))=r(\alpha(q))=\alpha^0(y)$ and therefore $x\approx y$.
\end{enumerate}

Note that by Lemma \ref{lemwell} we have $f\approx g$ if and only if $\alpha^1(f)=\alpha^1(g)$, for any $f,g\in F^1$. Define a graph homomorphism $\beta:F_\approx\to G$ by $\beta^0([v])=\alpha^0(v)$ and $\beta^1([f])=\alpha^1(f)$. Clearly $\psi\beta=\phi_{\approx}$ and therefore $\beta:(F_\approx,\phi_\approx)\to (G,\psi)$ is a morphism. In view of Proposition \ref{propsurj}, $\beta^0$ and $\beta^1$ are bijective and hence $\beta$ is an isomorphism.

$(\Leftarrow)$ Suppose now that $(G,\psi)\cong (F_\approx,\phi_\approx)$ for some admissible equivalence relation $\approx$ on $F^0$. In order to show that there is a morphism $\alpha:(F,\phi)\to (G,\psi)$ it suffices to show that there is a morphism $\beta:(F,\phi)\to (F_\approx,\phi_\approx)$. But this is obvious (define $\beta^0(v)=[v]$ and $\beta^1(f)=[f]$).
\end{proof}

\subsubsection{The subcategories $\RG(E,w)_C$} \label{subsecnn}

Let $(F,\phi)$ and $(G,\psi)$ be objects in $\RG(E,w)$. We write $(F,\phi)\leftrightharpoons(G,\psi)$ if there is a $u\in F^0$ and a $v\in G^0$ such that $\phi({}_u\!\Path(F_d))=\psi({}_v\!\Path(G_d))$. One checks easily that $\leftrightharpoons$ defines an equivalence relation on $\Ob(\RG(E,w))$. If $C$ is a $\leftrightharpoons$-equivalence class, then we denote by $\RG(E,w)_C$ the full subcategory of $\RG(E,w)$ such that $\Ob(\RG(E,w)_C)=C$. If $\alpha:(F,\phi)\to (G,\psi)$ is a morphism in $\RG(E,w)$, then $(F,\phi)\leftrightharpoons(G,\psi)$ by Proposition \ref{propcruc}. Thus $\RG(E,w)$ is the disjoint union of the subcategories $\RG(E,w)_C$, where $C$ ranges over all $\leftrightharpoons$-equivalence classes.

Let $F$ be a graph and $F_d$ its double graph. A path $p=x_1\dots x_n\in \Path(F_d)$ is called \emph{backtracking} if there is a $1\leq j\leq n-1$ such that $x_jx_{j+1}=ff^*$ or $x_jx_{j+1}=f^*f$ for some $f\in  F^1$. We say that $p$ is \emph{reduced} if it is not backtracking.

Fix a $\leftrightharpoons$-equivalence class $C$, a representation graph $(F,\phi)\in C$ and a vertex $u\in F^0$. %We call a path $p=x_1\dots x_n\in \Path(\hat E_d)$ {\it backtracking}, if there is a $1\leq j\leq n-1$ such that $x_j x_{j+1}=y y^*$ or $x_j x_{j+1}=y^*y$ for some $y\in \hat E^1$. 
We denote by $\phi({}_u\!\Path(F_d))_{\nb}$ the set of all paths in $\phi({}_u\!\Path(F_d))$ that are reduced. Define a representation graph $(T,\xi)=(T_C,\xi_C)$ for $(E,w)$ by
\begin{align*}
T^0&=\{v_p\mid p\in \phi({}_u\!\Path(F_d))_{\nb}\},\\
T^1&=\{e_p\mid p\in \phi({}_u\!\Path(F_d))_{\nb},p\neq \phi(u)\},\\
s(e_{x_1\dots x_n})&=\begin{cases}v_{x_1\dots x_{n-1}},&\text{ if }x_n\in \hat E^1,\\v_{x_1\dots x_n},&\text{ if }x_n\in (\hat E^1)^*,\end{cases}\\
r(e_{x_1\dots x_n})&=\begin{cases}v_{x_1\dots x_n},&\text{ if }x_n\in \hat E^1,\\v_{x_1\dots x_{n-1}},&\text{ if }x_n\in (\hat E^1)^*,\end{cases}\\
\xi^0(v_{x_1\dots x_n})&=\begin{cases}\phi(u),&\text{ if }n=1 \text{ and }x_1=\phi(u),\\r_{\hat E}(x_n),&\text{ if }x_n\in \hat E^1,\\s_{\hat E}(x_n^*),&\text{ if }x_n\in (\hat E^1)^*,\end{cases}\\
\xi^1(e_{x_1\dots x_n})&=\begin{cases}x_n,&\text{ if }x_n\in \hat E^1,\\x_n^*,&\text{ if }x_n\in (\hat E^1)^*.\end{cases}
\end{align*}
Here we use the convention that if $x_1\dots x_n \in \phi({}_u\!\Path(F_d))_{\nb}$, where $n=1$, then $x_1\dots x_{n-1}=\phi(u)$. % Clearly $T$ is nonempty and connected. 

\begin{example}\label{exuniv}
Suppose $(E,w)$ is the weighted graph $\xymatrix{
 v \ar@(dl,ul)^{e,2} \ar@(dr,ur)_{f,2} &
}$
and $(F,\phi)$ is the representation graph for $(E,w)$ given by
\[F:\xymatrix{
 u \ar@(dl,ul)^{g} \ar@(dr,ur)_{h} &
}\]
and $\phi^0(u)=v$, $\phi^1(g)=e_1$, $\phi^1(h)=f_2$. Then $\phi({}_u\!\Path(F_d))_{\nb}$ consists of all reduced paths in $\Path(\hat E_d)$ whose letters come from the alphabet $\{v,e_1,e_1^*,f_2,f_2^*\}$. Let $C$ be the $\leftrightharpoons$-equivalence class of $(F,\phi)$. Then $T_C$ is the graph
\[\xymatrix@C=7pt@R=20pt{
&&&&&&&v_v\ar[ddrr]^(.6){e_{e_1}}\ar[ddrrrrrr]^{e_{f_2}}&&&&&&&\\
&&&&&&&&&&&&&&\\
&v_{e_1^*}\ar[uurrrrrr]^{e_{e_1^*}}\ar[ddr]^(.5){e_{e_1^*f_2}}&&&&v_{f_2^*}\ar[uurr]^(.4){e_{f_2^*}}\ar[ddr]^(.5){e_{f_2^*e_1}}&&&&v_{e_1}\ar[dd]^(.7){e_{e_1e_1}}\ar[ddr]^(.5){e_{e_1f_2}}&&&&v_{f_2}\ar[dd]^(.7){e_{f_2e_1}}\ar[ddr]^(.5){e_{f_2f_2}}&\\
&&&&&&&&&&&&&&\\
v_{e_1^*e_1^*}\ar[uur]^(.4){e_{e_1^*e_1^*}}&v_{e_1^*f_2^*}\ar[uu]^(.3){e_{e_1^*f_2^*}}&v_{e_1^*f_2}&&v_{f_2^*e_1^*}\ar[uur]^(.4){e_{f_2^*e_1^*}}&v_{f_2^*f_2^*}\ar[uu]^(.3){e_{f_2^*f_2^*}}&v_{f_2^*e_1}&&v_{e_1f_2^*}\ar[uur]^(.4){e_{e_1f_2^*}}&v_{e_1e_1}&v_{e_1f_2}&&v_{f_2e_1^*}\ar[uur]^(.4){e_{f_2e_1^*}}&v_{f_2e_1}&v_{f_2f_2}.\\
&\vdots&&&&\vdots&&&&\vdots&&&&\vdots&
}\]
\end{example}

\begin{proposition}[{\cite[Proposition 14]{hazrat-preusser-shchegolev}}]\label{propuniv}
If $(G,\psi)\in C$, then there is a morphism $\alpha:(T_C,\xi_C)\to (G,\psi)$. 
\end{proposition}
\begin{proof}
Since $(G,\psi)\leftrightharpoons(F,\phi)$, there is a $v\in G^0$ such that $\phi({}_u\!\Path(F_d))=\psi({}_v\!\Path(G_d))$.
Define a homomorphism $\alpha:T_C\to G$ as follows. Let $x_1\dots x_n\in \phi({}_u\!\Path(F_d))_{\nb}$. Since $\phi({}_u\!\Path(F_d))=\psi({}_v\!\Path(G_d))$, there is a (uniquely determined) path $y_1\dots y_n\in {}_v\!\Path(G_d)$ such that $\psi(y_1\dots y_n)=x_1\dots x_n$. Define $\alpha^0(v_{x_1\dots x_n})=r(y_n)$, $\alpha^1(e_{x_1\dots x_n})=y_n$ if $y_n\in G^1$ and respectively $\alpha^1(e_{x_1\dots x_n})=y_n^*$ if $y_n\in (G^1)^*$. We leave it to the reader to check that $\alpha$ is a graph homomorphism and that $\psi\alpha=\xi_C$.
\end{proof}

\begin{corollary}[{\cite[Corollary 15]{hazrat-preusser-shchegolev}}]\label{corcat1}
Up to isomorphism the representation graphs in $C$ are precisely the quotients of $(T_C,\xi_C)$, and consequently  \[(S_C,\zeta_C):=((T_C)_\sim,(\xi_C)_\sim)\]  is the unique irreducible representation graph in $C$.
\end{corollary}
\begin{proof}
The first statement follows from Theorem \ref{thmadm} and Proposition \ref{propuniv}. The second statement now follows since a quotient $((T_C)_\approx,(\xi_C)_\approx)$ satisfies the condition in Definition \ref{defrepirres} if and only if $\approx$ equals $\sim$. 
\end{proof}

Recall that an object $X$ in a category {\bf $C$} is called {\it repelling} (resp. {\it attracting}) if for any object $Y$ in {\bf $C$} there is a morphism $X\rightarrow Y$ (resp. $Y\rightarrow X$). By Proposition \ref{propuniv}, $(T_C,\xi_C)$ is a repelling object in $C$. On the other hand, if $(G,\psi)$ is an object in $C$, then clearly $(S_C,\zeta_C)$ is isomorphic to a quotient of $(G,\psi)$. It follows from Theorem \ref{thmadm} that $(S_C,\zeta_C)$ is an attracting object in $C$.

\begin{example}\label{excatone}
Suppose $(E,w)$ is the weighted graph
$\xymatrix{
 v \ar@(dl,ul
 )^{e,2} \ar@(dr,ur)_{f,2} &
}$. Consider the representation graphs $(F_1,\phi_1),\dots, (F_7,\phi_7)$ for $(E,w)$ given below.\\
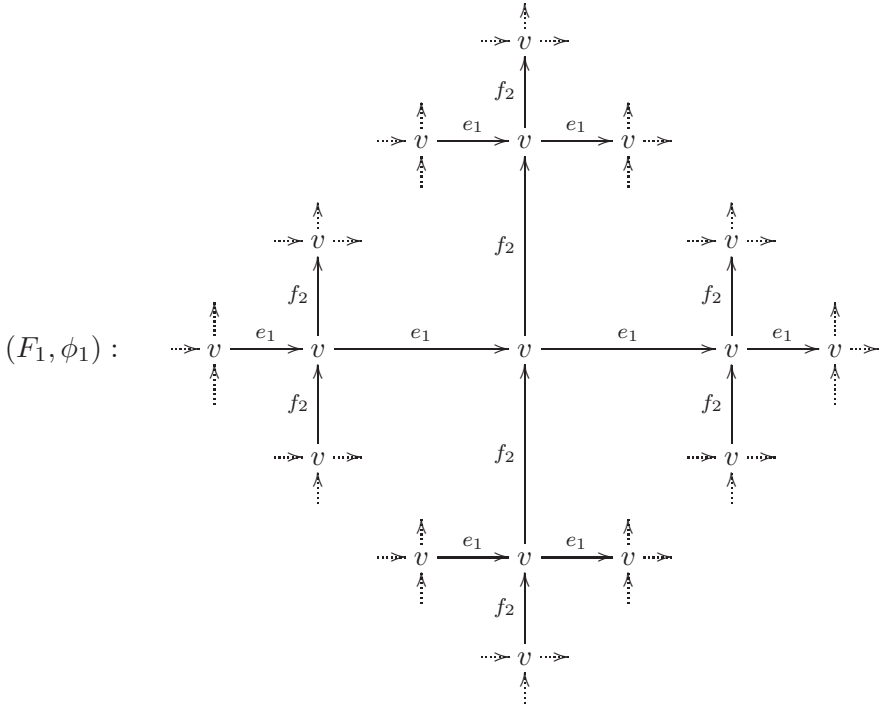
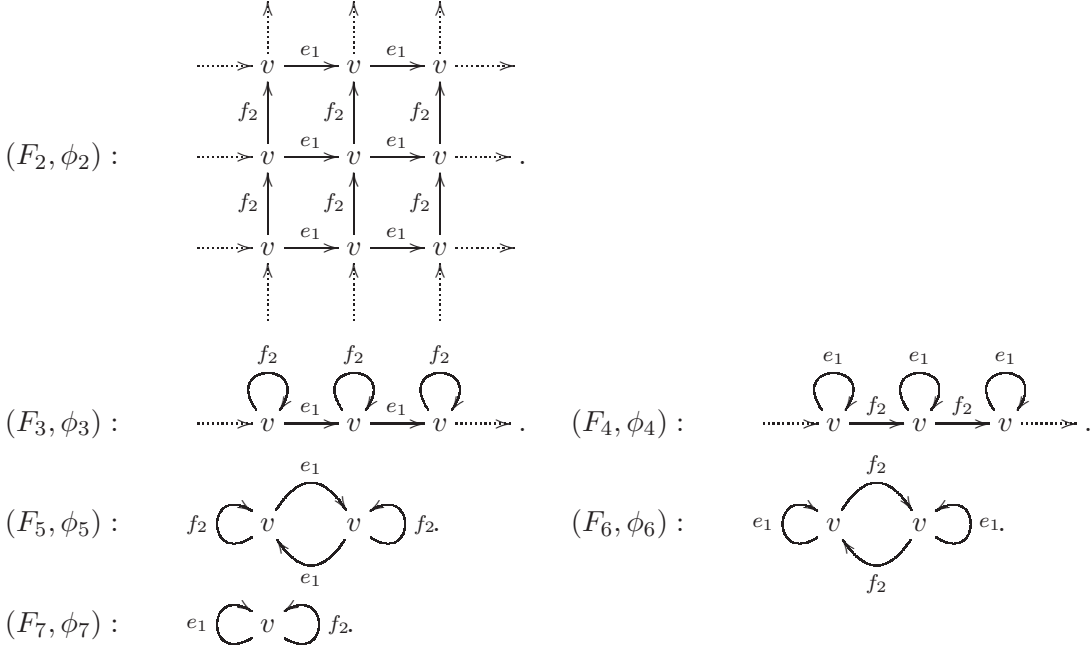
\begin{tabular}{l}
$\xymatrix@C=10pt@R=10pt{
&&&&&&&&&&&&&&&\\
&&&&&&&\ar@{..>}[r]&v\ar@{..>}[u]\ar@{..>}[r]&&&&&&&\\
&&&&&&&&&&&&&&&\\
&&&&&\ar@{..>}[r]&v\ar@{..>}[u]\ar^{e_1}[rr]&&v\ar[uu]^{f_2}\ar^{e_1}[rr]&&v\ar@{..>}[u]\ar@{..>}[r]&&&&&\\
&&&&&&\ar@{..>}[u]&&&&\ar@{..>}[u]&&&&&\\
&&&\ar@{..>}[r]&v\ar@{..>}[r]\ar@{..>}[u]&&&&&&&\ar@{..>}[r]&v\ar@{..>}[u]\ar@{..>}[r]&&&\\
&&&&&&&&&&&&&&&\\
(F_1,\phi_1):&\ar@{..>}[r]&v\ar@{..>}[u]\ar^{e_1}[rr]&&v\ar^{f_2}[uu]\ar^{e_1}[rrrr]&&&&v\ar^{f_2}[uuuu]\ar^{e_1}[rrrr]&&&&v\ar^{f_2}[uu]\ar^{e_1}[rr]&&v\ar@{..>}[u]\ar@{..>}[r]&.\\
&&\ar@{..>}[u]&&&&&&&&&&&&\ar@{..>}[u]&\\
&&&\ar@{..>}[r]&v\ar^{f_2}[uu]\ar@{..>}[r]&&&&&&&\ar@{..>}[r]&v\ar^{f_2}[uu]\ar@{..>}[r]&&&\\
&&&&\ar@{..>}[u]&&&&&&&&\ar@{..>}[u]&&&\\
&&&&&\ar@{..>}[r]&v\ar@{..>}[u]\ar^{e_1}[rr]&&v\ar^{f_2}[uuuu]\ar^{e_1}[rr]&&v\ar@{..>}[u]\ar@{..>}[r]&&&&&\\
&&&&&&\ar@{..>}[u]&&&&\ar@{..>}[u]&&&&&\\
&&&&&&&\ar@{..>}[r]&v\ar@{..>}[r]\ar^{f_2}[uu]&&&&&&&\\
&&&&&&&&\ar@{..>}[u]&&&&&&&
}$
\end{tabular}\\
\begin{tabular}{ll}
$\xymatrix@C=20pt@R=20pt{
&&&&&\\
&\ar@{..>}[r]&v\ar^{e_1}[r]\ar@{..>}[u]&v\ar^{e_1}[r]\ar@{..>}[u]&v\ar@{..>}[r]\ar@{..>}[u]&\\
(F_2,\phi_2):&\ar@{..>}[r]&v\ar^{e_1}[r]\ar^{f_2}[u]&v\ar^{e_1}[r]\ar^{f_2}[u]&v\ar^{f_2}[u]\ar@{..>}[r]&.\\
&\ar@{..>}[r]&v\ar^{e_1}[r]\ar^{f_2}[u]&v\ar^{e_1}[r]\ar^{f_2}[u]&v\ar^{f_2}[u]\ar@{..>}[r]&\\
&&\ar@{..>}[u]&\ar@{..>}[u]&\ar@{..>}[u]&
}$\\
$\xymatrix@C=20pt@R=20pt{
(F_3,\phi_3):&\ar@{..>}[r]&v\ar@(ul,ur)^{f_2}\ar^{e_1}[r]&v\ar@(ul,ur)^{f_2}\ar^{e_1}[r]&v\ar@(ul,ur)^{f_2}\ar@{..>}[r]&.}$
&$\xymatrix@C=20pt@R=20pt{
(F_4,\phi_4):&\ar@{..>}[r]&v\ar@(ul,ur)^{e_1}\ar^{f_2}[r]&v\ar@(ul,ur)^{e_1}\ar^{f_2}[r]&v\ar@(ul,ur)^{e_1}\ar@{..>}[r]&.
}$\\
$\xymatrix@C=20pt@R=20pt{
(F_5,\phi_5):&&v\ar@(dl,ul)^{f_2}\ar@/^1.3pc/^{e_1}[r]&v\ar@(dr,ur)_{f_2}\ar@/^1.3pc/^{e_1}[l]&.}$
&$\xymatrix@C=20pt@R=20pt{
(F_6,\phi_6):&&v\ar@(dl,ul)^{e_1}\ar@/^1.3pc/^{f_2}[r]&v\ar@(dr,ur)_{e_1}\ar@/^1.3pc/^{f_2}[l]&.}$\\
$\xymatrix@C=20pt@R=20pt{
(F_7,\phi_7):&&v \ar@(dl,ul)^{e_1} \ar@(dr,ur)_{f_2}&.}$&
\end{tabular}\\
\\
\\
All the representation graphs $(F_i,\phi_i)~(1\leq i\leq 7)$ lie in the same $\leftrightharpoons$-equivalence class $C$. One checks easily that $(F_1,\phi_1)\cong (T_C,\xi_C)$ and $(F_7,\xi_7)\cong (S_C,\zeta_C)$ (cf. Example \ref{exuniv}). Moreover, we have
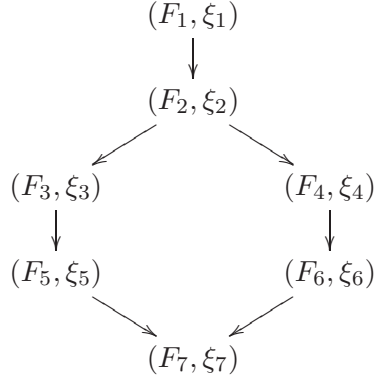
\[\xymatrix@C=10pt@R=15pt{
&(F_1,\xi_1)\ar[d]&\\
&(F_2,\xi_2)\ar[dl]\ar[dr]&\\
(F_3,\xi_3)\ar[d]&&(F_4,\xi_4)\ar[d]\\
(F_5,\xi_5)\ar[dr]&&(F_6,\xi_6)\ar[dl]\\
&(F_7,\xi_7)&
}\]
where an arrow $(F_i,\phi_i)\longrightarrow(F_j,\phi_j)$ means that $(F_j,\phi_j)$ is a quotient of $(F_i,\phi_i)$.
\end{example}

\subsection{Representations of weighted Leavitt path algebras via representation graphs}\label{weightedrbgh}

\subsubsection{The functor V}

For an object $(F,\phi)$ in $\RG(E,w)$, let $V_{(F,\phi)}$ be the $K$-vector space with basis $F^0$. For any $u\in E^0$, $e\in E^1$ and $1\leq i\leq w(e)$,  define endomorphisms $\sigma_u,\sigma_{e_i},\sigma_{e_i^*}\in \End_K(V_{(F,\phi)})$ by
\begin{align*}
\sigma_u(v)          &=\begin{cases}v,\quad\quad\quad             &\text{if }\phi^0(v)=u \\0,                                           & \text{otherwise} \end{cases},\\
\sigma_{e_i}(v)     &=\begin{cases}r_F(f),\quad                       &\text{if }\exists f\in s_F^{-1}(v):~\phi^1(f)=e_i\\0,     & \text{otherwise} \end{cases},\\
\sigma_{e_i^*}(v)  &=\begin{cases}s_F(f),\quad                        &\text{if }\exists f\in r_F^{-1}(v):~\phi^1(f)=e_i\\0,    & \text{otherwise} \end{cases},
\end{align*}
where $v\in F^0$. Then there is an $K$-algebra homomorphism $\pi:L(E,w)\rightarrow \End_K(V_{(F,\phi)})^{\op}$ such that $\pi(u)=\sigma_u$, $\pi(e_i)=\sigma_{e_i}$ and $\pi(e_i^*)=\sigma_{e_i^*}$. We call this representation the {\it representation of $L(E,w)$ defined by $(F,\phi)$}. Clearly $V_{(F,\phi)}$ becomes a right $L(E,w)$-module by defining $x \cdot a:= \pi(a)(x)$ for any $a\in L(E,w)$ and $x\in V_{(F,\phi)}$. A morphism $\alpha:(F,\phi)\rightarrow(G,\psi)$ in $\RG(E,w)$ induces a surjective $L(E,w)$-module homomorphism $V_\alpha:V_{(F,\phi)}\rightarrow V_{(G,\psi)}$ such that $V_\alpha(u)=\alpha^0(u)$ for any $u\in F^0$. We obtain a functor 
\[V:\RG(E,w)\to \Mod\text{-}L(E,w)\]
where $\Mod$-$L(E,w)$ denotes the category of unital right $L(E,w)$-modules.

The lemma below is easy to check. It describes the action of monomial elements of the weighted Leavitt path algebra $L(E,w)$ on the $K$-vector space $V_{(F,\phi)}$. Note that by Lemma \ref{lemwell}, for any $p\in \Path(\hat E_d)$ and $u\in F^0$ there is at most one $v\in F^0$ such that $p\in\phi({}_u\!\Path_v(F_d))$.
 
\begin{lemma}[{\cite[Lemma 21]{hazrat-preusser-shchegolev}}]\label{lemaction}
Let $(F,\phi)$ be an object in $\RG(E,w)$. If $p\in \Path(\hat E_d)$ and $u\in F^0$, then
\[u\cdot p=\begin{cases}
v,\quad&\text{if }p\in \phi({}_u\!\Path_v(F_d)),\text{ for some }v\in F^0,\\
0,\quad&\text{otherwise.}
\end{cases}\]
\end{lemma}

\begin{corollary}[{\cite[Corollary 22]{hazrat-preusser-shchegolev}}]\label{coraction}
Let $(F,\phi)$ be an object of $\RG(E,w)$. If $a=\sum_{p\in \Path(\hat E_d)}k_pp\in L(E,w)$ and $u\in F^0$, then
\[u\cdot a=\sum_{v\in F^0}\Big(\sum_{p\in\phi({}_u\!\Path_v(F_d))}k_p\Big)v.\]
\end{corollary}

The example below shows that the functor $V$ is not full, namely, there can be $L(E,w)$-module homomorphisms $V_{(F,\phi)}\rightarrow V_{(G,\psi)}$ that are not induced by a morphism $(F,\phi)\rightarrow(G,\psi)$. 

\begin{example}\label{excat11}
Suppose $(E,w)$ is the weighted graph $\xymatrix{
 v \ar@(dr,ur)_{e,1} &}$. Let $(F,\phi)$ and $(G,\psi)$ be the representation graphs for $(E,w)$ given by \[F:\xymatrix{
 u \ar@(dr,ur)_{f} &},\quad\quad G:\xymatrix{u_1 \ar@/^1.3pc/^{f^{(1)}}[r] &u_2\ar@/^1.3pc/^{f^{(2)}}[l]}\]
and $\phi^0(u)=v$, $\phi^1(f)=e_1$, $\psi^0(u_1)=\psi^0(u_2)=v$, $\psi^1(f^{(1)})=\psi^1(f^{(2)})=e_1$. 
%Clearly $(F,\phi)$ is isomorphic to a quotient of $(G,\psi)$. Hence, by Theorem \ref{thmadm}, there is a morphism $\alpha:(G,\psi)\rightarrow (F,\phi)$ and therefore a morphism $V_\alpha:V_{(G,\psi)}\rightarrow V_{(F,\phi)}$.
Since $G$ has more vertices than $F$, $(G,\psi)$ cannot be isomorphic to a quotient of $(F,\phi)$. Hence, by Theorem \ref{thmadm}, there is no morphism $(F,\phi)\rightarrow (G,\psi)$. But there is a $L(E,w)$-module homomorphism $\sigma:V_{(F,\phi)}\rightarrow V_{(G,\psi)}$ such that $\sigma(u)=u_1+u_2$. %Note that $V_\alpha$ and $\sigma$ are not inverse to each other.
\end{example}

Let $(F,\phi)$ and $(G,\psi)$ be objects in $\RG(E,w)$. If $(F,\phi)\cong(G,\psi)$, then clearly $V_{(F,\phi)}\cong V_{(G,\psi)}$. It is not known if the converse implication is also true. But we will see that $V_{(F,\phi)}\cong V_{(G,\psi)}$ implies at least that $(F,\phi)\leftrightharpoons(G,\psi)$ (i.e. that $(F,\phi)$ and $(G,\psi)$ lie in the same subcategory $\RG(E,w)_C$ of $\RG(E,w)$).

\begin{lemma}[{\cite[Lemma 31]{hazrat-preusser-shchegolev}}]\label{lemmodule}
Let $(F,\phi)$ and $(G,\psi)$ be objects in $\RG(E,w)$ and let $\sigma:V_{(F,\phi)}\rightarrow V_{(G,\psi)}$ be an $L(E,w)$-module homomorphism. Let $u\in F^0$ and $\sigma(u)=\sum_{s=1}^nk_sv_s$, where $n\geq 1$, $k_1,\dots,k_n\in K^\times$ and $v_1,\dots,v_n$ are pairwise distinct vertices from $G^0$. Then $\phi({}_u\!\Path(F_d))=\psi({}_{v_s}\!\Path(G_d))$ for any $1\leq s\leq n$.
\end{lemma}
\begin{proof} 
Let $p\in \Path(\hat E_d)$ such that $p\not\in\phi({}_u\!\Path(F_d))$. Then
\[0=\sigma(0)=\sigma(u\cdot p)=\sigma(u)\cdot p=\sum_{s=1}^nk_sv_s\cdot p=\sum_{s=1}^nk_s(v_s\cdot p)\]
by Lemma \ref{lemaction}. One more application of Lemma \ref{lemaction} gives that $v_s\cdot p =0$, for any $1\leq s\leq n$, whence $p\not\in \psi({}_{v_s}\!\Path(G_d))$ for any $1\leq s\leq n$. Hence we have shown that $\phi({}_u\!\Path(F_d))\supseteq\psi({}_{v_s}\!\Path(G_d))$ for any $1\leq s\leq n$. It follows from Lemma \ref{lempath} that $\phi({}_u\!\Path(F_d))=\psi({}_{v_s}\!\Path(G_d))$ for any $1\leq s\leq n$.
\end{proof}

\begin{proposition}[{\cite[Proposition 32]{hazrat-preusser-shchegolev}}]\label{prop42}
Let $(F,\phi)$ and $(G,\psi)$ be objects in $\RG(E,w)$. If there is a nonzero $L(E,w)$-module homomorphism $\sigma:V_{(F,\phi)}\rightarrow V_{(G,\psi)}$, then $(F,\phi)\leftrightharpoons (G,\psi)$.
\end{proposition}
\begin{proof}
The  proposition follows from Lemma \ref{lemmodule} and the definition of $\leftrightharpoons$.
\end{proof}

\begin{proposition}[{\cite[Proposition 33]{hazrat-preusser-shchegolev}}]\label{prop43}
Let $(F,\phi)$ and $(G,\psi)$ be irreducible representation graphs for $(E,w)$. Then $V_{(F,\phi)}\cong V_{(G,\psi)}$ as $L(E,w)$-modules if and only if $(F,\phi)\cong(G,\psi)$.
\end{proposition}
\begin{proof}
Clearly isomorphic objects in $\RG(E,w)$ yield isomorphic $L(E,w)$-modules. Hence we only have to show that $V_{(F,\phi)}\cong V_{(G,\psi)}$ implies $(F,\phi)\cong(G,\psi)$. Suppose that $V_{(F,\phi)}\cong V_{(G,\psi)}$. Then $(F,\phi)\leftrightharpoons(G,\psi)$ by Proposition \ref{prop42}, i.e. $(F,\phi)$ and $(G,\psi)$ are in the same $\leftrightharpoons$-equivalence class $C$. Since they are irreducible, it follows from Corollary \ref{corcat1} that $(F,\phi)\cong (S_C,\zeta_C)\cong (G,\psi)$.
\end{proof}

\subsubsection{Gradedness of the modules $V_{(F,\phi)}$}
Recall that the standard grading on $L(E,w)$ is a $\Z^{\lambda}$-grading where $\lambda=\sup\{w(e) \mid e \in E^{1}\}$ if this supremum is finite and otherwise $\lambda=\omega$ where $\omega$ is the smallest infinite ordinal. %Let $(F,\phi)$ be a representation graph for $(E,w)$. We define the map $d:\Path(F)\to\Z^{\lambda}$ by $d(p)=\deg(\phi(p))$ where $\deg(\phi(p))$ denotes the homogeneous degree of $\phi(p)$.

\begin{lemma}[{\cite[Lemma 24]{hazrat-preusser-shchegolev}}]\label{windynhgf}
Let $R$ be a $\Gamma$-graded ring and $M$ a right $R$-module, where $\Gamma$ is a totally ordered abelian group. If there is a homogeneous element $r\in R$ with $\deg(r)\not = 0$ and $0\not = m\in M$ such that $m\cdot r=m$, then $M$ cannot be $\Gamma$-graded.  
\end{lemma}

\begin{theorem}[{\cite[Theorem 25]{hazrat-preusser-shchegolev}}]\label{gradedrepres}
Let $(F,\phi)$ be an object in $\RG(E,w)$. Then the $L(E,w)$-module $V_{(F,\phi)}$ is graded with respect to the standard grading of $L(E,w)$ if and only if $\deg(\phi(p))=\deg(\phi(q))$ for any $u,v\in F^0$ and $p,q\in{}_u\!\Path_v(F_d)$.  
\end{theorem}
\begin{proof}
Set $L:=L(E,w)$ and $V:=V_{(F,\phi)}$. First suppose that $\deg(\phi(p))=\deg(\phi(q))$ for any $u,v\in F^0$ and $p,q\in{}_u\!\Path_v(F_d)$. Choose a $u\in F^0$. Define a map $\deg:F^0\rightarrow \mathbb Z^\lambda$ by $\deg(v)=\deg(\phi(p))$ where $p\in {}_u\!\Path_v(F_d)$. For $\alpha \in \mathbb Z^\lambda$ set $V_\alpha:=\bigoplus_{v\in F^0,\deg(v)=\alpha}Kv$. Then clearly $V=\bigoplus_{\alpha \in \mathbb Z^\lambda} V_\alpha$. Next we check that $V_\alpha L_\beta \subseteq V_{\alpha+\beta}$ for any $\alpha,\beta \in \mathbb Z^\lambda$. It suffices to show that $v\cdot p \in V_{\alpha+\beta}$ for any $v\in V_{\alpha} \cap F^0$ and $p\in\Path(\hat E_d)$ such that $\deg(p)=\beta$. Assume that $v\cdot p\neq 0$. Then, by Lemma~\ref{lemaction} there is an $x\in F^0$ and a $q\in{}_v\!\Path_x(F_d)$ such that $\phi(q)=p$. Since $\deg(v)=\alpha$, there is a $t\in {}_u\!\Path_v(F_d)$ such that $\deg(\phi(t))=\alpha$. Since $tq\in {}_u\!\Path_x(F_d)$, we have $\deg(x)=\deg(\phi(tq))=\deg(\phi(t))+\deg(\phi(q))=\alpha+\beta$. Thus $x\in V_{\alpha+\beta}$.  

Suppose now that $V$ is a graded $L$-module. Assume there are $u,v\in F^0$ and $p,q\in{}_u\!\Path_v(F_d)$ such that $\deg(\phi(p))\neq\deg(\phi(q))$. Set $t:=\phi(pq^*)$. Then clearly $\deg (t)\neq 0$ and $u\cdot t=u$. It follows from Lemma~\ref{windynhgf} that $V$ is not graded, which is a contradiction.  
\end{proof}

\begin{corollary}[{\cite[Corollary 27]{hazrat-preusser-shchegolev}}]\label{acyclicrep}
Let $C$ be a $\leftrightharpoons$-equivalence class. Then the module $V_{(T_C,\xi_C)}$ is graded. 
\end{corollary}

\begin{example}
Consider the weighted graph $ \xymatrix{
 v \ar@(dl,ul)^{e,2} \ar@(dr,ur)_{f,2} &
}$ from Example~\ref{excatone}. By Theorem~\ref{gradedrepres}, the $L(E,w)$-modules $V_{(F_1,\phi_1)}$ and $V_{(F_2,\phi_2)}$ are graded whereas $V_{(F_3,\phi_3)}, V_{(F_4,\phi_4)}, V_{(F_5,\phi_5)}, V_{(F_6,\phi_6)}$ and the simple module $V_{(F_7,\phi_7)}$ are not graded.   
\end{example}

\subsubsection{Simplicity of the modules $V_{(F,\phi)}$}

The lemma below is easy to check.
\begin{lemma}[{\cite[Lemma 63]{hazrat-preusser-shchegolev}}]\label{lembasis}
Let $W$ be a $K$-vector space and $B$ a linearly independent subset of  $W$. Let $k_i\in K$ and $u_i,v_i\in B$, where $1\leq i \leq n$. Then $\sum_{s=1}^nk_s(u_s-v_s)\not\in B$.
\end{lemma}

\begin{theorem}[{\cite[Theorem 28]{hazrat-preusser-shchegolev}}]\label{thmirr} 
Let $(F,\phi)$ be an object in $\RG(E,w)$. Then the following are equivalent.
\begin{enumerate}[(i)]
\item $V_{(F,\phi)}$ is simple.
\medskip
\item For any $x\in V_{(F,\phi)}\setminus\{0\}$ there is an $a\in L(E,w)$ such that $x \cdot a \in F^0$.
\medskip
\item For any $x\in V_{(F,\phi)}\setminus\{0\}$ there is a $k\in K$ and a $p\in \Path(\hat E_d) \text{ such that } x \cdot kp \in F^0$.
\medskip
\item $(F,\phi)$ is irreducible.
\end{enumerate}
\end{theorem}

\begin{proof} 
(i) $\Rightarrow$ (iv). Assume that there are $u\neq v\in F^0$ such that $\phi({}_u\!\Path(F_d))= \phi({}_v\!\Path(F_d))$. Consider the submodule $(u-v)\cdot L(E,w)\subseteq V_{(F,\phi)}$. Since $V_{(F,\phi)}$ is simple by assumption, we have $(u-v)\cdot L(E,w)=V_{(F,\phi)}$. Hence there is an $a\in L(E,w)$ such that $(u-v)\cdot a=v$. Clearly there is an $n\geq 1$, $k_1,\dots,k_n\in K^\times$ and pairwise distinct $p_1,\dots,p_n\in \Path(\hat E_d)$ such that $a=\sum_{s=1}^nk_sp_s$. We may assume that $(u-v)\cdot p_s\neq 0$ for any $1\leq s\leq n$. It follows from Lemma \ref{lemaction} that $p_s\in \phi({}_u\!\Path(F_d))= \phi({}_v\!\Path(F_d))$ for any $s$ and moreover, that $(u-v)\cdot p_s=u_s-v_s$ for some distinct $u_s,v_s\in F^0$. Hence 
\[v=(u-v)\cdot a=(u-v)\cdot(\sum_{s=1}^nk_sp_s)=\sum_{s=1}^nk_s(u_s-v_s)\]
which contradicts Lemma \ref{lembasis}.

\medskip 
(iv) $\Rightarrow$ (iii). Let $x\in V_{(F,\phi)}\setminus \{0\}$. Then there is an $n\geq 1$, pairwise disjoint $v_1,\dots, v_n\in F^0$ and $k_1,\dots,k_n\in K^\times$ such that $x=\sum_{s=1}^nk_sv_s$. If $n=1$, then $x\cdot k_1^{-1}\phi^{0}(v_1)= v_1$. Suppose now that $n>1$. By assumption, we can choose a $p_1\in \phi({}_{v_1}\!\Path(F_d))$ such that $p_1\not\in\phi({}_{v_2}\!\Path(F_d))$. Clearly $x\cdot p_1\neq 0$  is a linear combination of at most $n-1$ vertices from $F^0$. Proceeding this way, we obtain paths $p_1,\dots,p_m$ such that $x\cdot p_1\dots p_m=kv$ for some $k\in K^\times$ and $v\in F^0$. Hence $x\cdot k^{-1}p_1\dots p_m=v$.

\medskip 

(iii) $\Rightarrow$ (ii). Trivial.

\medskip 

(ii) $\Rightarrow$ (i). Let $U\subseteq V_{(F,\phi)}$ be a nonzero $L(E,w)$-submodule and $x\in U\setminus\{0\}$. By assumption, there is an $a\in L(E,w)$ and a $v\in F^0$ such that $v=x\cdot a\in U$. Let now $v'$ be an arbitrary vertex in $F^0$. Since $F$ is connected, there is a $p\in {}_{v}\!\Path_{v'}(F_d)$. It follows from Lemma \ref{lemaction} that $v'=v\cdot \phi(p)\in U$. Hence $U$ contains $F^0$ and thus $U=V_{(F,\phi)}$.
%(iv) $\Longleftrightarrow$ (v). This follows from the fact that $\phi(e^*)=\phi(e)^*$, $e\in F_d^1$ and ${}_u\Path(F_d)^*={}\Path_u(F_d)$ for any $u\in F^0$. 
\end{proof}

\subsubsection{Indecomposability of the modules $V_{(F,\phi)}$}
Recall that for a ring $R$, an $R$-module is called \emph{indecomposable} if it is non-zero and cannot be written as a direct sum of two non-zero submodules. It is easy to see that an $R$-module $M$ is indecomposable if and only if the only idempotent elements of the endomorphism ring $\End_R(M)$ are $0$ and $1$. 

Let $C$ be a $\leftrightharpoons$-equivalence class and define the representation graphs $(S_C,\zeta_C)$ and $(T_C,\xi_C)$ as in \S 11.1.4. Then $V_{(S_C,\zeta_C)}$ is indecomposable since it is simple. One can show that $V_{(T_C,\xi_C)}$ is also  indecomposable, see the theorem below. In general the indecomposibility of $V_{(F,\phi)}$, for a representation graph $(F,\phi)$, depends on the ground field $K$, see \cite[Example 34]{hazrat-preusser-shchegolev}.

\begin{theorem}[{\cite[Theorem 43]{hazrat-preusser-shchegolev}}]\label{thmdecomp}
Let $C$ be a $\leftrightharpoons$-equivalence class. Then the $L(E,w)$-module $V_{(T_C,\xi_C)}$ is indecomposable.
\end{theorem}
\begin{proofs}
Recall from \S 11.1.4 that $T_C^0=\{v_p\mid p\in \phi({}_u\!\Path(F_d))_{\nb}\}$ where $(F,\phi)$ is some fixed representation graph in $C$ and $u$ is a vertex in $F$. Define the set
\[
G:=\big \{p\in \phi({}_u\!\Path(F_d))_{\nb}\mid v_p\sim v_{\phi(u)}\big \}
\]
($\sim$ is defined in \S 11.1.3). For $p,p'\in G$ we define a reduced path $p\ast p'\in \Path(\hat E_d)$ as follows. If $p,p'\in G\setminus \{\phi(u)\}$ let $p\ast p'$ be the element of $\Path(\hat E_d)$ one gets by removing all subwords of the form $e_ie_i^*$ and $e_i^*e_i$ from the juxtaposition $pp'$ (if $p'=p^*$, then $p\ast p':=\phi(u))$. Moreover, define $\phi(u)\ast p:=p$, $p\ast \phi(u):=p$ and $\phi(u)\ast\phi(u):=\phi(u)$. One can show that $(G,\ast)$ is a free group, see \cite[Proposition 39]{hazrat-preusser-shchegolev}.

Let $A$ be the subalgebra of $L(E,w)$ generated by the image of the group $G$ in $L(E,w)$, and set $W:=\bigoplus_{p\in G}Kv_p\subseteq V_{(T_C,\xi_C)}$. Then $W$ is a right $A$-module where the action of $A$ on $W$ is induced by the action of $L(E,w)$ on $W$. Set $\bar A:=A/\ann(W)$ where $\ann(W)$ denotes the annihilator of the $A$-module $W$. Then $W$ is also a right $\bar A$-module where the action of $\bar A$ on $W$ is induced by the action of $A$ on $W$. One can show that the $K$-algebra $\bar A$ is isomorphic to the group algebra $K[G]$, see \cite[Proposition 41]{hazrat-preusser-shchegolev}. Moreover, the $A$-module $W$ is free of rank $1$ as an $\bar A$-module, see \cite[Proposition 42]{hazrat-preusser-shchegolev}. 

Let $\epsilon$ be an idempotent endomorphism of the $L(E,w)$-module $V_{(T_C,\xi_C)}$. It follows from Lemma \ref{lemmodule} that $\epsilon(W)=W$ and hence $\epsilon|_W\in \End_K(W)$. Clearly we have also $\epsilon|_W\in \End_{\bar A}(W)$. By the previous paragraph we have $\End_{\bar A}(W)\cong \bar A\cong K[G]$. Since $G$ is free, the group ring $K[G]$ has no zero divisors by \cite[Theorem 12]{higman}. It follows that $0$ and $1$ are the only idempotents of $K[G]$ whence $\epsilon|_W=0$ or $\epsilon|_W=\id_W$. Clearly $\epsilon(v_p)=\epsilon( v_{\phi(u)}\cdot p)=\epsilon(v_{\phi(u)})\cdot p$ for any basis element $v_p\in T^0$. Hence $\epsilon=0$ if $\epsilon|_W=0$ and $\epsilon=\id$ if $\epsilon|_W=\id_W$. Thus $V_{(T_C,\xi_C)}$ is indecomposable. 
\end{proofs}

\subsection{Representations for weighted Leavitt path algebras via branching systems}

\subsubsection{Branching systems for weighted graphs}

\begin{definition}
Let $X$ be a set, $\{R_{e_i} \mid e_i \in \hat{E}^1 \}$ and $\{ D_{v} \mid v \in E^0 \}$ families of subsets of $X$ and $\{ g_{e_i} \colon R_{e_i} \hookrightarrow D_{r(e)}
    \mid e_i \in \hat{E}^1 \}$ a family of injections such that conditions (i)-(iii) below are satisfied.
    \begin{enumerate}[(i)]
    	\item $\{ D_{v} \mid v \in E^0 \}$ is a partition of $X$ (i.e. $D_v \cap D_u = \emptyset$ whenever
    		$v \neq u$ and $\bigcup_{v \in E^0} D_v = X$).
    		\smallskip
    	\item For any $v \in E^0$ and $1 \leq i \leq w(v)$ the family
    		$\{ R_{e_i} \mid e \in s^{-1}(v), w(e) \geq i \}$ forms a partition of $D_v$.
    		    		\smallskip
    	\item Set $D_{e_i}:= g_{e_i}(R_{e_i})$ for any $e_i \in \hat{E}^1$. 
    		Then for any $e \in E^1$ the family $\{ D_{e_i} \mid 1 \leq i \leq w(e) \}$ forms a partition of $D_{r(e)}$.
    \end{enumerate}
    We call the quadruple $X=(X, \{R_{e_i} \}, \{ D_{v} \}, \{g_{e_i} \})$ an {\it $(E,w)$-branching system}.
\end{definition}

\begin{example}
Suppose that $(E,w)$ is an at most countable weighted graph, i.e. $E^0$ and $E^1$ are both 
finite or countably infinite. By fixing some linear order on $E^0$ we may write
$E^0 = \{ v^1, v^2, \dots \}$. For each
$i$ set $D_{v^i} = [i-1,i)$. Clearly, such sets are disjoint. Put $X = \bigcup_i D_{v^i}$.

Now fix a vertex $v^i \in E^0$ and $1 \le j \le w(v^i)$.
The set $S^i_j = \{ e \in s^{-1}(v^i) \mid w(e) \geq j \}$ is finite (as $E$ is row-finite). 
By ordering this set we can rewrite it as
$S^i_j = \{e^{i,j,1}, e^{i,j,2}, \dots \}$. For each $1\leq k \leq |S^i_j|$ set
\[
	R_{e^{i,j,k}_j} = [i-1+\frac{k-1}{|S^i_j|},i-1+\frac{k}{|S^i_j|}).
\]
It is clear that the set $\{ R_{e_j} \mid e \in S^i_j \}$ forms a partition of $D_{v^i}$.

In a similar fashion fix some $e \in E^1$ and let $v^i = r(e)$. For each $1 \leq j \leq w(e)$
set
\[
	D_{e_j} = [i-1+\frac{j-1}{w(e)}, i-1+\frac{j}{w(e)}).
\]
Clearly, the family of sets $\{ D_{e_j} \mid 1 \le j \le w(e) \}$ forms a partition of $D_{v^i}$.

Finally, the bijections $g_{e_j} \colon R_{e_j} \rightarrow D_{e_j}$ may be
chosen arbitrary, for example, a composition of a translation, scaling and another translation. One checks easily that $(X, \{R_{e_j} \}, \{ D_{v} \}, \{g_{e_j} \})$ is an $(E,w)$-branching system.
\end{example}

Let $X=(X, \{R_{e_i} \}, \{ D_{v} \}, \{g_{e_i} \})$ be an $(E,w)$-branching system. Let $M_X$ be the $K$-module of all functions $X \rightarrow K$ with finite support. We are going to define a structure of a right $L(E,w)$-module on $M_X$. 
In order to simplify notations, we will abuse the notation as follows. Let $Z \subseteq Y \subseteq X$ 
and $\psi : Y \rightarrow K$. By $\chi_Z \cdot \psi$ we denote the function $X \rightarrow K$
\[
	x \mapsto \begin{cases}
		\psi(x) & \text{ if } x \in Z \\
		0 & \text{ otherwise }.
	\end{cases}
\]
Using this convention, set for any $\phi \in M_X$, any $e_i \in \hat{E}^1$ and any $v \in E^0$
\begin{equation*}
\label{eq:brs:action}
\begin{aligned}
    \phi.e_i &= \chi_{D_{e_i}} \cdot (\phi \circ g^{-1}_{e_i}), \\
    \phi.e_i^* &= \chi_{R_{e_i}} \cdot (\phi \circ g_{e_i}), \\
    \phi.v &= \chi_{D_v} \cdot \phi.
\end{aligned}
\end{equation*}

By linearly extending the scalar multiplication defined above, $M_X$ becomes a right $L(E,w)$-module, see \cite[Theorem 53]{hazrat-preusser-shchegolev}. %We denote by $M^{\fin}_X$ the $L(E,w)$-submodule of $M_X$ that consists of all functions with finite support. 

\subsubsection{Branching systems versus representation graphs}
The $(E,w)$-branching sytems and the representations graphs for $(E,w)$ yield precisely the same modules for the weighted Leavitt path algebra $L(E,w)$ (up to isomorphism). This can be shown as follows. 

Let $(F, \phi)$ a representation graph for $(E,w)$. Put $X = F^0$ and $D_v = (\phi^0)^{-1}(v)$ for each $v \in E^0$. Clearly,
$\{ D_v \mid v \in E^0\}$ is a partition of $X$. Fix $v \in E^0$ and a tag $1\leq i\leq w(v)$. For each $e \in s^{-1}(v)$ such that $i \le w(e)$ set
\[R_{e_i} = \{ u \in D_v \mid \text{ there exists } f \in s_F^{-1}(u) : \phi^1(f)=e_i \}.\]
Condition (i) of Definition \ref{defwp} of a representation 
graph translates in this setting as follows:
each vertex $u$ in $D_v$ is contained in one and only one $R_{e_i}$, where
$e$ ranges over all edges of weight at least $i$ emitted by $v$. In other words,
the family $\{ R_{e_i} \mid e \in s^{-1}(v), w(e) \geq i \}$ forms a
partition for $D_v$. Now fix a vertex $v \in F^0$ and an edge $e \in r^{-1}(v)$. For each $1 \leq i \leq w(v)$
set 
\[D_{e_i} = \{ u \in D_v \mid \text{ there exists } f \in r_F^{-1}(u) \colon
\phi^1(f)=e_i\}.\]
Condition (ii) of Definition \ref{defwp} of a representation 
graph guarantees that the family of sets 
$\{ D_{e_i} \mid 1 \leq i \leq w(E) \}$ forms a partition for $D_v$. Finally, fix some $e_i \in \hat{E}^1$. We need to define a bijection
$g_{e_i} \colon R_{e_i} \rightarrow D_{e_i}$. For each $u \in R_{e_i}$
there exists precisely one edge $f \in s^{-1}(u) \cap (\phi^1)^{-1}(e_i)$.
Set $g_{e_i}(u) = r(f)$. 

\begin{theorem}\cite[Theorem 55]{hazrat-preusser-shchegolev}
	The quadruple $X = (X, \{ R_{e_i} \}, \{ D_{e_i} \}, \{ g_{e_i} \})$ defined above
	is an $(E,w)$-branching system. Moreover, the $L(E,w)$-modules $M_X$ and $V_{(F,\phi)}$ are isomorphic.
\end{theorem}

Conversely, to any $(E,w)$-branching system $X$ one can associate a representation graph $(F,\phi)$ for $(E,w)$ such that the $L(E,w)$-modules $M_X$ and $V_{(F,\phi)}$ are isomorphic, see \cite[Theorem 56 and Corollary 57]{hazrat-preusser-shchegolev}.

\section{Open problems}

Below we mention some open problems regarding weighted Leavitt path algebras.

\begin{itemize}
\item The Reduction Theorem and the Uniqueness Theorems belong to the most important results for unweighted Leavitt path algebras (see e.g. \cite[Chapter 2]{abrams-ara-molina}). Can one find analogous results for weighted Leavitt path algebras?
\smallskip
\item While the (graded) ideal structure of unweighted Leavitt path algebras is well understood (see \cite[Chapter 2]{abrams-ara-molina}), not much is known about the (graded) ideal structure of weighted Leavitt path algebras. Can one find a ``nice'' classification of the ideals and graded ideals of a weighted Leavitt path algebra?
\smallskip 
\item It is known that an unweighted Leavitt path algebra is semiprimitive, i.e. its Jacobson radical is zero. Moreover, the socle of an unweighted Leavitt path algebra has been computed (see again \cite[Chapter 2]{abrams-ara-molina}). For weighted Leavitt path algebras it is only known in some cases that the Jacobson radical is zero (cf. Section 8). Can one find a weighted Leavitt path algebras that is not semiprimitive? Can one determine the left or right socle of a weighted Leavitt path algebra?
\smallskip 
\item
As mentioned in Section 7, the zero component of an unweighted Leavitt path algebra (with respect to its standard grading) is an ultramatricial algebra, i.e. a union of an increasing chain of finite products of matrix algebras over a field. Can one find a ``good'' description of the zero component of a weighted Leavitt path algebra?
\smallskip 
\item It is known which weighted graphs $(E,w)$ have the property that $L(E,w)$ is isomorphic to an unweighted Leavitt path algebra, see Section 5. Leavitt path algebras of separated graphs and Kumjian-Pask algebras are graph algebras generalising the unweighted Leavitt path algebras, see \cite{aragoodearl} respectively \cite{pchr}. Which weighted Leavitt path algebras are isomorphic to a Leavitt path algebra of a separated graph, or to a Kumjian-Pask algebra?
\end{itemize}
$~$\\

\end{document}